\documentclass[review, 11pt, 3p, authoryear]{elsarticle}
\newdimen\Lmargin
\newdimen\Rmargin
\Rmargin=\paperwidth \advance\Rmargin by-\Lmargin \advance\Rmargin by-\textwidth
\ifdim\Lmargin>\Rmargin \Rmargin=\Lmargin \fi

\usepackage{amssymb, mathrsfs, amsthm, amsmath, bm}
\usepackage{graphicx, subfigure}
\usepackage[colorlinks=true]{hyperref}
\usepackage{multirow, booktabs, tabularx, makecell, longtable}
\usepackage{natbib}
\usepackage{lscape}
\usepackage{utfsym}
\usepackage{caption, subcaption}
\usepackage[section]{placeins}
\usepackage{lineno}

\begin{document}

\newtheorem{theorem}{Theorem}
\newtheorem{proposition}{Proposition}
\newtheorem{corollary}{Corollary}
\newtheorem{lemma}{Lemma}
\newdefinition{example}{Example}
\newdefinition{definition}{Definition}
\newdefinition{remark}{Remark}
\newdefinition{property}{Property}
\newdefinition{assumption}{Assumption}

\begin{frontmatter}
    
    \title{Preference Robust Ordinal Priority Approach with Preference Elicitation under Incomplete Information for Multi-Attribute Robust Ranking and Selection \tnoteref{t1,t2}}
    \tnotetext[t1]{Disclosure of interest: The authors declare that they have no known competing financial interests or personal relationships that could have appeared to influence the work reported in this paper.}
    \tnotetext[t2]{Data availability: Data will be made available on request.}

    \author[inst1]{Renlong Wang\corref{cor1}}
    \ead{13127073530@163.com}
    \affiliation[inst1]{organization={School of Emergency Management Science and Engineering, University of Chinese Academy of Sciences},
        city={Beijing},
        postcode={100049},
        country={China}}
    \cortext[cor1]{Corresponding author}

    \newpageafter{author}

    \begin{abstract}
        Ordinal Priority Approach (OPA) has recently been proposed to determine the weights of experts, attributes, and alternatives using ordinal preference without precise information for multi-attribute ranking and selection (MARS). This study extends OPA with preference elicitation under incomplete information to counter the parametric and preference uncertainty within MARS. Specifically, we propose Preference Robust Ordinal Priority Approach (OPA-PR) within a two-stage optimization framework to generalize marginal utility structure and resolve ambiguity in ranking parameters and utility preferences. In the first stage, the worst-case marginal utility functions are elicited from utility preference ambiguity sets, characterized by monotonicity, normalization, concavity, and Lipschitz continuity for global information, and moment-type preference elicitation for the local. In the second stage, decision weights are optimized based on the elicited marginal utility functions, considering the ranking parameters within norm-, budget-, and conditional value-at-risk-based ambiguity sets. We derive tractable reformulations of OPA-PR, especially through piecewise linear approximation for the marginal utility preference ambiguity sets for the first stage. This approximation is verified by the error bounds for both stages, establishing the foundation of preference elicitation strategy design. The proposed approach is demonstrated through a numerical experiment on the emergency supplier selection problem, including the case, sensitivity, and comparison tests.
    \end{abstract}



    \date{ }
\end{frontmatter}


\section{Introduction}
\label{section-1}

Over the past few decades, multi-attribute ranking and selection (MARS) has emerged as a critical area of operations research, particularly for tackling complex discrete decision problems characterized by conflicting objectives and diverse data \citep{C24}. A classic MARS problem involves a decision-maker (DM) determining the optimal alternative from a set of alternatives or establishing global rankings based on multiple attributes evaluated by various experts \citep{G.S.W24}. Given expert set $\mathcal{I} := \{1, \dots,I\}$, attribute set $\mathcal{J} := \{1,\dots, J\}$, and alternative set $\mathcal{K} :=\{1, \dots, K\}$, the multi-attribute group evaluation score $Z_{k}$ for alternative $k \in \mathcal{K}$ can be determined by a function $F: \mathbb{R} \to \mathbb{R}$ that has an associated collection of the weights for all expert $i \in \mathcal{I}$ and attribute $j \in \mathcal{J}$ and corresponding utility function:
\[
Z_{k} = F(v_{jk}) = \sum_{i \in \mathcal{I}} w_{i} \sum_{j \in \mathcal{J}} w_{j} u_{ij}(v_{jk}), \quad \forall k \in \mathcal{K},
\]
where $v_{jk}$ denotes the performance of alternative $k$ under attribute $j$, $u_{ij} : \mathbb{R} \to \mathbb{R}$ represents the utility function for attribute $j$ of expert $i$ that mapping the performance score $v_{jk}$ of alternative $k$ to the utility value of expert $i$ under attribute $j$, $w_{i}$ and $w_{j}$ are the weight for expert $i$ and attribute $j$, respectively.

Classical research is mainly based on the weights and utility functions obtained through sophisticated heuristic methods with complete information assumption, which means one can always elicit precise weights and utility functions from collected data \citep{A17}. If DM can provide all the necessary information to resolve MARS problems, the prior (sophisticated) methods based on precise data are advisable.
However, challenges arise when the information is incomplete in specifying the weights and utility functions \citep{G.X.Z24}. Correspondingly, there are two different sources of uncertainty that can plague a MARS problem arising from incomplete information: parameter uncertainty and preference uncertainty \citep{L.S21}. In parameter uncertainty, the DM is unaware of the exact nature of the model parameters influencing their decision, a common issue in collecting alternative performance. Preference uncertainty occurs when a decision must be made but the preferences of DMs regarding trade-offs among incomparable outcomes are not completely determined, which is typical in eliciting the utility functions for alternative performance and weights for experts and attributes. This preference uncertainty in the utility structure, referred to as utility preference ambiguity in decision science literature, is an endogenous uncertainty that arises from variability, contradiction, and diﬀiculties in accurately defining preferences with incomplete information \citep{H.Z.X.Z24}. Overall, the above uncertain situations with incomplete information have general implications in the contexts characterized by time constraints, inadequate data, and limited domain knowledge and cognitive burden of experts \citep{Z.K.Y.Z.D23}.

In this study, based on Ordinal Priority Approach (OPA), a value function-based MARS methods under ordinal preference, we propose Preference Robust Ordinal Priority Approach (OPA-PR) to counter the parametric and preference uncertainty. To effectively solve OPA-PR for practical usage, we also develop a piecewise linear approximation (PLA) scheme to derive the tractable reformulation with error bound guarantees and preference elicitation strategy. 

\subsection{Related Literature}

In this section, we aim to position our contributions in the literature by briefly reviewing relevant studies on MARS with incomplete preference information. The study considers scenarios where DM cannot provide complete information for all experts, attributes, and alternatives, whether subjective or objective, except for rankings, the minimal data required for decision-making. To provide contextual structure, we identify two primary streams of MARS with incomplete preference information since 21st century: optimization-based methods and extreme point-based methods, as well as OPA, which serves as the baseline for this study. Notably, since the 1960s, classical MARS methods and their corresponding extended models have gradually emerged to address incomplete preference information. These models have their own assumptions and axioms, achieving significant development in their respective research areas. For further information into these representative methods and their evolution, we recommend readers refer to the review by \citet{G.S.W24}, which covers the development of MARS over the past fifty years.

The optimization-based methods aim to solve the optimal weight assignment under the constraints of incomplete preference information, represented by data envelopment analysis-preference voting model (DEA-PVM) and robust ordinal regression (ROR). DEA-PVM is a typical approach based on social choice theory. It maximizes the total score of each alternative defined as the weighted sum of their votes across all ranking positions \citep{A24a}. In this process, it considers the voting performance of other alternatives, essentially serving as a relative benchmarking method. The classical model of DEA-PVM is the one proposed by \citet{C.K90}, which evaluates the position of each alternative relative to the entire set of alternatives and employs a discrimination intensity function to assess its properties, which has been shown to be equivalent to the Borda-Kendall consensus model under certain conditions. Currently, DEA-PVM has been expanded to include various types of preference information and constraints, such as ratio scales \citep{I.F19}, decreasing and convex sequences \citep{L16}, cross-evaluation \citep{S.S.L22}, and exclusion of inefficient candidates \citep{B.V16}. ROR, based on multi-attribute utility theory, identifies and utilizes a complete set of compatible instances, represented by necessary and possible weak preference relations, of value functions that reflect the preferences of all DMs involved \citep{G.M.S08}. It employs piecewise-linear marginal value functions, characterized by breakpoints when partial preference information is present, such as in pairwise and intensity comparisons. Subsequently, ROR identifies the most robust feasible weight disparities based on the partial preference information, which can be used to identify the relations between alternatives. Currently, ROR has been adapted in various ways to integrate with other methods, including stochastic multiobjective acceptability analysis (SMAA) \citep{C.G.K.S16}, non-additive value functions (represented by the Choquet integral and Sugeno integral \citep{B.W.D20}), and outranking relation preference models (such as ELECTRE \citep{C.F.G.S17} and PROMETHEE \citep{K.G.S12}).

The extreme point-based methods aim to find the extreme points that characterize the weights incorporating in a set of incomplete preference information \citep{A15}. Once identified, the final rankings of alternatives are determined by multiplying the extreme points by the attribute values of the alternatives.
\citet{A15} proposed a straightforward method for finding the extreme points of common types of incomplete preference information in the literature, including weakly ordered relations, ratio scales, absolute differences, and lower bounds on weights. \citet{A17} transformed the coefficient matrices of incomplete preference information into a class of $\mathbf{M}$-matrices to identify extreme points, subsequently minimizing the squared deviations from the extreme points to approximate the weights. Furthermore, \citet{A24} derived a dual linear programming problem to obtain closed-form solutions, identifying extreme points derived from a set of (strictly) ranked preference information. Additionally, a prevalent weight elicitation approach involves rank-based surrogate weights, where each surrogate weight can be uniquely represented by a set of extreme points \citep{B.N23}. Notable, there reveals a connection between DEA-PVM and extreme point-based methods \citep{A24a}. For instance, \citet{L24} proposed explicit expressions for weights of various simplex centroids in ranking voting systems inspired by specific simplex centroids of ROC weights.

Overall, the studies discussed above seek to engage constructively and transparently with DMs to accurately elicit and represent their evolving preferences while effectively managing imperfect preference information that may be partial, inconsistent, unstable, or uncertain. However, some studies subjectively assume a piecewise linear utility function based on multi-attribute utility theory, which is a simplified approximation to the true utility function. Theoretically, this subjective assumption fails to provide error bounds and performance guarantees for utility function approximation, going further to provide an effective preference elicitation scheme that reduces these errors.

OPA, proposed by \citet{A.M.F.L20}, is an optimization-based method for addressing MARS problems with incomplete preference information. It frames the weight elicitation problem as a linear programming model within a normalized weight space with strong dominance relations (or refers to ordinal preference), allowing for the simultaneous determination of weights for experts, attributes, and alternatives \citep{W24}. OPA utilizes ordinal data as model inputs, which are more readily available and stable compared to the cardinal values and pairwise comparisons used in other MARS methods. OPA eliminates the need for data standardization, expert opinion aggregation, and prior weight determination. Recently, several extensions of OPA have emerged, including fuzzy OPA (OPA-F) \citep{M.J.M22,P.D.G.D.K.P23}, rough set OPA (OPA-RS) \citep{D.L.J.G.C23, K.P.D.E.D23}, grey OPA (OPA-G) \citep{S.M.D.M24}, and robust OPA (OPA-R) \citep{M.A.D22} to address data uncertainty; partial OPA (OPA-P) \citep{W.S.C.S.C.G24} for managing Pareto dominance relations; TOPSIS-OPA \citep{M.D.J.Y21}, DEMATEL-OPA \citep{Z.H.L24}, and DGRA-OPA-P \citep{W24a} for large-scale group decision-making; and DEA-OPA \citep{M.A.D22a, C.W.L24} for relative efficiency analysis. However, despite its practical applications, research on the fundamental properties of OPA is limited \citep{M.J23a, M.J23}. Consequently, the lack of theoretical analysis of structural characteristics of OPA prevents current research from effectively addressing preference ambiguity and the broader utility forms it represents. Notably, \citet{W24} appears to be the only study that systematically examined the fundamental properties of the original OPA model, derived its equivalent expression, and demonstrated its closed-form solution, decomposability, and relationship with common rank-based surrogate weights, such as rank order centroid and rank reciprocal weights. Building on these properties, the Generalized Ordinal Priority Approach (GOPA) was introduced within a bilevel optimization framework, where the lower level employs cross-entropy utility minimization for preference elicitation with incomplete information, and the upper level optimizes the weights. Similarly, our discussion is driven by the insights gained from the structural characteristics and properties of OPA, forming the basis for extending OPA to handle parameter and utility preference ambiguity, which differentiates our work from that of \citet{W24}.

\subsection{Contributions}

This study proposes OPA-PR that extends OPA to the MARS scenarios involving parameter uncertainty and utility preference ambiguity under incomplete preference information. It utilizes an estimate-then-optimize two-stage procedure, where the first-stage elicit the worst-case utility functions cross all experts and attributes from the preference ambiguity set of all plausible marginal utility functions, and the second-stage optimizes the weights within the ranking parameter ambiguity sets. One of the most important components of OPA-PR are the design of preference ambiguity sets, incorporating the properties of monotonicity, normalization, concavity, and Lipschitz continuity for global information, and moment-type preference elicitation for location information, and its PLA-based tractable reformulation. The main contributions are summarized as follows: 

\begin{itemize}
    \item \textbf{Modeling:} Despite \citet{W24} proposing GOPA in bilevel optimization framework based on the cross-entropy utility minimization estimator for preference elicitation, his model actually considers the deterministic incomplete preference information, where DM only provides the preference information they certained, including the form of weak ordered relations, ratio scales, absolute differences, and lower bounds. In contrast, our approach considers undeterministic incomplete preference information represented by preference ambiguity sets, which naturally encompass deterministic cases. The preference forms discussed in \citet{W24}, along with other common forms in decision analysis literature, can be mathematically represented using the proposed preference ambiguity sets (see Examples \ref{example-01}-\ref{example-03}). Additionally, the OPA-R proposed by \citet{M.A.D22} can be viewed as a specific instance of our model when the worst-case utility employs rank-order centroid weights and attribute rankings are based on discrete scenarios with box ambiguity sets.
    \item \textbf{Methodology:} We employ a linear envelope of a piecewise linear concave marginal utility function to determine the worst-case utility function, with fixed values at finite points derived from the proposed preference elicitation strategy. This leads to a linear tractable reformulation of the first-stage problem with infinite dimensions (Proposition \ref{proposition-04}). We then derive a tractable reformulation for the second-stage problem through dual theory, incorporating ranking parameter ambiguity sets based on norm, budget, and CVaR (Propositions \ref{proposition-05}-\ref{proposition-07}), representing a novel extension for OPA-R with various ambiguity sets.
    \item \textbf{Theory:} In contrast to assuming a piecewise linear form for the marginal utility subjectively, we provide the theoretical foundation for the PLA of the marginal utility function, i.e., approximation error bounds for both stages of OPA-PR. We demonstrate that when the true utility function is concave, the step-like approximation of the partial preference information in moment-type preference elicitation is equivalent to the PLA of the marginal utility function, with no approximation error in the first stage (Propositions \ref{proposition-08} and \ref{proposition-09}), which indicates that information forms, such as deterministic utility comparisons and stochastic lottery comparisons, introduce no errors. Using a pseudo-metric, we quantify the disparities between the approximated and unapproximated arguments of the first stage (i.e., the elicited marginal utility functions) (Lemma \ref{lemma-02}), followed by error bounds for the second stage solutions (Theorem \ref{theorem-02}). Additionally, we show the decomposability of the optimal weights in OPA-PR, given the worst-case ranking parameters (Theorem \ref{theorem-01}).
\end{itemize}

\subsection{Organization}

The structure of this paper is organized as follows: Section \ref{section-2} outlines the preliminaries of OPA. Section \ref{section-3} proposes the unified framework, tractable reformulation, and approximation error bounds of OPA-PR. Section \ref{section-4} demonstrates the proposed approach through a numerical experiments on emergency supplier selection during the 7.21 mega-rainstorm disaster in Zhengzhou, China, including the case, sensitivity, comparison tests. Section \ref{section-5} gives the conclusions and future directions.

\subsection{Notation}

Throughout this paper, we use the following notation. By convention, $\mathbb{R}^{n}$ represents an $n$-dimensional Euclidean space, and $\mathbb{R}^{n\times m}$ denotes the space of $n \times m$ matrices. Vectors are denoted by bold letters, for example, $\bm{a} \in \mathbb{R}^{n}$. Matrices are denoted by boldface letters, for example, $\mathbf{A} \in \mathbb{R}^{n\times m}$. The norm form of a letter is used to represent the elements of a vector or matrix, for example, $a_{ij}$ denotes the element in the $i$-th row and $j$-th column of $\mathbf{A}$. The floral form represents a set, for example, $\mathcal{A}$. We use $k \in [K]$ to denote $k = 1, \dots, K$.

\section{Preliminaries}
\label{section-2}

Consider a classical MARS problem where DM needs to select the optimal alternative from $K$ alternatives, $\mathcal{K} :=\{1, \dots, K\}$, based on $J$ attributes, $\mathcal{J} := \{1,\dots, J\}$, as evaluated by $I$ experts, $\mathcal{I} := \{1, \dots,I\}$.
In OPA, DM initially assigns importance ranking $t_{i} \in [I]$ to each expert $i \in \mathcal{I}$.
Each expert $i \in \mathcal{I}$ is then required to provide the ranking $s_{ij} \in [J]$ for each attribute $j \in \mathcal{J}$ and the ranking $r_{ijk} \in [K]$ for each alternative $k \in \mathcal{K}$ under attribute $j \in \mathcal{J}$.
Expert evaluations are conducted independently without group discussions to ensure rankings reflect their personal preferences.
By convention, the most important attribute is ranked as 1, the next as 2, and so forth.
To streamline subsequent discussions, we first define the following three sets
\[
    \begin{aligned}
         & {{\mathcal{X}}^{1}} := \left\{(i,j,k,l)\in \mathcal{I}\times \mathcal{J}\times \mathcal{K}\times \mathcal{K} : {{r}_{ijl}}={{r}_{ijk}}+1,{{r}_{ijk}}\in [K-1] \right\}, \\
         & {{\mathcal{X}}^{2}} := \left\{(i,j,k)\in \mathcal{I}\times \mathcal{J}\times \mathcal{K}:{{r}_{ijk}}=K \right\},                                                        \\
         & \mathcal{X} := \left\{(i,j,k)\in \mathcal{I}\times \mathcal{J}\times \mathcal{K} \right\}.
    \end{aligned}
\]
Intuitively, $\mathcal{X}$ defines the indices of all experts, attributes, and alternatives, while $\mathcal{X}_{1}$ represents the set of indices of alternatives with consecutive rankings under each expert and attribute, and $\mathcal{X}_{2}$ denotes the set of indices of alternatives ranked last under each expert and attribute. The following presents the original OPA model proposed by \citet{A.M.F.L20}, referred to as OPA-\uppercase\expandafter{\romannumeral1} in this study.
\begin{equation}
    \begin{aligned}
        \text{[OPA-\uppercase\expandafter{\romannumeral1}]} \quad \max_{\bm{w}, z} \text{ } & z,  \\
        \mathrm{s.t.} \text{ }    & z \leq t_{i}s_{ij}r_{ijk}(w_{ijk}-w_{ijl}), \quad && \forall (i,j,k,l)\in \mathcal{X}^{1}, \\
        & z \leq t_{i}s_{ij}r_{ijk}(w_{ijk}), \quad && \forall (i,j,k)\in \mathcal{X}^{2}, \\
        & \sum\limits_{i=1}^{I} \sum\limits_{j=1}^{J} \sum\limits_{k=1}^{K} w_{ijk}=1,  \\
        & w_{ijk} \geq 0, \quad  && \forall (i,j,k)\in \mathcal{X},
    \end{aligned}
    \label{eq-01}
\end{equation}
where $z$ is a weight disparity scalar. After solving for the optimal solution $(z^{*}, \bm{w}^{*})$, the weights of experts, attributes, and alternatives, denoted as $W^{\mathcal{I}}$, $W^{\mathcal{J}}$, and $W^{\mathcal{K}}$, are then given by
\begin{equation}
    \begin{aligned}
            & W_{i}^{\mathcal{I}}=\sum\limits_{j=1}^{J}\sum\limits_{k=1}^{K} w_{ijk}^{*}, \quad  && \forall i\in \mathcal{I}, \\
            & W_{j}^{\mathcal{J}}=\sum\limits_{i=1}^{I}\sum\limits_{k=1}^{K}w_{ijk}^{*}, \quad   && \forall j\in \mathcal{J}, \\
            & W_{k}^{\mathcal{K}}=\sum\limits_{i=1}^{I}\sum\limits_{j=1}^{J}w_{ijk}^{*}, \quad   && \forall k\in \mathcal{K}. \\
    \end{aligned}
    \label{eq-02}
\end{equation}

OPA searches for a set of weights that maximize the weight disparity scalar for alternatives with consecutive rankings while reflecting the ordinal preferences of experts, attributes, and alternatives, within the normalized weight space. Although \citet{A.M.F.L20} refers to the outcomes as ``weights,'' they actually represent the value of the ranked object derived from the value function. However, to maintain consistency, we retain the original terminology. This difference maximization process in OPA-\uppercase\expandafter{\romannumeral1} is common in MARS methods, such as ROR, which identifies relations between alternatives regarding necessary and possible weak preference by maximizing the allowable representation error \citep{G.M.S08}. Furthermore, OPA-\uppercase\expandafter{\romannumeral1} employs a practical yet unconventional approach by multiplying the weight disparities of consecutively ranked alternatives by the rankings of experts, attributes, and alternatives. Intuitively, this leads to the intuition that as rankings increase, the marginal effect of weight disparity between consecutive alternatives decreases. \citet{W24} shows the rationale behind the manipulation in OPA-\uppercase\expandafter{\romannumeral1} by deriving an equivalent formulation and analyzing the decomposability of its solutions into the product of commonly used ranking-based surrogate weights in decision theory, which are determined by the corresponding rankings.

The following proposition presents the equivalent formulation of OPA-\uppercase\expandafter{\romannumeral1} proposed by \citet{W24}, referred to as OPA-\uppercase\expandafter{\romannumeral2} in this study.
\begin{proposition}[\citet{W24}]
Map the alternative index $k \in \mathcal{K}$ to the ranking index $r \in \mathcal{R}$ corresponding to their ranking position $r_{ijk}$ and define $\mathcal{Y} := \left\{(i,j,r) \in \mathcal{I} \times \mathcal{J} \times \mathcal{R}\right\}$. OPA-\uppercase\expandafter{\romannumeral1} has the following equivalent formulation
\begin{equation}
    \begin{aligned}
        \text{[OPA-\uppercase\expandafter{\romannumeral2}]} \quad  \max_{\bar{\bm{w}}, \bar{z}} \text{ } & \bar{z}, \\
        \mathrm{s.t.} \text{ } & R U_{r}^{ROC} \bar{z} \leq t_{i}s_{ij}\bar{w}_{ijr}, \quad &  & \forall (i,j,r)\in \mathcal{Y}, \\
        & \sum\limits_{i=1}^{I}\sum\limits_{j=1}^{J}\sum\limits_{r=1}^{R} \bar{w}_{ijr} = 1,  \\
        & \bar{w}_{ijr} \geq 0, \quad && \forall (i,j,r) \in \mathcal{Y},
    \end{aligned}
    \label{eq-03}
\end{equation}
where $U_{r}^{ROC} = (\sum_{h = r}^{R} \frac{1}{h}) / R$ is the rank order centroid weight for the alternative ranked $r \in \mathcal{R}$. After solving for the optimal solution $(z^{*}, \bar{\bm{w}}^{*})$, $\bar{w}_{ijr}^{*}$ are mapped to $w_{ijk}^{*}$ to calculate the weights of experts, attributes, and alternatives based on Equation \eqref{eq-02}.
\label{proposition-02}
\end{proposition}

OPA-\uppercase\expandafter{\romannumeral2} further reveals the logic behind weight assignment of OPA. It searches for a set of weights that maximize the weight disparity scalar with the rank order centroid weights for ranked alternatives while reflecting the ordinal preferences of experts and attributes, within the normalized weight space.

The following corollary illustrates that the optimal weights of OPA-\uppercase\expandafter{\romannumeral1} can be decomposed into rank-based surrogate weights independently determined by the ranking indices of experts, attributes, and alternatives, thereby forming a diminishing marginal effect commonly observed in rank-based surrogate weights \citep{B.N23}.
\begin{corollary}[\citet{W24}]
    Let $w_{l}^{ROC}$ and $w_{l}^{RR}$ denote the rank order centroid weight and rank reciprocal weight of the object ranked $l$, defined as $w_{l}^{ROC} = (\sum_{h=l}^{L} \frac{1}{h}) / L$ and $w_{l}^{RR} = {1}/(l \sum_{h = 1}^{L} \frac{1}{h})$ for all $l =1, \dots, L$. The optimal weights for OPA-\uppercase\expandafter{\romannumeral1} are equivalent to
    \[
        w_{ijk}^{*} = w_{t_{i}}^{RR}w_{s_{ij}}^{RR}w_{r_{ijk}}^{ROC}, \quad \forall (i,j,k) \in \mathcal{X}.
    \]
    \label{corollary-01}
\end{corollary}

From Corollary \ref{corollary-01}, it can be inferred that the identical constraint coefficients $t_{i}s_{ij}r_{ijr}$ in OPA-\uppercase\expandafter{\romannumeral1} do not necessarily indicate identical weights. Identical weights only occur when both $r_{ijk}$ and $t_{i}s_{ij}$ are correspondingly the same. Corollary \ref{corollary-01} further demonstrates that solving the multi-expert OPA model is equivalent to solving the single-expert OPA model and then applying rank reciprocal weights for aggregation, which is adopted in the rest of this study. For all $i \in \mathcal{I}$, define
\[
\begin{aligned}
    & \bm{s}_{i} = (\underbrace{s_{i1}, \dots, s_{i1}}_{R \text{ elements}} , \dots, \underbrace{s_{iJ}, \dots, s_{iJ}}_{R \text{ elements}})^{\top}, \\
    & \bm{U}^{ROC} = \big(\underbrace{U_{1}^{ROC}, \dots, U_{R}^{ROC}, \dots, U_{1}^{ROC}, \dots, U_{R}^{ROC}}_{JR \text{ elements}} \big)^{\top}, \\
    & \bar{\bm{w}}_{i} = (\underbrace{\bar{w}_{i11}, \dots, \bar{w}_{i1R}, \bar{w}_{i21}, \dots, \bar{w}_{iJR}}_{JR \text{ elements}})^{\top}, \\
    & \bar{\mathbf{W}}_{i} = \mathrm{diag}(\underbrace{\bar{w}_{i11}, \dots, \bar{w}_{i1R}, \bar{w}_{i21}, \dots, \bar{w}_{iJR}}_{JR \text{ elements}})^{\top}.
\end{aligned}
\]

We can rewrite the OPA-\uppercase\expandafter{\romannumeral2} for any $i \in \mathcal{I}$ (i.e., the single-expert model) in matrix form as
\begin{equation}
    \begin{aligned}
        \max_{\bar{\bm{w}}_{i} \in \bar{\mathcal{W}}_{i}, \bar{z}_{i}} \text{ } & \bar{z}_{i},  \\
        \mathrm{s.t.} \text{ } & R \bm{U}^{ROC} \bar{z}_{i} \leq \bar{\mathbf{W}}_{i}\bm{s}_{i}, \\
        & \bar{\mathcal{W}}_{i} := \left\{ \bar{\bm{w}}_{i} \in \mathbb{R}^{JR} \left| \begin{array}{l} 
            \bm{e}^{\top}\bar{\bm{w}}_{i} = 1, \\
            \bar{\bm{w}}_{i} \geq 0,
        \end{array}\right. \right\},
    \end{aligned}
    \label{eq-04}
\end{equation}
where $\bm{e}$ denotes a $JR$-dimensional vector of all ones.

\section{Preference Robust Ordinal Priority Approach}
\label{section-3}

This section introduces the preference robust ordinal priority approach (OPA-PR) based on a bilevel optimization framework to counter the ambiguity within both the utility structure and ranking parameter.

\subsection{Modeling Framework}
\label{section-3-1}

The components that may introduce uncertainty into OPA in Equation \eqref{eq-04} include the alternative utilities and attribute ranking parameters. For alternative utilities, it is desired to extend the rank order centroid weights for alternatives to a more general utility structure that accounts for ambiguity in the preference information. This can significantly expand the flexibility of the proposed approach by customizing it to reflect different preferences, rather than restricting it to a specific alternative utility structure as in OPA. The focus of OPA-PR in preference ambiguity is the situation where expert $i \in \mathcal{I}$ does not have complete information to uniquely specify their marginal utility function $u_{ij}: \mathbb{R} \to \mathbb{R}$ for attribute $j \in \mathcal{J}$ that maps alternative rankings to utility values. However, partial information can be gathered to construct the ambiguity sets of marginal utility functions, denoted as $\mathcal{U}_{ij}$, such that the true marginal utility function reflecting the expert preference falls within the ambiguity set with high likelihood. Similarly, for ranking parameter ambiguity, a common robust optimization approach is applied to construct the ambiguity sets for attribute rankings $\tilde{\bm{s}}_{i}$, denoted as $\mathcal{V}_{i}$. 

This study formulates OPA-PR in an estimate-then-optimize two-stage procedure. Consider a measurable space $(\Omega, \mathcal{F})$, where $\Omega$ is a finite non-empty set and $\mathcal{F}$ is a $\sigma$-algebra on $\Omega$. The finite sample assumption that $\Omega := \{ \bm{\xi}_{1}, \dots, \bm{\xi}_{E} \}$ with $1 \leq E < \infty$ and $\mathcal{E}:= \{1, \dots, E\}$ is common in literature related to risk preference, which can be regarded as a discrete sample approximation of the continuous sample space \citep{W.H.H.X22}. The first stage involves eliciting the worst-case marginal utility function $u^{*}_{ij}$ for all $(i,j) \in \mathcal{I} \times \mathcal{J}$ from a set of plausible marginal utility functions $\mathcal{U}_{ij}$ over the expectation of random return $h(x, \bm{\xi})$ of lottery $\bm{x} \in \mathbb{R}^{m}$ indexed by the finite scenario $\bm{\xi} \in \mathbb{R}^{n}$ with associated probabilities $p_{e} = \mathbb{P}[\bm{\xi} = \bm{\xi}_{e}]$ for $e \in \mathcal{E}$. In the second stage, informed by the estimated worst-case marginal utility function. It optimize for the optimal decision weights under the worst-case ranking parameters from the corresponding ambiguity set $\mathcal{V}_{i}$. Let $U_{ijr}^{*} = \int_{0}^{R-r+1} \mathrm{d}u_{ij}^{*}(x)$ for all $(i,j,r) \in \mathcal{Y}$ and 
\[
    \bm{U}_{i}^{*} = \big(\underbrace{U_{i11}^{*}, \dots, U_{i1R}^{*}, U_{i21}^{*}, \dots, U_{iJR}^{*}}_{JR \text{ elements}} \big)^{\top}.
\]

For any $i \in \mathcal{I}$, we consider the modeling framework in Equation \eqref{eq-05} for OPA-PR.
\begin{equation}
    \begin{aligned}
        & \begin{aligned}
            \max_{\bar{\bm{w}}_{i} \in \bar{\mathcal{W}}_{i}, z_{i}} \text{ } &  \bar{z}_{i} \\
            \mathrm{s.t.} \text{ } &  R \bm{U}_{i}^{*} \bar{z}_{i} \leq \bar{\mathbf{W}}_{i} \tilde{\bm{s}}_{i} \quad\quad\quad\quad\quad\quad\quad\quad\text{ } \forall \tilde{\bm{s}}_{i} \in \mathcal{V}_{i}
        \end{aligned} && \quad \quad \text{(Second Stage: Optimize)} \\
        & \quad\quad\quad\quad u^{*}_{ij} = \underset{u_{ij} \in \mathcal{U}_{ij}}{\arg\min} \text{ } \mathbb{E}_{\mathbb{P}}[u_{ij}(h(\bm{x}, \bm{\xi}))] \quad\quad \forall j \in \mathcal{J} && \quad \quad \text{(First Stage: Estimate)}
    \end{aligned}
    \label{eq-05}
\end{equation}

The proposed OPA-PR framework clarifies the optimal weights when experts exhibit conservative (risk-averse) behavior against uncertainty. It optimizes the worst-case marginal utility function from the ambiguity set that reflects the preference attitude of expert in the first stage, providing a more interpretable formulation. The key to success is designing the ambiguity sets for both marginal utility functions and attribute rankings, particularly the ambiguity set for marginal utility functions based on partial preference information of experts, and further determining the tractable reformulation of Equation (\ref{eq-05}).

\subsection{First-Stage Estimation: Worst-Case Marginal Utility Functions}
\label{section-3-2}

This section presents the ambiguity set design of the marginal utility functions and develops a tractable reformulation for the first-stage estimation of OPA-PR.

\subsubsection{Ambiguity Set Design for Marginal Utility Functions}
\label{section-3-2-1}

We begin with designing the ambiguity set for the marginal utility functions related to the first-stage problem of OPA-PR, which is formed by the intersection of the following properties of marginal utility function: monotonicity, normalization, concavity, Lipschitz continuity, and moment-type preference elicitation.

Let $\mathscr{U}$ represent a class of real-valued utility functions, where any $u \in \mathscr{U}$ is piecewise continuously differentiable with a finite number of non-differentiable turning points at some of the rankings:
\[
    0 = \tau_{0} < \cdots < \tau_{H_{ij}} = R,
\]
and $\mathcal{H}_{ij} := \{1, \dots, H_{ij}\}$ for all $(i,j) \in \mathcal{I} \times \mathcal{J}$.
Suppose that the input $x$ is bounded on $\Theta := (0, \theta]$ and $\theta = R$.

\begin{assumption}[Monotonicity] 
    The set of utility functions satisfying the monotonicity is denoted as $\mathcal{U}^{\text{mon}}$, where $u \in \mathscr{U}$ is monotonic if $x \preceq y$ implies $u(x) \leq u(y)$ for all $x,y \in \Theta$.
    \label{assumption-01}
\end{assumption}
\begin{assumption}[Normalization]
    The set of utility functions satisfying the normalization is denoted as $\mathcal{U}^{\text{nor}}$, where
    $u \in \mathscr{U}$ is normalized if $u(0) = 0$ and $u(\theta) = 1$.
    \label{assumption-02}
\end{assumption}

The monotonicity and normalization of utility functions are common in literature related to decision theory \citep{W.Z.W.X23}.

\begin{assumption}[Concavity]
    The set of utility functions satisfying the concavity is denoted as $\mathcal{U}^{\text{conc}}$, where
    $u \in \mathscr{U}$ is concave if $u(\lambda x + (1-\lambda) y) \geq \lambda u(x) + (1-\lambda)u(y)$ for all $x,y \in \Theta$ and $\lambda \in [0,1]$.
    \label{assumption-03}
\end{assumption}

The concavity of the utility function reflects the risk aversion preference, implying that for any random lottery $X$, one prefers the certain outcome $\mathbb{E}[X]$ over the lottery $X$ itself \citep{A.D15}.
Although this concavity assumption may not apply to all scenarios, it is generally suitable for most decision-making situations \citep{T13}. 

\begin{assumption}[Lipschitz continuity]
    The set of utility functions satisfying the Lipschitz continuity is denoted as $\mathcal{U}^{\text{lip}}$, where $u \in \mathscr{U}$ is Lipschitz continuous with modulus being bounded by $G$ if $|u(x) - u(y)| \leq G \left\| x - y \right\|$ for all $x,y \in \Theta$.
    \label{assumption-04}
\end{assumption}

Lipschitz continuity can be interpreted as the restriction that a finite input cannot lead to an infinite improvement, meaning the preferences cannot change too rapidly around any specific input \citep{G.X.Z24}.
It establishes an upper bound $G$ for the the first-order derivative of utility function over $\Theta$.
Thus, when the expert specifies their nominal utility function, we can derive $G$ from its first-order derivative, such as the maximum of $\nabla u(\tau) = 1/(\theta)$ for risk-neutral utility function and $\nabla u(\tau) = \left(\gamma e^{-\gamma \tau} \right) / \left( 1 - e^{-\gamma} \right)$ for constant absolute risk-aversion utility function over $\Theta$.
Additionally, Lipschitz continuity aids in establishing error bounds for ambiguity set approximations, which will be shown in Section \ref{section-3-4}.

\begin{assumption}[Moment-type preference elicitation]
    The set of utility functions satisfying the moment-type preference elicitation is defined as
    \begin{equation}
        \mathcal{U}^{\text{mpre}} := \left\{ u \in \mathscr{U}: -\infty < \int_{0}^{\theta} \psi_{l}(\tau)\mathrm{d}u(\tau)\leq c_l \text{ for } l \in \mathcal{L} \right\},
        \label{eq-06}
    \end{equation}
     where $\psi_l: \Theta \to \mathbb{R}$ are Lebesgue integrable functions and $c_{l}$ are some given constants for $l \in \mathcal{L} := \{1, \dots, L\}$.
    \label{assumption-05}
\end{assumption}

The term $\psi_{l}$ for $l \in \mathcal{L}$ reflects partial preference information, which varies across attributes among different experts. It is easy to verify that $\mathcal{U}^{\text{mpre}}$ is a convex set. The moment-type preference elicitation set encompasses various forms of partial preference information discussed in the literature related to decision theory, with several sensible examples presented below.

\begin{example}[Value comparison \citep{W24}]
    The value comparison primarily involves partial preference information for a value prospect, including ratio scales, absolute differences, and lower bounds, with the corresponding set of utility function defined as:
    \[
        \begin{aligned}
            \mathcal{U}^{\text{rs}} & := \left\{ u \in \mathscr{U} : \int_{0}^{r} \mathrm{d}u(\tau) - \alpha_{l_{1}} \int_{0}^{r-1} \mathrm{d}u(\tau) = 0 \text{ for } l_{1} \in \mathcal{L}_{1} \right\}, \\
            \mathcal{U}^{\text{ad}} & := \left\{ u \in \mathscr{U} : \int_{0}^{r} \mathrm{d}u(\tau) - \int_{0}^{r-1} \mathrm{d}u(\tau) = \beta_{l_{2}} \text{ for } l_{2} \in \mathcal{L}_{2} \right\},    \\
            \mathcal{U}^{\text{lb}} & := \left\{ u \in \mathscr{U} : \int_{0}^{r} \mathrm{d}u(\tau) = \gamma_{l_{3}} \text{ for } l_{3} \in \mathcal{L}_{3} \right\},
        \end{aligned}
    \]
    where $\mathcal{U}^{\text{rs}}$, $\mathcal{U}^{\text{ad}}$, and $\mathcal{U}^{\text{lb}}$ represent the set of utility functions with ratio scales, absolute differences, and lower bounds, respectively.
    In cases of complete preference information, $\mathcal{U}^{\text{rs}}$ is commonly applied in pairwise-comparison-based MARS methods like AHP, ANP, and BWM, while $\mathcal{U}^{\text{ad}}$ is prevalent in quasi-distance-based methods like TOPSIS, VIKOR, and EDAS.
    We can express the above sets using indicator functions as follows:
    \[
        \begin{aligned}
            \mathcal{U}^{\text{rs}} & := \left\{ u \in \mathscr{U} : \int_{0}^{\theta} (\mathbf{1}_{(0,r]}(\tau) - \alpha_{l_{1}} \mathbf{1}_{(0,r-1]}(x)) \mathrm{d}u(\tau) = 0 \text{ for } l_{1} \in \mathcal{L}_{1} \right\}, \\
            \mathcal{U}^{\text{ad}} & := \left\{ u \in \mathscr{U} : \int_{0}^{\theta} (\mathbf{1}_{(0,r]}(\tau) - \mathbf{1}_{(0,r-1]}(\tau)) \mathrm{d}u(\tau) = \beta_{l_{2}} \text{ for } l_{2} \in \mathcal{L}_{2} \right\},    \\
            \mathcal{U}^{\text{lb}} & := \left\{ u \in \mathscr{U} : \int_{0}^{\theta} \mathbf{1}_{(0,r]}(\tau) \mathrm{d}u(\tau) = \gamma_{l_{3}} \text{ for } l_{3} \in \mathcal{L}_{3} \right\},
        \end{aligned}
    \]
    where $I_{r}(\tau)$ is indicator functions with domain $[0,r]$ with $\tau \in \Theta$. Without loss of generality, we can unify the above sets into the ambiguity set with value comparison:
    \[
        \mathcal{U}^{\text{duc}} := \left\{ u \in \mathscr{U} : \int_{0}^{\theta} \eta_{l}(\tau) \mathrm{d}u(\tau) = c_{l} \text{ for } l \in \mathcal{L}  \right\},
    \]
    which is consistent with the form of $\mathcal{U}^{\text{mpre}}$. Notably, $\eta_{l}$ is a step function with jumps at given constants $c_{l}$ for $l \in \mathcal{L}$.
    \label{example-01}
\end{example}
\begin{example}[Stochastic lottery comparison \citep{A.D15}]
    Given any lottery with cumulative distribution $F(\tau)$ on the domain $\Theta$ with the expected utility $u(r_{1}) \leq \mathbb{E}[F(\tau)] \leq u(r_{2})$.
    By integration by parts, we have
    \[
        \int_{0}^{\theta} u(\tau) \mathrm{d} F(\tau) = u(\tau) F(\tau) \Big|_{0}^{\theta} - \int_{0}^{\theta} F(\tau) \mathrm{d}u(\tau) \ge \int_{0}^{r_{1}} \mathrm{d} u(\tau) \Rightarrow \int_{0}^{\theta} (F(\tau) + \mathbf{1}_{(0,r_{1}]}(\tau)) \mathrm{d} u(\tau) \leq 1,
    \]
    and
    \[
        \int_{0}^{\theta} u(\tau) \mathrm{d} F(\tau) = u(\tau) F(\tau) \Big|_{0}^{\theta} - \int_{0}^{\theta} F(\tau) \mathrm{d}u(\tau) \leq \int_{0}^{r_{2}} \mathrm{d} u(\tau) \Rightarrow \int_{0}^{\theta} (F(\tau) + \mathbf{1}_{(0,r_{2}]}(\tau)) \mathrm{d} u(\tau) \ge 1,
    \]
    where $I_{r}(\tau)$ is the indicator function and $u(\tau) F(\tau) \Big|_{0}^{\theta} = u(\tau) G(\tau) \Big|_{0}^{\theta} = 1$ is from the normalization condition.
    Then, we have the ambiguity set with stochastic lottery comparisons:
    \[
        \begin{aligned}
            \mathcal{U}^{\text{slc}} := \left\{ u \in \mathscr{U}: \int_{0}^{\theta} (F(\tau) + \mathbf{1}_{(0,r_{1}]}(\tau)) \mathrm{d} u(\tau) \leq 1, \int_{0}^{\theta} (- F(\tau) - \mathbf{1}_{(0,r_{2}]}(\tau)) \mathrm{d} u(\tau) \leq 1 \right\},
        \end{aligned}
    \]
    which is consistent with the form of $\mathcal{U}^{\text{mpre}}$.
    If the lottery involves pairwise comparisons with corresponding singleton probability, then $\psi_{l}$ for $l \in \mathcal{L}$ in $\mathcal{U}^{\text{slc}}$ is a step function.
    \label{example-02}
\end{example}
\begin{example}[Stochastic dominance relation \citep{H.M15}]
    Given any two lotteries with cumulative distribution $F(\tau)$ and $G(\tau)$ with $\tau \in \Theta$. Suppose that DM prefers $F(\tau)$ to $G(\tau)$.
    By the expected utility theory, we have $\mathbb{E}[F(\tau)] \ge \mathbb{E}[G(\tau)]$, which can be rewriten by integration by parts as:
    \[
        \int_{0}^{\theta} u(\tau) \mathrm{d} F(\tau) \ge \int_{0}^{\theta} u(\tau) \mathrm{d} G(\tau) \Rightarrow \int_{0}^{\theta} [G(\tau) - F(\tau)] \mathrm{d} u(\tau) \ge 0.
    \]
    Then, we have the ambiguity set with stochastic dominance relation:
    \[
        \mathcal{U}^{\text{sdr}} := \left\{ u \in \mathscr{U} : \int_{0}^{\theta}  [F(\tau) - G(\tau)] \mathrm{d} u(\tau) \leq 0 \right\},
    \]
    which is consistent with the form of $\mathcal{U}^{\text{mpre}}$.
    \label{example-03}
\end{example}

The ambiguity set for marginal utility function of OPA-PR can be derived from the intersection of the previously discussed sets
\begin{equation}
    \mathcal{U} := \mathcal{U}^{\text{mon}} \cap \mathcal{U}^{\text{nor}} \cap \mathcal{U}^{\text{conc}} \cap \mathcal{U}^{\text{lip}} \cap \mathcal{U}^{\text{mpre}},
    \label{eq-07}
\end{equation}
where monotonicity, normalization, concavity, and Lipschitz continuity represent global information, while moment-type preference elicitation serves as local customized information of marginal utility function.
Notably, the ambiguity sets of marginal utility functions differ from each attribute given by experts, forming the basis for eliciting customized preferences.
Thus, for any $(i,j) \in \mathcal{I} \times \mathcal{J}$, we have the first-stage problem of OPA-PR
\begin{equation}
    \rho_{ij} = \underset{u_{ij} \in \mathcal{U}_{ij}}{\min} \sum_{e \in \mathcal{E}} p_{e} u_{ij}(h(\bm{x}, \bm{\xi}_{e})).
    \label{eq-08}
\end{equation}

\subsubsection{First-Stage Tractable Reformulation}
\label{section-3-2-2}

This section introduces the piecewise linear approximation (PLA) of the ambiguity set of marginal utility functions to derive a tractable reformulation of the first-stage elicitation of the worst-case marginal utility function.

We begin by presenting the PLA of the ambiguity set of marginal utility functions, specifically $\mathcal{U}^{\text{mpre}}$ in $\mathcal{U}$, which is nonlinear within the space of utility function $\mathscr{U}$.
PLA is of interest for two main reasons.
First, the ordinal preferences are elicited at discrete rankings in OPA.
Connecting utility values at these points to form a piecewise linear marginal utility function is a straightforward way to obtain an approximated marginal utility function  \citep{W.Z.W.X23}.
Second, the piecewise linear assumption is widely used in decision theory research, as seen in various MARS methods such as UTA, MACBETH, and ROR \citep{G.S.W24}.

\begin{definition}
    Let $\tilde{\mathscr{U}} \subset \mathcal{U}$ where any $\tilde{u} \in \tilde{\mathscr{U}}$ is a piecewise linear function with non-differentiable turning points on $\Theta$.
    The set of piecewise linear utility functions satisfying moment-type preference elicitation is defined as
    \begin{equation}
        \tilde{\mathcal{U}}^{\text{mpre}} := \left\{ \tilde{u} \in \tilde{\mathscr{U}}: -\infty < \int_{0}^{\theta} \psi_{l}(\tau)\mathrm{d}\tilde{u}(\tau)\leq c_l \text{ for } l \in \mathcal{L}\right\}.
        \label{eq-09}
    \end{equation}
    \label{definition-01}
\end{definition}

For any $u_{ij} \in \mathcal{U}_{ij}^{\text{mpre}}$, the corresponding piecewise linear marginal utility function $\tilde{u}_{ij} \in \tilde{\mathcal{U}}_{ij}^{\text{mpre}}$ can be constructed by connecting the function values at the endpoints of each interval $(\tau_{h_{ij}-1},\tau_{h_{ij}}]$ for any $h_{ij} \in \mathcal{H}_{ij}$, which is given by the following proposition.
\begin{proposition}
    For any $(i,j) \in \mathcal{I} \times \mathcal{J}$, suppose that $\psi_{l_{ij}}(\tau)$ for any $l_{ij} \in \mathcal{L}_{ij}$ is a step function on $\Theta$ with jumps at $\tau_{h_{ij}}$ for all $h_{ij} \in \mathcal{H}_{ij}$.
    Then, for any $u_{ij} \in \mathcal{U}^{\text{mpre}}_{ij}$, there exists a piecewise linear marginal utility function $\tilde{u}_{ij} \in \tilde{\mathscr{U}}_{ij}$ with $\tilde{u}(\tau_{h_{ij}}) = u(\tau_{h_{ij}})$ for all $h_{ij} \in \mathcal{H}_{ij}$, such that $\tilde{u}_{ij} \in \tilde{\mathcal{U}}_{ij}^{\text{mpre}}$.
    Specifically, such $\tilde{u}$ can be constructed by
    \begin{equation}
        \tilde{u}_{ij}(t) = \sum_{h_{ij} \in \mathcal{H}_{ij}} \left(u_{ij}(\tau_{h_{ij}-1}) + \frac{u_{ij}(\tau_{h_{ij}}) - u_{ij}(\tau_{h_{ij}-1})}{\tau_{h_{ij}} - \tau_{h_{ij}-1}}(t - \tau_{h_{ij}-1})\right) \mathbf{1}\left\{t \in (\tau_{h_{ij}-1},\tau_{h_{ij}}]\right\}, \forall t \in \Theta.
        \label{eq-10}
    \end{equation}
    \label{proposition-03}
\end{proposition}

By Proposition \ref{proposition-03}, for any $(i,j) \in \mathcal{I} \times \mathcal{J}$, we can introduce the approximated first-stage problem of OPA-PR 
\begin{equation}
    \tilde{\rho}_{ij} = \underset{\tilde{u}_{ij} \in \tilde{\mathcal{U}}_{ij}}{\min} \sum_{e \in \mathcal{E}} p_{e} \tilde{u}_{ij}(h(\bm{x}, \bm{\xi}_{e})).
    \label{eq-11}
\end{equation}
where $\tilde{\mathcal{U}}_{ij} := \mathcal{U}^{\text{mon}} \cap \mathcal{U}^{\text{nor}} \cap \mathcal{U}^{\text{conc}} \cap \mathcal{U}^{\text{lip}} \cap \tilde{\mathcal{U}}^{\text{mpre}}_{ij}$.
In Section \ref{section-3-4}, we will analyze the error bounds on the optimal value and solution of both stages of OPA-PR introduced by PLA.

We then derive a tractable reformulation of the approximated first-stage problem in Equation \eqref{eq-11}, which has infinite dimensions. The following lemma states that, given a set of points in $\mathbb{R}^{n} \times \mathbb{R}$, a concave function mapping from $\mathbb{R}^{n}$ to $\mathbb{R}$ can be expressed as the upper envelope of linear functions, forming a piecewise linear concave function that passes through some of those points.
\begin{lemma}[\citet{H.H.X18}]
    Let $f: \mathbb{R}^{n} \to \mathbb{R}$ and $\left< \cdot \right>$ denote the Euclidean inner product. The following statement holds.

    (1) $f$ is concave if and only if 
    \[
    f(x) = \inf_{i \in \mathcal{I}} g_{i}(\bm{x}), \quad \forall \bm{x} \in \mathrm{dom} f,
    \]
    where $I$ is an index set, potentially infinite, and $g_{i}(\bm{x}) = \left< \bm{a}_{i}, \bm{x} \right> + \bm{b}_{i}$ for all $i \in \mathcal{I}$, representing the support function of $f$ at $\bm{x}$ for any $a \in \partial f(\bm{x})$.

    (2) For any finite set $\mathcal{O} \subset \mathbb{R}^{n}$ and values $\{ \bm{v}_{o} \}_{o \in \mathcal{O}} \subset \mathbb{R}$, the function $\hat{f}: \mathbb{R}^{n} \to \mathbb{R}$ defined by
    \[
        \hat{f}(x) = \min_{\bm{a},\bm{b}}\{ \left< \bm{a}, \bm{x} \right> + \bm{b}: \left< \bm{a}, \bm{o} \right> + \bm{b} \geq \bm{v}_{o}, \forall \bm{o} \in \mathcal{O} \}
    \]
    is concave. Furthermore, $\hat{f} \leq \tilde{f}$ holds over $\mathbb{R}^{n}$ for all increasing and concave functions $\tilde{f}$ with $\tilde{f}(\bm{o}) \geq \bm{v}_{o}$.
    \label{lemma-01}
\end{lemma}

The following proposition gives the tractable reformulation of Equation \eqref{eq-11} based on Lemma \ref{lemma-01}.
Our result shows that Equation \eqref{eq-11} can be solved by a finite-dimensional linear programming, involving $2(R + E)$ variables and $4R + ER + L + E$ constraints without nonnegativity constraints.
\begin{proposition}
    For any $(i,j) \in \mathcal{I} \times \mathcal{J}$, Equation \eqref{eq-11}, satisfying Assumptions \ref{assumption-01}-\ref{assumption-05}, can be reformulated as the following minimization problem in Equation \eqref{eq-12}, which is a finite linear programming problem given $\bm{x}$ and $\bm{\xi}$, where $h(\bm{x}, \bm{\xi})$ is affine in $\bm{x}$.
    \begin{subequations}\label{eq-12}
        \begin{align}
            \min_{\bm{a}, \bm{b}, \bm{y}, \bm{\mu}} \text{ } & \sum_{e \in \mathcal{E}} p_{e} (a_{e} h(\bm{x}, \bm{\xi}_{e}) + b_{e}) \label{eq-12-01} \\
            \mathrm{s.t.} \text{ } & y_{h_{ij}} - y_{h_{ij}-1} = \mu_{h_{ij}-1}(\tau_{h_{ij}} - \tau_{h_{ij}-1}) \quad && \forall h_{ij} \in \mathcal{H}_{ij}          \label{eq-12-02} \\
            & y_{h_{ij}} - y_{h_{ij}-1} \geq \mu_{h_{ij}}(\tau_{h_{ij}} - \tau_{h_{ij}-1}) \quad && \forall h_{ij} \in \mathcal{H}_{ij}\backslash\{H_{ij}\} \label{eq-12-03} \\
            & 0 \leq \mu_{h_{ij}} \leq G \quad && \forall h_{ij} \in \mathcal{H}_{ij} \cup \{0\} \backslash\{H_{ij}\} \label{eq-12-04} \\
            & \sum_{h_{ij} \in \mathcal{H}_{ij}} \mu_{h_{ij}-1} \int_{\tau_{h_{ij}-1}}^{\tau_{h_{ij}}} \psi_{l_{ij}}(\tau) \mathrm{d}t \leq c_{l} \quad && \forall l_{ij} \in \mathcal{L}_{ij}  \label{eq-12-05} \\
            &  a_{e} \tau_{h_{ij}} + b_{e} \ge y_{h_{ij}} \quad  && \forall e \in \mathcal{E}, h_{ij} \in \mathcal{H}_{ij}  \label{eq-12-06} \\
            & y_{0} = 0, y_{R} = 1 \label{eq-12-07} \\
            & a_{e} \ge 0 \quad && \forall e \in \mathcal{E} \label{eq-12-08} 
        \end{align}
    \end{subequations}
    Given the optimal solution $(\bm{a}^{*}, \bm{b}^{*}, \bm{y}^{*}, \bm{\mu}^{*})$, the worst-case marginal utility function is given by
    \begin{equation}
        \tilde{u}_{ij}^{*}(t) =  \sum_{h_{ij} \in \mathcal{H}_{ij}} \left\{ \frac{y_{h_{ij}}^{*} - y_{h_{ij}-1}^{*}}{\tau_{h_{ij}} - \tau_{h_{ij}-1}} t + \frac{\tau_{h_{ij}} y_{h_{ij}-1}^{*} - \tau_{h_{ij}-1} y_{h_{ij}}^{*}}{\tau_{h_{ij}} - \tau_{h_{ij}-1}} \right\} \mathbf{1}\left\{t \in (\tau_{h_{ij}-1},\tau_{h_{ij}}]\right\}, \quad \forall t \in \Theta,
        \label{eq-13}
    \end{equation}
    with $\tilde{u}_{ij}^{*}(0) = 0$ and $\tilde{u}_{ij}^{*}(R) = 1$.
    \label{proposition-04}
\end{proposition}

By Proposition \ref{proposition-04}, we can derive the worst-case marginal utility functions in piecewise linear form for all experts and attributes, which are then incorporated into the second-stage problem of OPA-PR to determine the optimal decision weights.

\subsection{Second-Stage Optimization: Decision Weights}
\label{section-3-3}

This section introduces several ambiguity sets for attribute ranking parameters and derives the corresponding tractable reformulation of the second-stage optimization of OPA-PR.

Let $\tilde{\bm{U}}_{i}^{*} \in \mathbb{R}^{J \times R}$, where $\tilde{U}_{ijr} = \int_{0}^{R-r+1} \mathrm{d}\tilde{u}_{ij}^{*}(x)$. The second-stage optimization deals with the following problem
\begin{equation}
    \begin{aligned}
        \max_{\tilde{z}_{i}, \tilde{\bm{w}}_{i} \in \bar{\mathcal{W}}_{i}} \text{ } & \tilde{z}_{i}, \\
        \mathrm{s.t.} \text{ } & R \tilde{\bm{U}}_{i}^{*} \tilde{z}_{i} \leq \tilde{\mathbf{W}}_{i} \tilde{\bm{s}}_{i}, \quad \forall \tilde{\bm{s}}_{i} \in \mathcal{V}_{i}. \\
    \end{aligned}
    \label{eq-14}
\end{equation}

We primarily consider the following ambiguity sets for the attribute ranking parameters of the second-stage problem in OPA: norm-, budget-, and conditional value-at-risk (CVaR)-based ambiguity sets, which are the mainstream ambiguity sets considered in decision analysis.

\begin{assumption}[Norm-based ambiguity set] The norm-based ambiguity set for the attribute ranking $\tilde{\bm{s}}$ is defined by
\[
    \mathcal{V}^{\text{norm}} := \left\{ \tilde{\bm{s}} \in \mathbb{R}^{JR}:\|\bm{\Sigma}^{-\frac{1}{2}}(\tilde{\bm{s}} - \bm{\mu}) \|_{2} \leq \delta \right\},
\]
where $\bm{\mu}$ is the mean of $\bm{s}$ across experts, and $\bm{\Sigma}$ is the covariance matrix of $\bm{s}$ across experts, which is assumed to be $\bm{\Sigma} \succ 0$.
\label{assumption-06}
\end{assumption}

The norm-based ambiguity set represents the variation regions defined by the deviation of $\tilde{\bm{s}}$ from its mean $\bm{\mu}$ across experts, transformed by the covariance matrix $\mathbf{M}$, measured by distance metrics $\|\cdot\|_{2}$. The mean represents the central position of the attribute ranking, while the covariance matrix characterizes the dispersion.

\begin{assumption}[Budget-based ambiguity set] The budget-based ambiguity set for the attribute ranking $\tilde{\bm{s}}$ is defined by
\[
    \mathcal{V}^{\text{budget}} := \left\{ \tilde{\bm{s}} \in \mathbb{R}^{JR}: |\tilde{s}_{g} - \mu_{g}| \leq \gamma_{g}, \forall g \in \mathcal{G}, \sum_{g \in \mathcal{G}}\frac{|\tilde{s}_{g} - \mu_{g}|}{\gamma_{g}} \leq \Gamma, |\mathcal{G}| = JR \right\}.
\]
\label{assumption-07}
\end{assumption}

The budget-based ambiguity set can be regarded as the intersection of $\|\cdot\|_{\infty}$ and $\|\cdot\|_{1}$, which is less conservative as it excludes rare events compared to interval ambiguity set.

\begin{assumption}[CVaR-based ambiguity set] The CVaR-based ambiguity set for attribute ranking $\tilde{\bm{s}}$ is defined by
\[
    \mathcal{V}^{\text{CVaR}} := \left\{ \tilde{\bm{s}} \in \mathbb{R}^{JR}: \exists \bm{\eta} \in \mathbb{R}^{I}, \tilde{\bm{s}} = \sum_{i \in \mathcal{I}} \eta_{i} \bm{s}_{i}, \sum_{i \in \mathcal{I}} \eta_{i} = 1, 0 \leq \eta_{i} \leq \frac{1}{\alpha I}, \forall i \in \mathcal{I} \right\},
\]
where $\alpha \in (0,1]$.
\label{assumption-08}
\end{assumption}

In the CVaR-based ambiguity set, attribute rankings are determined by the convex combination of ranking parameters provided by experts. The upper bound of the weights is related to CVaR, with $\alpha$ controlling the tail risk. A smaller $\alpha$ results in a larger ambiguity set and a more robust solution.

We next derive the tractable reformulation of the second-stage problem of OPA-PR under the above discussed ambiguity sets for attribute rankings. 

\begin{proposition}
    For any $i \in \mathcal{I}$, the second-stage problem of OPA-PR with an ambiguity set for attribute rankings that satisfies Assumption \ref{assumption-06} is equivalent to the following second order cone programming
    \begin{equation}
        \begin{aligned}
            \max_{\tilde{\bm{w}}_{i} \in \bar{\mathcal{W}}_{i}, \tilde{z}_{i}, \bm{\lambda}_{i}} \text{ } &  \tilde{z}_{i}, \\
            \mathrm{s.t.} \text{ } & ( \|(\bm{\Sigma}^{\frac{1}{2}})^{\top}\tilde{\bm{w}}_{i}\|_{2}) \bm{e} \leq \frac{\tilde{\mathbf{W}}_{i}\bm{\mu} - R \tilde{\bm{U}}_{i}^{*} \tilde{z}_{i}}{\delta_{i}}.
        \end{aligned}
        \label{eq-15}
    \end{equation}
    \label{proposition-05}
\end{proposition}

\begin{proposition}
    For any $i \in \mathcal{I}$, the second-stage problem of OPA-PR with an ambiguity set for attribute rankings that satisfies Assumption \ref{assumption-07} is equivalent to the following linear programming
    \begin{equation}
        \begin{aligned}
            \max_{\tilde{\bm{w}}_{i} \in \bar{\mathcal{W}}_{i}, \tilde{z}_{i}, \bm{y}_{1}, \bm{y}_{2}, \bm{\lambda}, \nu} \text{ } &  \tilde{z}_{i}, \\
            \mathrm{s.t.} \text{ } & (\bm{y}_{1} + \bm{y}_{2})^{\top} \bm{\mu} + (\bm{\lambda}^{\top}\bm{e} + \Gamma \nu) \bm{e} \leq - R \tilde{\bm{U}}_{i}^{*} \tilde{z}_{i}, \\
            & \mathbf{Y}_{1} + \mathbf{Y}_{2} = -\tilde{\mathbf{W}}_{i}, \\
            & \bm{\lambda}\geq-\bm{\Delta}_{i}^{-1}\bm{y}_{1}, \\
            & \bm{\lambda}\geq\bm{\Delta}_{i}^{-1}\bm{y}_{1}, \\
            & \nu \bm{e}\geq-\bm{\Delta}_i^{-1}\bm{y}_2, \\
            & \nu \bm{e}\geq\bm{\Delta}_i^{-1}\bm{y}_2,
        \end{aligned}
        \label{eq-16}
    \end{equation}
    where $\mathbf{Y}_{1} = \mathrm{diag}(y_{111}, \dots, y_{11R}, y_{121}, \dots, y_{1JR})$, $\mathbf{Y}_{2} = \mathrm{diag}(y_{211}, \dots, y_{21R}, y_{221}, \dots, y_{2JR})$, and
    \[
        \bm{\Delta}_{i} = \mathrm{diag}(\underbrace{1/\gamma_{i1}, \dots, 1/\gamma_{i1}}_{R \text{ elements}} , \dots, \underbrace{1/\gamma_{iJ}, \dots, 1/\gamma_{iJ}}_{R \text{ elements}}).
    \]
    \label{proposition-06}
\end{proposition}

\begin{proposition}
    For any $i \in \mathcal{I}$, the second-stage problem of OPA-PR with an ambiguity set for attribute rankings that satisfies Assumption \ref{assumption-08} is equivalent to the following linear programming
    \begin{equation}
        \begin{aligned}
            \max_{\tilde{\bm{w}}_{i} \in \bar{\mathcal{W}}_{i}, \tilde{z}_{i}, \lambda, \bm{\beta}} \text{ } &  \tilde{z}_{i}, \\
            \mathrm{s.t.} \text{ } & \left(\lambda - \frac{1}{\alpha I} \sum_{i^{\prime} \in \mathcal{I}} \beta_{i^{\prime}}\right)\bm{e} \geq R \tilde{\bm{U}}_{i}^{*} \tilde{z}_{i}, \\ 
            & \tilde{\mathbf{W}}_{i} \bm{s}_{i^{\prime}} + (- \lambda + \beta_{i^{\prime}})\bm{e} \geq 0,  &&\quad \forall i^{\prime} \in \mathcal{I}, \\
            & \bm{\beta} \geq 0.
        \end{aligned}
        \label{eq-17}
    \end{equation}
    \label{proposition-07}
\end{proposition}

The optimal weights $\tilde{w}_{ijr}^{*}$ are then mapped to $w_{ijk}^{*}$ according to alternative rankings under attributes by experts, with Equation (\ref{eq-05}) calculating the final weights for experts, attributes, and alternatives. Notably, when the worst-case utilities $\tilde{U}_{ijr}^{*}$ for all experts and attributes $(i,j) \in \mathcal{I}\times\mathcal{J}$ reduce to rank order centroid weights, and the ambiguity set for attribute rankings are discrete scenario or box-based, OPA-PR recovers the OPA-R proposed by \citet{M.A.D22}.

The following theorem gives the decomposability of the optimal solution of OPA-PR, given the worst-case attribute ranking parameters. Let $\bm{s}_{i}^{*} = \arg\min_{\tilde{\bm{s}}_{i} \in \mathcal{V}_{i}} \bar{\mathbf{W}}_{i} \tilde{\bm{s}}_{i}$, which can be obtained from the dual optimal solution. For further details about recoverying the primal via the dual, we refer readers to Chapter 5 of \citet{B.V04}.

\begin{theorem}
    For any $i \in \mathcal{I}$, the optimal weights for the second-stage problem of OPA-PR is given by
    \begin{equation}
        w_{ijr}^{*} = w_{ij}^{WR}w_{ijr}^{WU}, \quad \forall (j, r) \in \mathcal{J} \times \mathcal{R},
        \label{eq-18}
    \end{equation}
    where 
    \[
        w_{ij}^{WR} = \frac{1}{s_{ij}^{*} \sum_{j \in \mathcal{J}} \frac{1}{s_{ij}^{*}}}, \quad \forall j \in \mathcal{J},
    \]
    and 
    \[
        w_{ijr}^{WU} = \left(\tilde{U}_{ijr}^{*} \sum_{j \in \mathcal{J}} \frac{1}{s_{ij}^{*}} \right) \Bigg/ \left(\sum_{j\in \mathcal{J}}\sum_{r \in \mathcal{R}} \frac{\tilde{U}_{ijr}^{*}}{s_{ij}^{*}} \right), \quad \forall (j, r) \in \mathcal{J} \times \mathcal{R},
    \]
    with $\sum_{j \in \mathcal{J}} w_{ij}^{WR} = 1$ and $\sum_{j \in \mathcal{J}} \sum_{r \in \mathcal{R}} w_{ijr}^{WU} = 1$.
    \label{theorem-01}
\end{theorem}

Theorem \ref{theorem-01} demonstrates the decomposability of the second-stage solution of OPA-PR, where the two weight components are normalized and depend solely on the worst-case attribute rankings and worst-case utilities. It further suggests that, once the worst-case attribute ranking and worst-case utilities are determined, the optimal weights for the second-stage problem of OPA-PR can be directly computed using Equation \eqref{eq-18}, providing a closed-form solution.

\subsection{Error Bound}
\label{section-3-4}

This section analyzes the error bounds of the ambiguity set approximation of marginal utility functions on the optimal value for the both stages of OPA-PR to provide a theoretical basis for the application. 

Compared to the first-stage error bound of OPA-PR, our primary interest lies in identifying which forms of partial preference information can minimize the first-stage error, as this can inform the design of more effective preference elicitation strategies for DMs. The following proposition shows that the PLA of the first-stage problem in OPA-PR introduces no additional approximation errors when the utility functions in the ambiguity set are concave and the non-differentiable turning points on $\Theta$ corresponds to discrete outcomes.

\begin{proposition}
    For any $(i,j) \in \mathcal{I} \times \mathcal{J}$, if $\psi_{l_{ij}}(\tau)$ for all $l_{ij} \in \mathcal{L}_{ij}$ are step functions over $\Theta$ with jumps at $\tau_{h_{ij}}$ for $h_{ij} \in \mathcal{H}_{ij}$, then $\tilde{\rho}_{ij} = \rho_{ij}$.
    \label{proposition-08}
\end{proposition}

The following proposition gives a equivalent representation of the PLA of the ambiguity sets for marginal utility functions, which is obtained from the step-like approximation of the partial preference information in moment-type preference elicitation set.
\begin{proposition}
    For any $(i,j) \in \mathcal{I} \times \mathcal{J}$, let $u_{ij} \in \mathcal{U}_{ij}$ and $\tilde{u}_{ij}$ be the corresponding piecewise linear approximation defined as Equation \eqref{eq-10}. Then, the step-like approximation of $\psi_{l_{ij}}$ for $l_{ij} \in \mathcal{L}_{ij}$, denoted as $\tilde{\psi}_{l_{ij}}$, resulting from the turning points at $\tau_{h_{ij}}$ for $h_{ij} \in \mathcal{H}$ such that $\tilde{\psi}_{l_{ij}}(\tau_{h_{ij}}) = \psi_{l}(\tau_{h_{ij}})$ is equivalent to the piecewise linear approximation of $u_{ij}$.
    \label{proposition-09}
\end{proposition}

Propositions \ref{proposition-05} and \ref{proposition-06} provide insights for designing effective preference elicitation strategies in OPA-PR: we can avoid approximation errors by properly eliciting the preferences with the step-like characteristic. Notably, as we discussed before, the partial preference information derived from deterministic utility comparison (Example \ref{example-01}) and stochastic lottery comparison with two lotteries (Example \ref{example-02}) is the step function. These findings will aid in developing preference elicitation strategies for OPA-PR, presented in Section \ref{section-3-5}.

We next quantify the impact of the elicited worst-case marginal utility functions with PLA in the first-stage problem on the optimal solution of the second-stage problem of OPA-PR. To achieve this, it suffices to analyze the difference between $\tilde{u}^{*}$ and $u^{*}$ and its impact on the optimal solution, as the PLA solely affects the structure of ambiguity set for utility functions, especially moment-type preference elicitation. Notably, PLA has impact on the second-stage problem even when no extra error is introduced in the first-stage problem. This is because the second-stage problem relies on the arguments from the first-stage minimization problem (i.e., $\tilde{u}_{ij}^{*}$ for all $(i,j) \in \mathcal{Y}$), rather than the optimal value $\tilde{\rho}_{ij}$. We begin by introducing the pseudo-metric for measuring the disparity between utility functions.
\begin{definition}
    Let $\mathscr{F}$ be a set of measurable functions defined over $\Theta$. For any $u_{1},u_{2} \in \mathscr{U}$, the pseudo-metric between $u_{1}$ and $u_{2}$ is defined as
    \[
    d_{\mathscr{F}}(u_{1},u_{2}) := \sup_{f \in \mathscr{F}} | \left<f, u_{1}\right> - \left<f, u_{2}\right> |,
    \]
    where $\left< \cdot \right>$ denotes the Euclidean inner product.
    \label{definition-02}
\end{definition}

It is easy to verify that $d_{\mathscr{F}}(u_{1},u_{2}) = 0$ if and only if $\left<f, u_{1}\right> = \left<f, u_{2}\right>$ for all $f \in \mathscr{F}$. Based on the pseudo-metric, the following lemma provides the upper bound for the PLA errors between marginal utility functions.
\begin{lemma}
    Let $\mathscr{F} := \left\{ f = I_{\Theta}(\cdot) \right\}$, where $I_{\Theta}(\cdot)$ is the indicator function with domain $\Theta$. For any $(i,j) \in \mathcal{I} \times \mathcal{J}$, the optimal $u^{*}_{ij} \in \mathcal{U}_{ij}$ is $G$-Lipschitz continuous over $\Theta$, and $\tilde{u}^{*}_{ij}$ represents the corresponding piecewise linear approximation defined as Equation \eqref{eq-10}. Then 
    \begin{equation}
        d_{\mathscr{F}}(u^{*}_{ij}, \tilde{u}^{*}_{ij}) = \sup_{\tau \in \Theta} \left| \tilde{u}_{ij}^{*}(\tau) - u_{ij}^{*}(\tau) \right| \leq G \zeta_{ij}, 
        \label{eq-19}
    \end{equation}
    where $\zeta_{ij} = \max_{h_{ij} \in \mathcal{H}_{ij}} (\tau_{h_{ij}} - \tau_{h_{ij}-1})$ and $G \geq 1/R$.
    \label{lemma-02}
\end{lemma}

With Lemma \ref{lemma-02}, we can quantify the error bound on the optimal solution of the second-stage problem in OPA-PR originating from the PLA of the elicited worst-case marginal utility functions. Consider the following dual problem of the second-stage problem with utilities derived from the unapproximated marginal utility functions for any $i \in \mathcal{I}$ with the dual variable $\theta_{i} \in \mathbb{R}$ and $\bm{\lambda}_{i} \in \mathbb{R}^{J \times R}$
\begin{equation}
    \begin{aligned}
        \min_{\theta_{i}, \bm{\lambda}_{i}} \text{ } & \theta_{i}, \\
        \mathrm{s.t.} \text{ } & \bm{\Lambda}_{i} \tilde{\bm{s}}_{i}^{*} \leq \theta_{i} \bm{e}^{\top}, \\
        & R (\bm{U}_{i}^{*})^{\top} \bm{\lambda}_{i} = 1, \\
        & \bm{\lambda}_{i} \geq 0,
    \end{aligned}
    \label{eq-20}
\end{equation}
where $\bm{\Lambda}_{i} = \mathrm{diag}(\lambda_{i11}, \dots, \lambda_{i1R}, \lambda_{i21}, \dots, \lambda_{iJR})$.

\begin{theorem}
    For any $(i,j) \in \mathcal{I} \times \mathcal{J}$, assume that $\mathscr{F} := \left\{ f = I_{\Theta}(\cdot) \right\}$ and $u_{ij}^{*} \in \mathcal{U}_{ij}$ is $G$-Lipschitz continuous over $\Theta$ and $\tilde{u}_{ij}^{*}$ is corresponding piecewise linear approximation defined as Equation \eqref{eq-10}. Let $(\theta_{i}^{*}, \bm{\lambda}_{i}^{*})$ denote the optimal solution of the dual problem of Equation \eqref{eq-20}. Let $(\tilde{z}_{i}^{*}, \tilde{\bm{w}}_{i}^{*})$ and $({z}_{i}^{*}, \bar{\bm{w}}_{i}^{*})$ denote the optimal solutions of the second-stage problem of OPA-PR, corresponding to unapproximated and approximated utilities, respectively. Then, 
    \begin{equation}
         (1 - R G \zeta_{ij}(\bm{\lambda}_{i}^{*})^{\top} \bm{e}){z}_{i}^{*} \leq \tilde{z}_{i}^{*} \leq \frac{1}{1 + R G \zeta_{ij}(\bm{\lambda}_{i}^{*})^{\top} \bm{e}} {z}_{i}^{*},
        \label{eq-21}
    \end{equation}
    and for all $(j,r) \in \mathcal{J} \times \mathcal{R}$,
    \begin{equation}
        \frac{(1 - RG\zeta_{ij}(\bm{\lambda}_{i}^{*})^{\top}\bm{e})R \tilde{U}_{ijr}^{*}}{\tilde{s}_{ij}^{*}}z_{i}^{*} \leq \tilde{w}_{ijr}^{*} \leq \frac{R \tilde{U}_{ijr}^{*}}{(1 + RG\zeta_{ij}(\bm{\lambda}_{i}^{*})^{\top}\bm{e})\tilde{s}_{ij}^{*}}z_{i}^{*},
        \label{eq-22}
    \end{equation}
    where $\zeta_{ij} = \max_{h_{ij} \in \mathcal{H}_{ij}} (\tau_{h_{ij}} - \tau_{h_{ij}-1})$.
    \label{theorem-02}
\end{theorem}

\subsection{Preference Elicitation Strategy Design}
\label{section-3-5}

This section focuses on the preference elicitation strategy relevant to the first-stage problem of OPA-PR. 
As discussed in Section \ref{section-3-4}, the way to elicit expert preferences directly influences the PLA errors of the ambiguity set for marginal utility functions. Thus, we utilize the stochastic pairwise lottery comparisons, which provide step-like preference information covering the deterministic utility comparisons, without introducing additional errors in the PLA as detailed in Propositions \ref{proposition-08} and \ref{proposition-09}. Specifically, experts are asked to determine the certainty equivalent of pairwise lotteries
\[
    Z_1=\left\{\begin{array}{lll}r_1,&\text{with probability}&1-p,\\r_3,&\text{with probability}&p,\end{array}\right.\quad\text{and}\quad Z_2=r_2.
\]
For any $(i,j) \in \mathcal{I} \times \mathcal{J}$, the question $l_{ij}$ is characterized by four parameters $r_{l_{ij}}^{1} \leq r_{l_{ij}}^{2} \leq r_{l_{ij}}^{3}$ and $p_{ij}$ with the corresponding $\psi_{l_{ij}}$ given by 
\[
    \psi_{l_{ij}}(\tau) = (1-p_{ij}) \mathbf{1}_{\left[r_{l_{ij}}^{1}, R\right]}(\tau) + p_{ij} \mathbf{1}_{\left[r_{l_{ij}}^{3}, R\right]}(\tau) - \mathbf{1}I_{\left[r_{l_{ij}}^{2}, R\right]}(\tau).
\]
If the utility function is normalized such that $\tilde{u}_{ij}(r_{l_{ij}}^{1}) = 0$ and $\tilde{u}_{ij}(r_{l_{ij}}^{3}) = 1$, the question simplifies to determining whether $\tilde{u}_{ij}(r_{l_{ij}}^{2}) \geq p_{ij}$ or not. Thus, the key is how to determine the probability $p_{ij}$.

In the following, we utilize the random utility split scheme to select these four parameters \citep{A.D15}. We first determine the number of questions $L_{ij}$, from which the number of breakpoints is $3L_{ij} + 2 \leq R$.
\begin{itemize}
    \item Initialization: Set $l_{ij} = 0$.
    \item Step 1: Choose $r_{l_{ij}}^{1}$ and $r_{l_{ij}}^{3}$ uniformly from the ranking set $\mathcal{R} = \{1, 2, \dots, R\}$ and $r_{l_{ij}}^{2}$ uniformly from the rankings in $[r_{l_{ij}}^{1}, r_{l_{ij}}^{3}]$. 
    \item Step 2: Normalize the utility function such that $\tilde{u}_{ij}(r_{l_{ij}}^{1}) = 0$ and $\tilde{u}_{ij}(r_{l_{ij}}^{3}) = 1$. Then, let $C_{ij} := [\underline{C}_{ij}, \bar{C}_{ij}] := \left[\min_{\tilde{u}_{ij} \in \tilde{\mathcal{U}}_{ij}} \tilde{u}_{ij}(r_{l_{ij}}^{2}), \max_{\tilde{u}_{ij} \in \tilde{\mathcal{U}}_{ij}} \tilde{u}_{ij}(r_{l_{ij}}^{2})\right]$ and $p_{ij} = \left( \underline{C}_{ij} + \bar{C}_{ij} \right)/2$. Taking the elicitation of the minimum $\underline{C}_{ij}$ as an example, which is solved by the following minimization problem
    \[
    \begin{aligned}
        \min_{\bm{a} \succeq 0, \bm{b}} \text{ } & \sum_{h_{ij} \in \mathcal{H}_{ij}} (a_{h_{ij}}r_{l_{ij}}^{2} + b_{h_{ij}}) \mathbf{1}_{(\tau_{h_{ij}-1}, \tau_{h_{ij}}]}(r_{l_{ij}}^2), \\
        \mathrm{s.t.} \text{ } & a_{r_{l_{ij}}^{1}-1} r_{l_{ij}}^{1} + b_{r_{l_{ij}}^{1}-1} = 0, \\
        & a_{r_{l_{ij}}^{3}-1} r_{l_{ij}}^{3} + b_{r_{l_{ij}}^{3}-1} = 1, \\
        & a_{h_{ij}-1}\tau_{h_{ij}} + b_{h_{ij}-1} = a_{h_{ij}} \tau_{h_{ij}} + b_{h_{ij}}, \quad && \forall h_{ij} \in \mathcal{H}_{ij}, \\
        & a_{h_{ij}} - a_{h_{ij}-1} \leq 0, \quad && \forall h_{ij} \in \mathcal{H}_{ij}, \\
        & \sum_{h_{ij} \in \mathcal{H}_{ij}} a_{h_{ij}-1} \int_{\tau_{h_{ij}}}^{\tau_{h_{ij}-1}} \psi_{l_{ij}}(\tau) \mathrm{d}t \leq c_{l_{ij}}, \quad && \forall l_{ij} \in \mathcal{L}_{ij}, \\
        & a_{h_{ij}-1} \leq G \quad && \forall h_{ij} \in \mathcal{H}_{ij}.
    \end{aligned}
    \]
    Similarly, $\bar{C}_{ij}$ can be determined by maximizing the above problem. 
    \item Step 3: Let $l_{ij} = l_{ij}+1$ and $\psi_{l_{ij}}(\tau) = (1-p_{ij}) \mathbf{1}_{[r_{l_{ij}}^{1}, R]}(\tau) + p_{ij} \mathbf{1}_{[r_{l_{ij}}^{3}, R]}(\tau) - \mathbf{1}_{[r_{l_{ij}}^{2}, R]}(\tau)$.
    If $p_{ij} \geq \tilde{u}_{ij}(r_{l_{ij}}^{2})$, which is equivalent to $\int_{\Theta} - \psi_{l_{ij}}(\tau) \mathrm{d} \tilde{u}_{ij}(\tau) \leq 0$, then the expert prefers $Z_{1}$ to $Z_{2}$. We then include the constraint $\int_{\Theta} - \psi_{l_{ij}}(\tau) \mathrm{d} \tilde{u}_{ij}(\tau) \leq 0$ in the ambiguity set $\tilde{\mathcal{U}}_{ij}$. Otherwise, add $\int_{\Theta} \psi_{l_{ij}}(\tau) \mathrm{d} \tilde{u}_{ij}(\tau) \leq 0$. Return to Step 1 until $l_{ij} = L_{ij}$.
\end{itemize}

Step 1 generates the lottery and certainty equivalent for pairwise comparisons. Step 2 determines the probability $p_{ij}$ at the midpoint of interval $C_{ij}$, which means the true utility function value at $r_{l_{ij}}^{2}$ lies in either the upper or lower half of interval $C_{ij}$, reducing the ambiguity set for possible utility functions by half. Step 3 requires the expert to select the lottery and certainty equivalent to add the constraints. Notably, when the expert specifies the nominal utility function, the modulus $G$ of Lipschitz continuity modulus can be derived from its first-order derivative.

\begin{remark}
    It is important to acknowledge that preference inconsistencies may occur during the elicitation process, arising from factors such as violations of expected utility theory axioms or contamination of preference information \citep{G.X.Z24}. In \ref{Asection-01}, we present extended formulations for two types of preference inconsistencies commonly discussed in literature \citep{G.X.Z24, W.H.H.X22}: weight disparity error and erroneous elicitation.
\end{remark}

\section{Numerical Experiments}
\label{section-4}

\subsection{Experiment Setup}
\label{section-4-1}

In this section, we conduct numerical experiments on OPA-PR, applied to the emergency supplier selection problem in the context of the 7.21 mega-rainstorm disaster in Zhengzhou, China. We focus on the in-disaster scenario, which contrasts with pre-disaster conditions. Unlike the pre-disaster setting, the in-disaster phase demands rapid decision-making due to the high urgency and uncertainty that typically follow such events. This scenario challenges traditional supplier selection as the unpredictability and complexity of disasters often make preselected suppliers inadequate \Citep{W.W.L.L22}. In response, the in-disaster emergency supplier selection problem must adapt by relying on partial preference information and accommodating the risk preferences of DM, significantly influencing outcomes under these conditions. Overall, the in-disaster emergency supplier selection problem, marked by urgency and complexity in decision-making, demonstrates the applicability and effectiveness of OPA-PR, aligning with the proposed decision framework and providing an ideal testing ground.

There are ten emergency suppliers (labeled A1 to A10) available, each presenting unique characteristics suited for flood relief efforts. Some provide stable supply chains and rapid response capabilities at higher costs, while others leverage strategic locations and efficient transport networks to improve emergency effectiveness. Six attributes are identified to evaluate these suppliers: response speed (C1), delivery reliability (C2), geographic coverage (C3), operational sustainability (C4), collaborative experience and credibility (C5), and supply cost (C6). A single-expert decision-making process is adopted in this study without loss of generality, as Corollary \ref{corollary-01} shows that the multi-expert OPA model can be equivalently reduced to a single-expert OPA model followed by aggregation using rank reciprocal weights. This not only simplifies the analysis but also yields results that are more intuitive and interpretable. The expert ranks each attribute and supplier, with the ranking data available in the Online Supplemental Material. The selection of the ambiguity set parameter for OPA-PR follows the method outlined in relevant classic textbooks. For the norm-based ambiguity set, let $\delta = 0.00795$; for the budget-based ambiguity set, set $\gamma = 1$ and $\Gamma = 0.875$; and for the CVaR-based ambiguity set, set $\alpha = 0.95$. For simplicity, OPA-PR with norm-, budget-, and CVaR-based ambiguity sets are referred to as OPA-PR(N), OPA-PR(B), and OPA-PR(C), respectively, in the numerical experiments.

The numerical experiments are structured into three parts: the case test, the sensitivity test, and the comparison test. The case test is designed to demonstrate the implementation of our proposed approach and provide an interpretable explanation of the case results. The sensitivity test assesses the effect of parameter variations on the outcomes. The comparison test contrasts the results with the original OPA model to assess the effects of the extensions on the outcomes.

\subsection{Experiment Results}
\label{section-4-2}

\begin{figure}[h]
    \centering
    \includegraphics[width=0.4\textwidth]{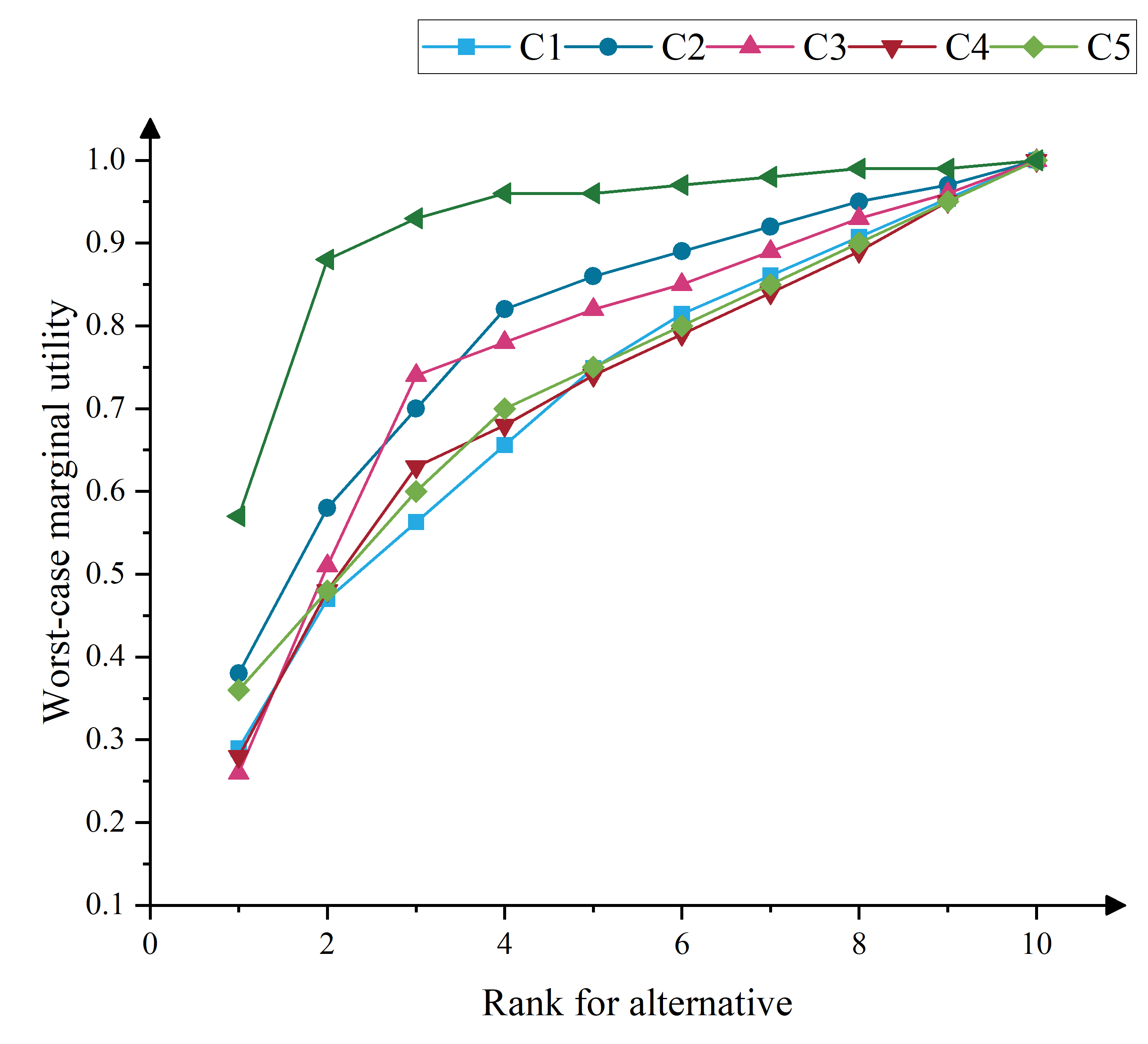}
    \caption{Worst-case marginal utility function (reversed rankings) for ranked alternatives across attributes}
    \label{fig-01}
\end{figure}

This section analyzes the case results of OPA-PR using the collected data. Figure \ref{fig-01} illustrates the worst-case marginal utilities of ranked alternatives under various attributes based on partial preference information obtained through the scheme in Section \ref{section-3-5}, which is the first-stage outcome of OPA-PR. The results indicate that the worst-case marginal utilities of ranked alternatives varies across attributes based on the partial preference information, even when the expert responds to identical utility elicitation questionnaires for different attributes. This highlights the need to customize the preferences of expert in decision-making processes. Furthermore, the worst-case marginal utility function reflects a risk-averse preference, contrasting with the original OPA model using rank order centroid weights, which suggests a risk-seeking preference, extending beyond the less common single risk-seeking preference in practice.

\begin{figure}[h]
    \centering
    \includegraphics[width=0.64 \textwidth]{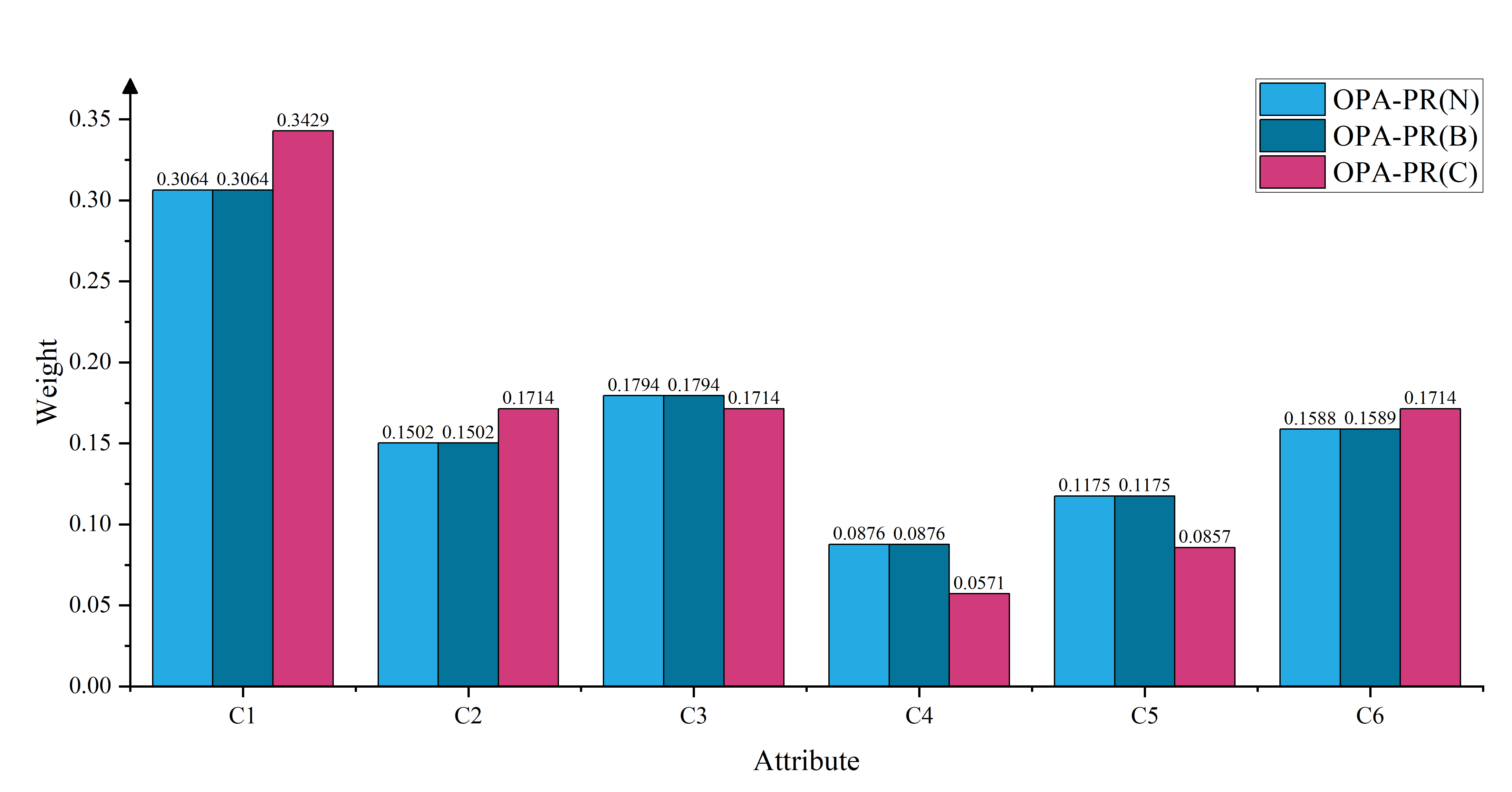}
    \caption{Optimal weights for attributes}
    \label{fig-02}
\end{figure}

\begin{figure}[h]
    \centering
    \includegraphics[width= 0.97\textwidth]{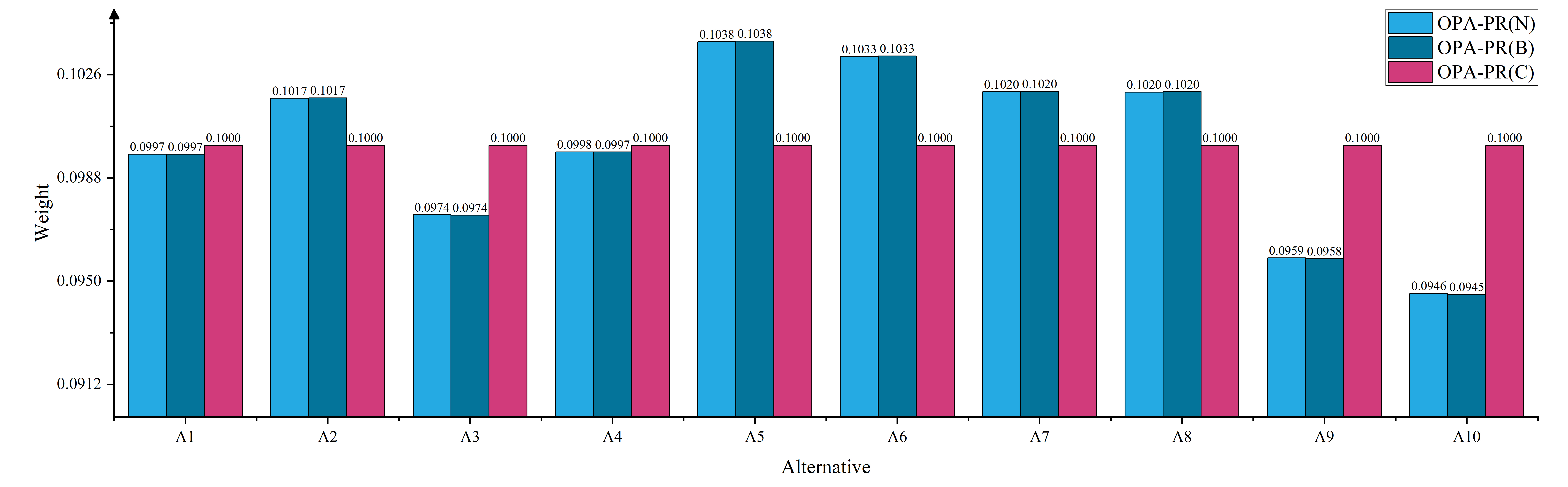}
    \caption{Optimal weights for alternatives}
    \label{fig-03}
\end{figure}

The second-stage problem of OPA-PR under the three ambiguity sets is solved using the elicited marginal utility functions. The optimal weight disparities for OPA-PR(N), OPA-PR(B), and OPA-PR(C) are consistent, with a value of 0.0025. Figures \ref{fig-02} and \ref{fig-03} present the optimal attribute and alternative weights, showing nearly identical results result from the same optimal weight disparities. Regarding attribute weights, OPA-PR(C) exhibits greater differences across attributes than OPA-PR(N) and OPA-PR(B), which show the same attribute weights. C1 is the most significant for all three, with a weight of 0.3064. For OPA-PR(N) and OPA-PR(B), C3 ranks second with a weight of 0.1794, followed by C6, C2, C5, and C4, with weights of 0.1588, 0.1502, 0.1175, and 0.0876, respectively. In contrast, OPA-PR(C) assigns weights of 0.1714 to C2, C3, and C6, followed by 0.0857 for C5 and 0.0571 for C4. The alternative weights for OPA-PR(C) are all identical, each valued at 0.1000, indicating a degenerate solution and suggesting its inapplicability to this scenario. Due to its poor performance, OPA-PR(C) will be excluded from the sensitivity and comparative analysis for further testing. For OPA-PR(N) and OPA-PR(B), the alternative ranking remains consistent. The top six alternatives are A5 (0.1038), A6 (0.1033), A7 (0.1020), A8 (0.1020), and A2 (0.1017), with A5 being the most favorable due to its performance in C1, C2, C4, and C5.

\subsection{Sensitivity Analysis}
\label{section-4-3}

This section performs a sensitivity analysis on OPA-PR, focusing on the perturbation of alternative rankings and size parameters of attribute ranking ambiguity sets, particularly examining the sensitivity of alternative rankings to assess the value of the prescribed alternative in the perturbed version of the perceived ranking.

For the sensitivity analysis of alternative rankings, we generate the samples from a normal distribution using the expert-provided alternative rankings as the mean, with standard deviations of $1/3$, $1/2$, and $2/3$. These standard deviations are based on the fact that approximately 99.7\% of the data within a normal distribution falls within three standard deviations of the mean, allowing estimation of the standard deviations using the maximum and minimum values. A standard deviation of $c/3$ results in a maximum range of $2c$, i.e., $(\text{mean} - c, \text{mean} + c)$, with a radius of $c$. Thus, a standard deviation of $1/3$ represents alternative ranking perturbation without reversal, while $1/2$ and $2/3$ indicate perturbation with ranking reversal. The stopping condition is defined as
\[
    \| (\arg\max_{k \in \mathcal{K}} p_{k1}^{m}, \dots, \arg\max_{k \in \mathcal{K}} p_{kR}^{m})^{\top} - (k_{1}^{*}, \dots, k_{R}^{*})^{\top} \|_{2} = 0 \quad \text{and} \quad \|(\max_{k\in \mathcal{K}} p_{k1}^{m}, \dots, \max_{k\in \mathcal{K}} p_{kR}^{m})^{\top}\|_{2} \geq d,
\] 
with a maximum of 1000 iterations, where $m$ denotes the iteration step, $d$ is a given constant, $k_{r}^{*}$ denotes the index of the $r$-th ranked alternative from the unperturbed OPA-PR, and $p_{kr}^{m}$ represents the probability of alternative $k$ being ranked $r$ for all $k \in \mathcal{K}$ and $r \in \mathcal{R}$ under iteration step $m$. The first condition ensures that the final ranking result converges to the mean-based alternative ranking result after a specified number of iterations. The second condition requires the probability of convergence to meet a minimum threshold. Let $d = 1.5811$, derived from $\|K(0.5, \dots, 0.5)^{\top}\|_{2}$. To make the analysis statistically meaningful, we set the minimum iterations to 50.

\begin{figure}[h]
    \centering
    \subfigure[OPA-PR(N) with radius 1] {
        \includegraphics[width=0.32\columnwidth]{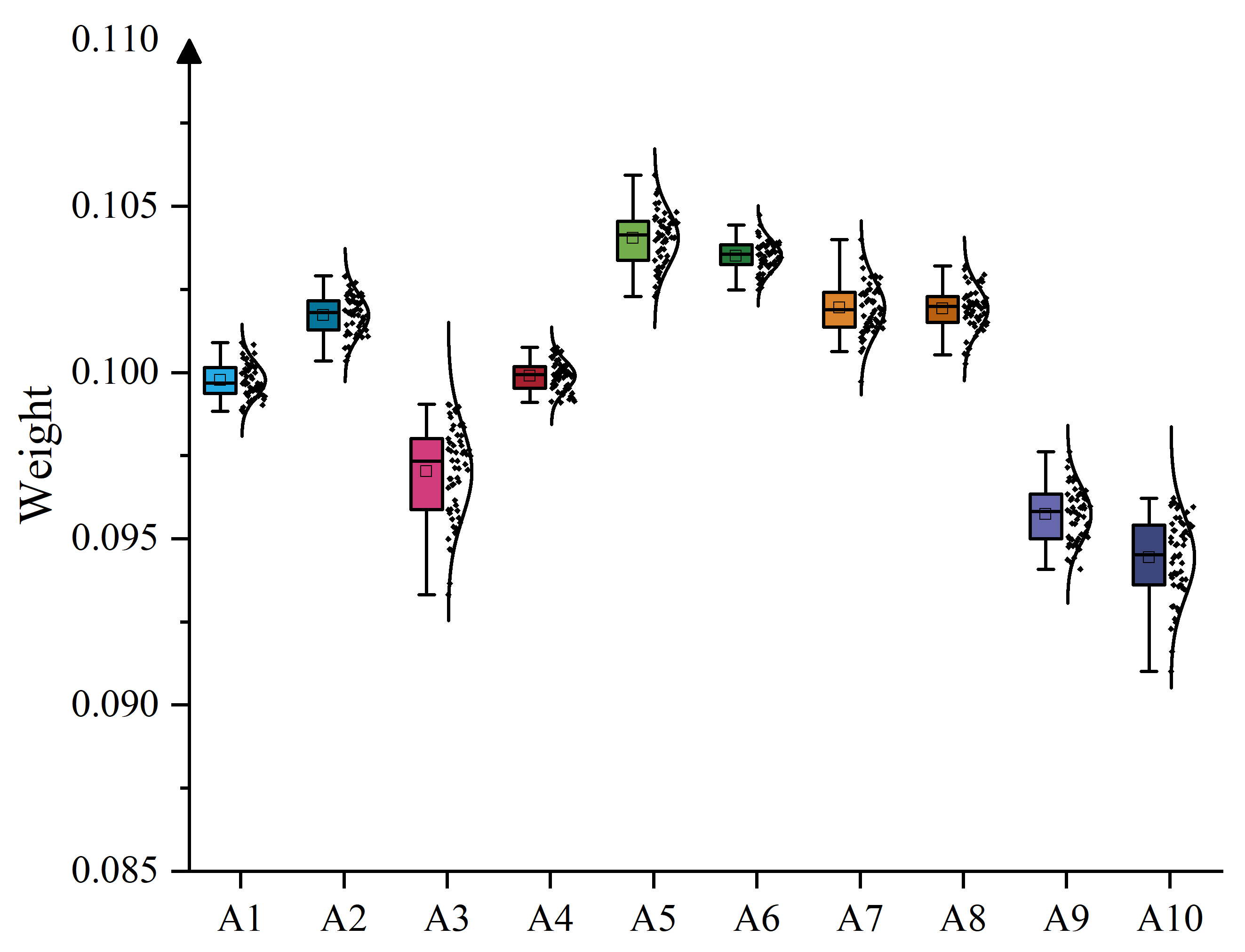}}
    \subfigure[OPA-PR(N) with radius 1.5] {
        \includegraphics[width=0.32\columnwidth]{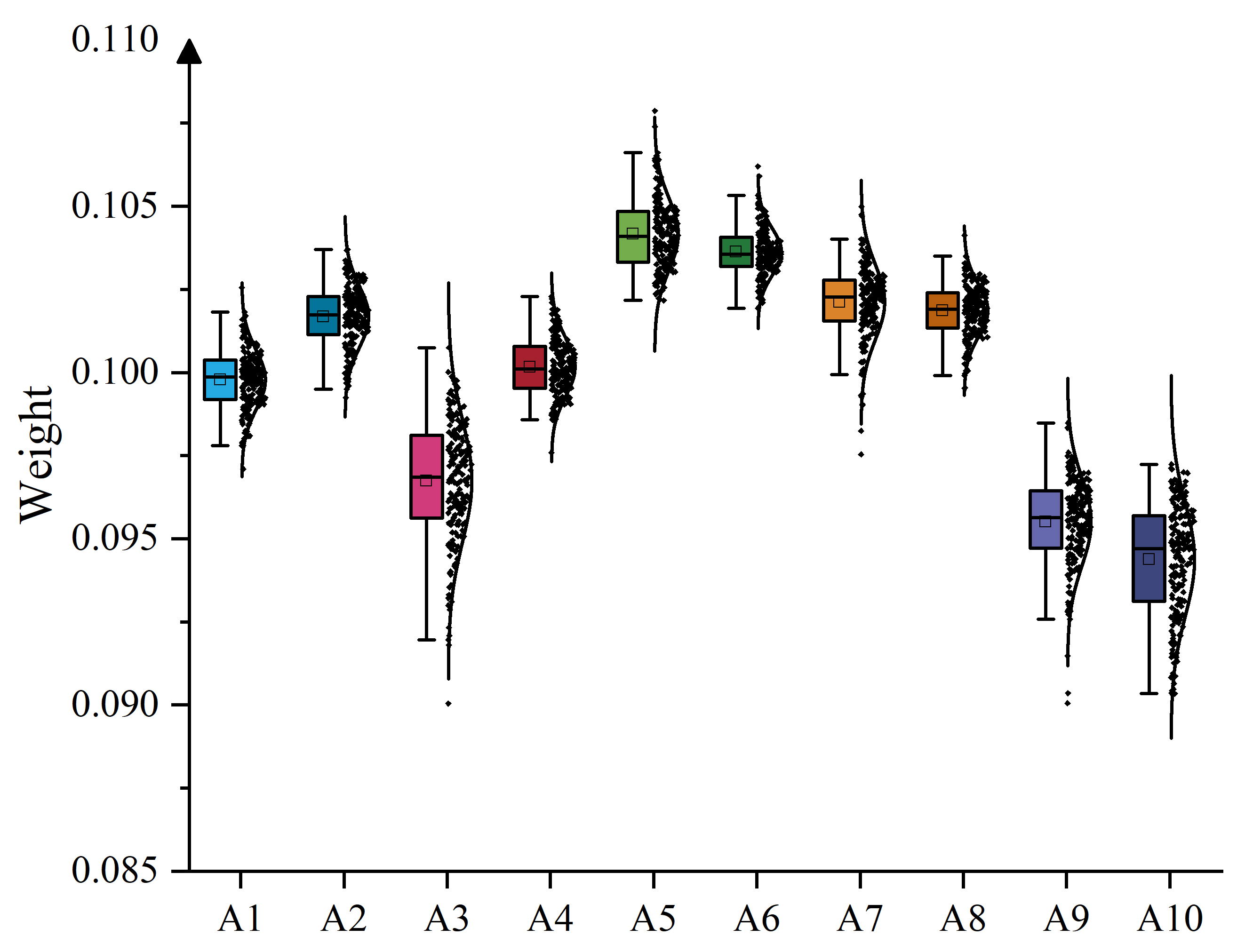}}        
    \subfigure[OPA-PR(N) with radius 2] {
        \includegraphics[width=0.32\columnwidth]{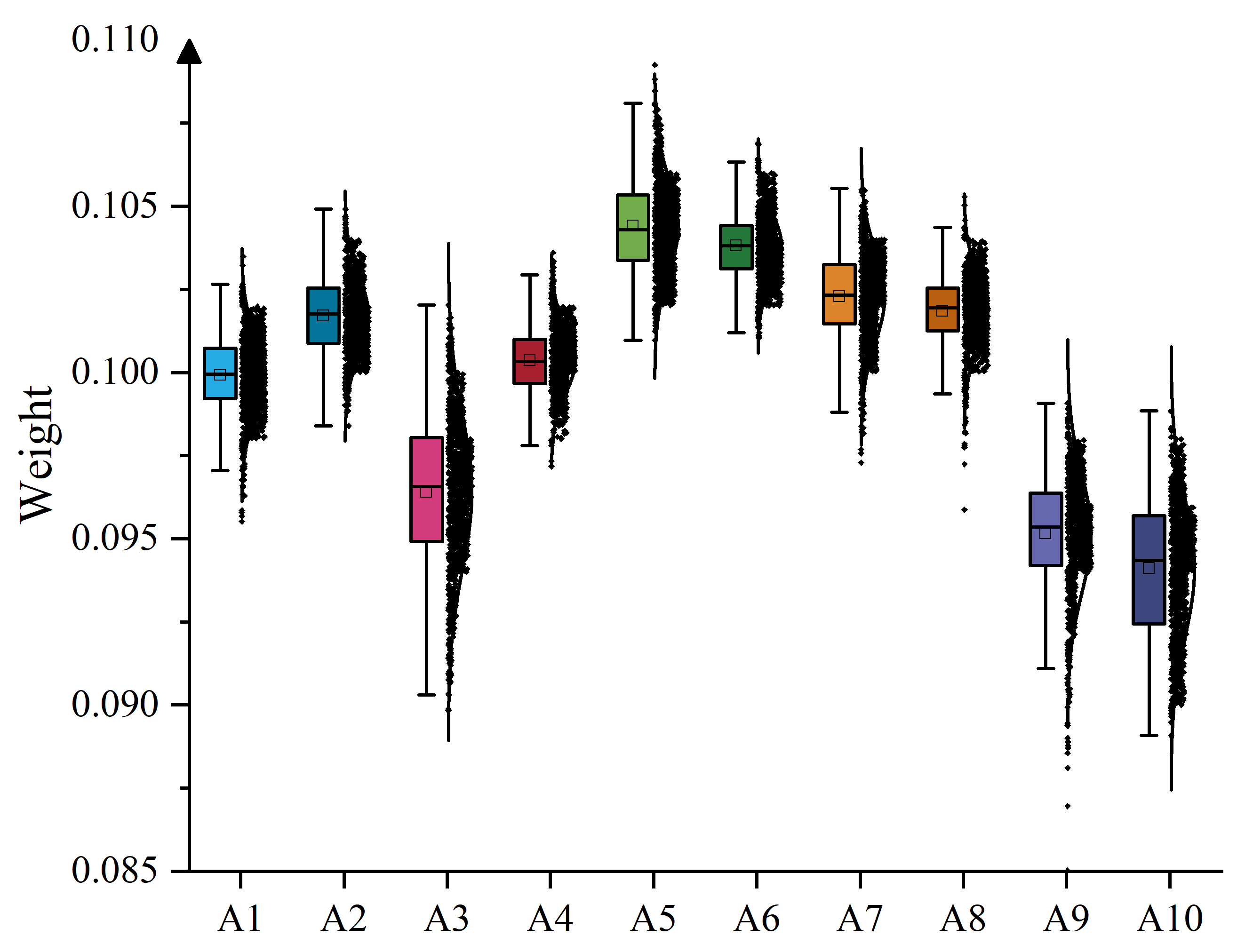}}       \\
    \subfigure[OPA-PR(B) with radius 1] {
        \includegraphics[width=0.32\columnwidth]{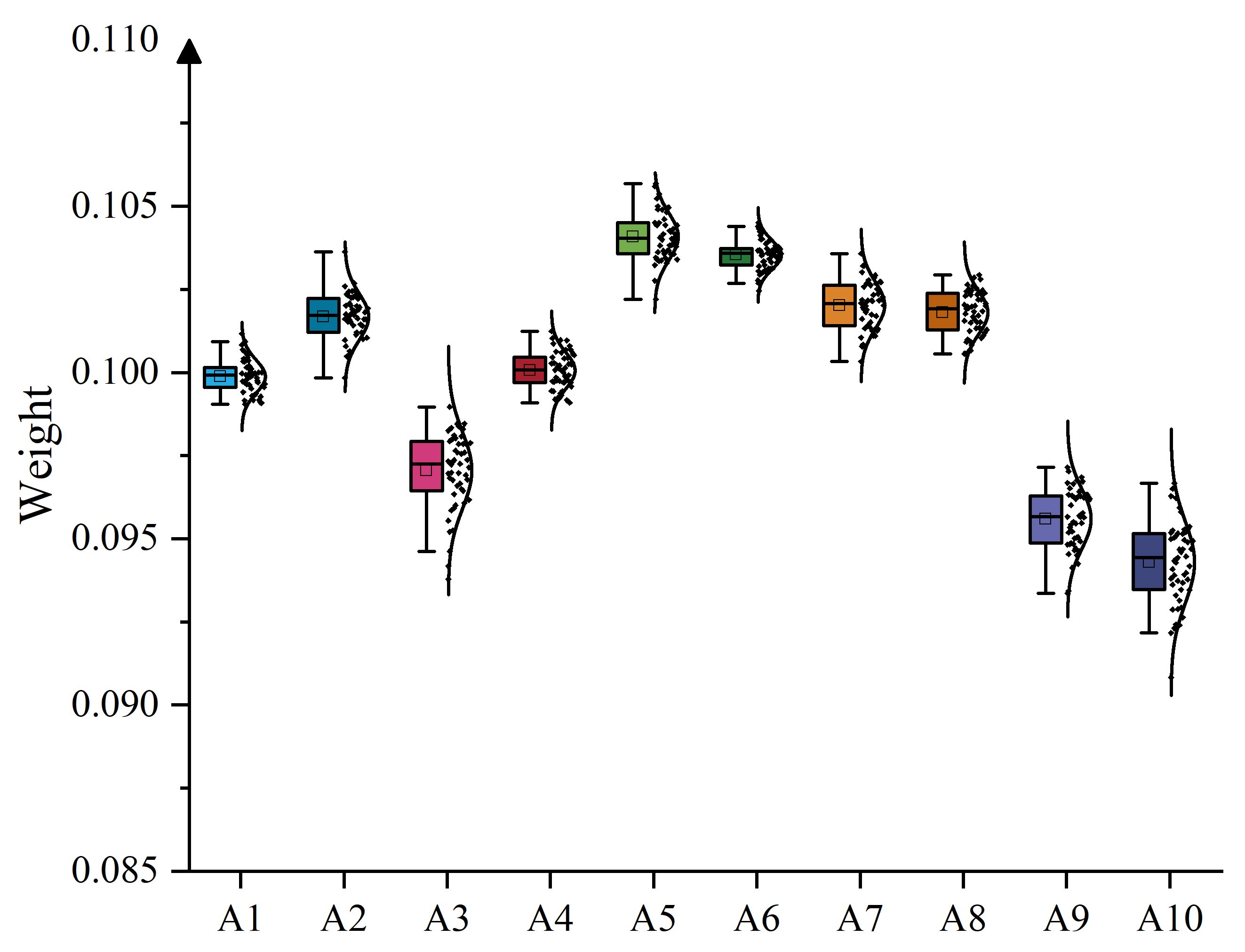}}
    \subfigure[OPA-PR(B) with radius 1.5] {
        \includegraphics[width=0.32\columnwidth]{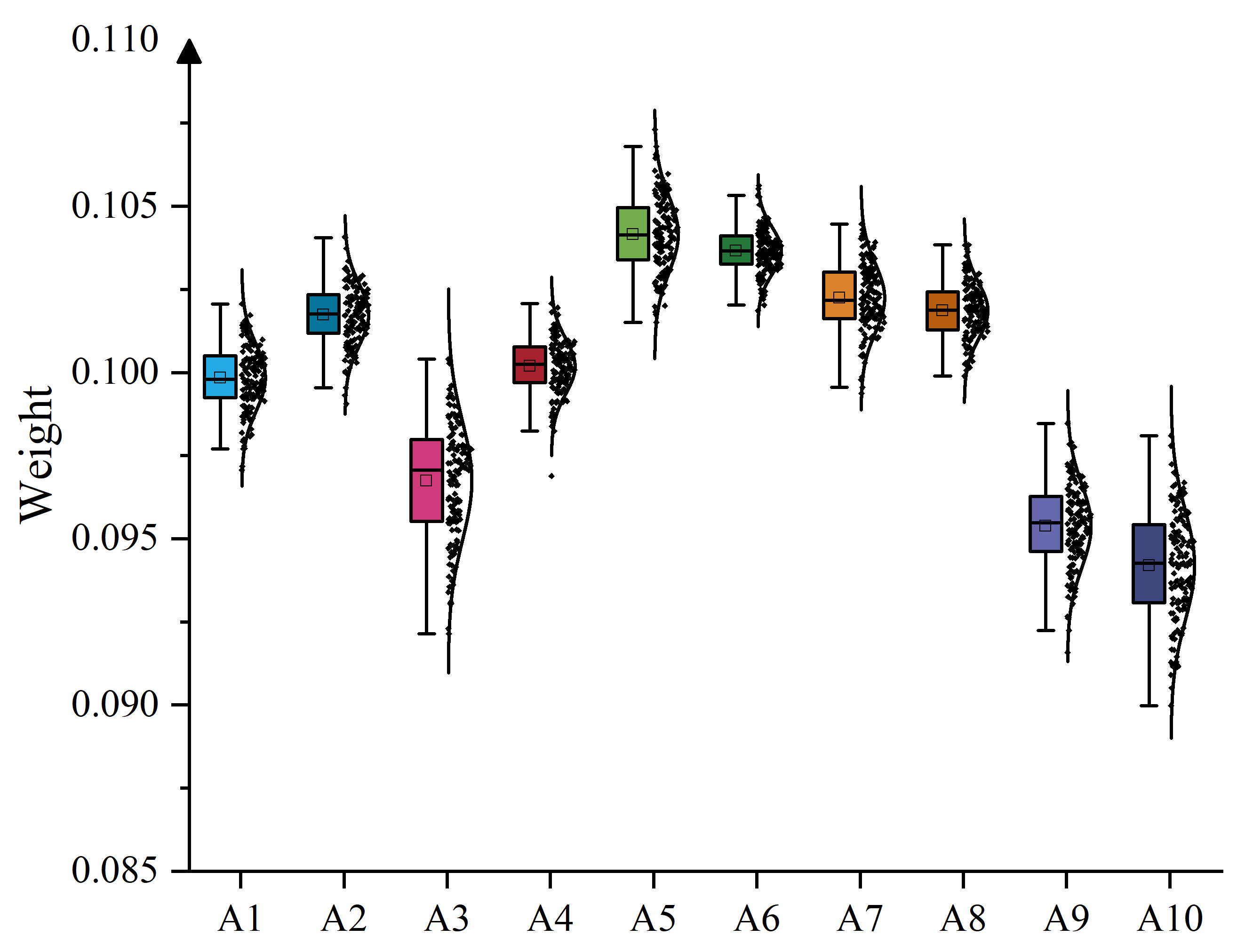}}        
    \subfigure[OPA-PR(B) with radius 2] {
        \includegraphics[width=0.32\columnwidth]{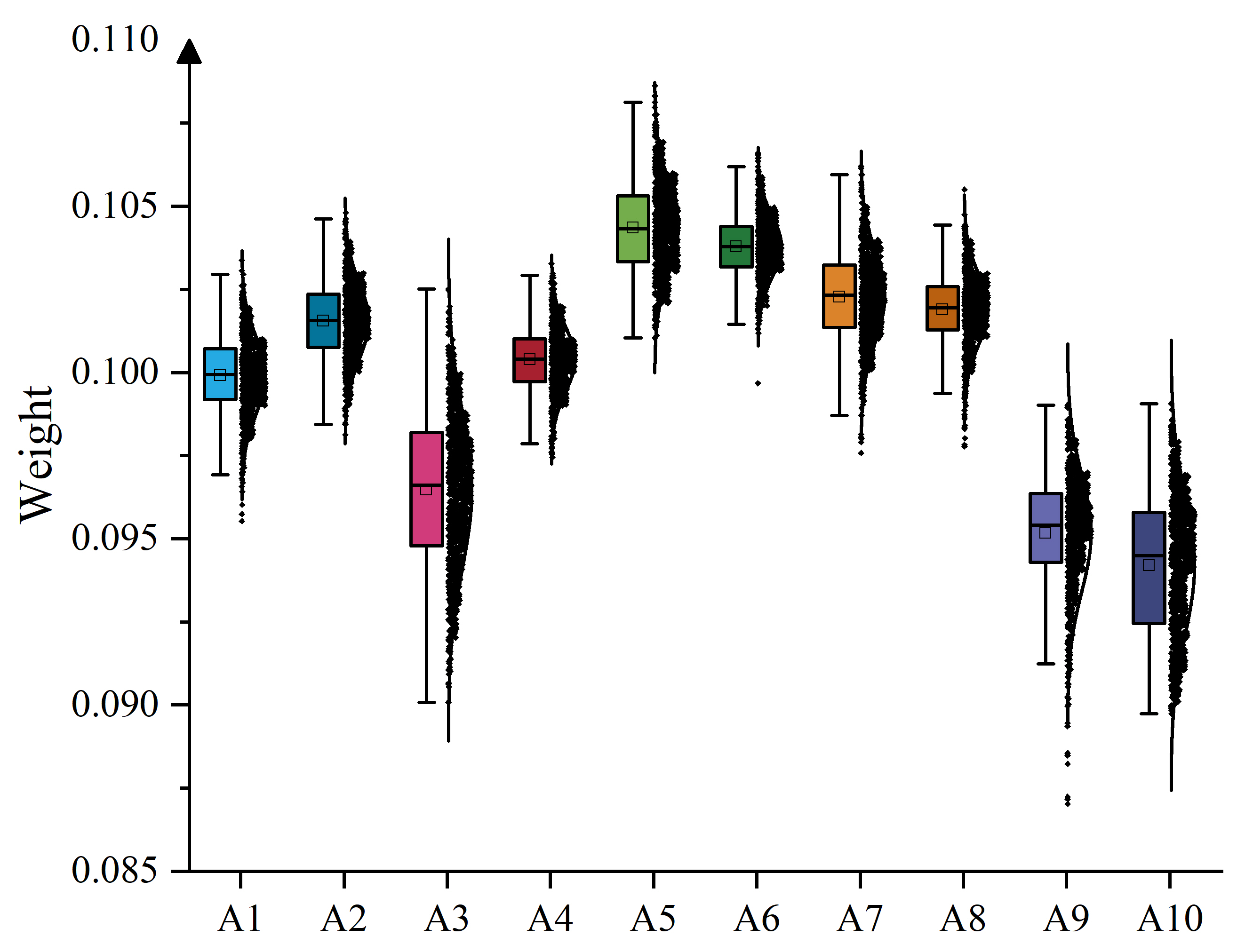}}
    \caption{Weight results of alternative ranking perturbations}
    \label{fig-04}
\end{figure}

Table \ref{tab-01} displays the number of iterations required to meet the stopping condition under different perturbation radii for alternative rankings, along with the final ranking result and the corresponding probabilities. At a radius of 1, both OPA-PR(N) and OPA-PR(B) converge within approximately 50 iterations. With a radius of 1.5, convergence occurs within 200 iterations. However, at a radius of 2, the maximum iteration limit is reached without convergence, primarily because the norm of probability vector does not meet the given minimum threshold. Regarding ranking results, both OPA-PR(N) and OPA-PR(B) produce consistent rankings across radii, all converging to the mean-based ranking result. Additionally, as the perturbation radius increases, the probability of an alternative achieving the optimal ranking decreases. Figure \ref{fig-04} shows the alternative weights across various perturbation radii, with a consistent vertical axis for comparison. The weight distributions follow a normal distribution, in line with the sampling rule. As the radius increases, weight variability becomes more noticeable. It is important to note that due to the stochastic nature of the simulation, slight variations may occur between different runs. Nevertheless, as iterations increase, the final ranking result of OPA-PR converge to the mean-based alternative ranking result, reinforcing the validity of OPA-PR.

\begin{table}[h]
    \centering
    \renewcommand{\arraystretch}{1}
    \caption{Ranking results of alternative ranking perturbations}
    \label{tab-01}
    \begin{tabular}{@{}lllllllllllll@{}}
    \toprule
    \multirow{2}{*}{Model}       & \multirow{2}{*}{Radius} & \multirow{2}{*}{Iteration} & 1st  & 2nd  & 3rd  & 4th  & 5th  & 6th  & 7th  & 8th  & 9th  & 10th \\ \cmidrule(l){4-13} 
                                  &                        &                            & A5   & A6   & A7   & A8   & A2   & A4   & A1   & A3   & A9   & A10  \\ \midrule
    \multirow{3}{*}{OPA-PR(N)}   & 1                      & 56                         & 75\% & 72\% & 33\% & 40\% & 33\% & 51\% & 54\% & 77\% & 60\% & 67\% \\
                                  & 1.5                    & 198                        & 63\% & 57\% & 40\% & 36\% & 31\% & 46\% & 51\% & 56\% & 51\% & 60\% \\
                                  & 2                      & 1000                      & 59\% & 46\% & 28\% & 34\% & 25\% & 39\% & 44\% & 49\% & 48\% & 54\% \\
    \multirow{3}{*}{OPA-PR(B)} & 1                      & 50                         & 71\% & 71\% & 45\% & 41\% & 31\% & 47\% & 57\% & 78\% & 57\% & 75\% \\
                                  & 1.5                    & 143                        & 60\% & 51\% & 36\% & 36\% & 28\% & 49\% & 51\% & 59\% & 58\% & 63\% \\
                                  & 2                      & 1000                      & 58\% & 45\% & 27\% & 34\% & 30\% & 40\% & 48\% & 48\% & 47\% & 52\% \\ \bottomrule
    \end{tabular}
\end{table}

\begin{figure}[h]
    \centering
    \subfigure[Attribute weights of OPA-PR(N)] {
        \includegraphics[width=0.32\columnwidth]{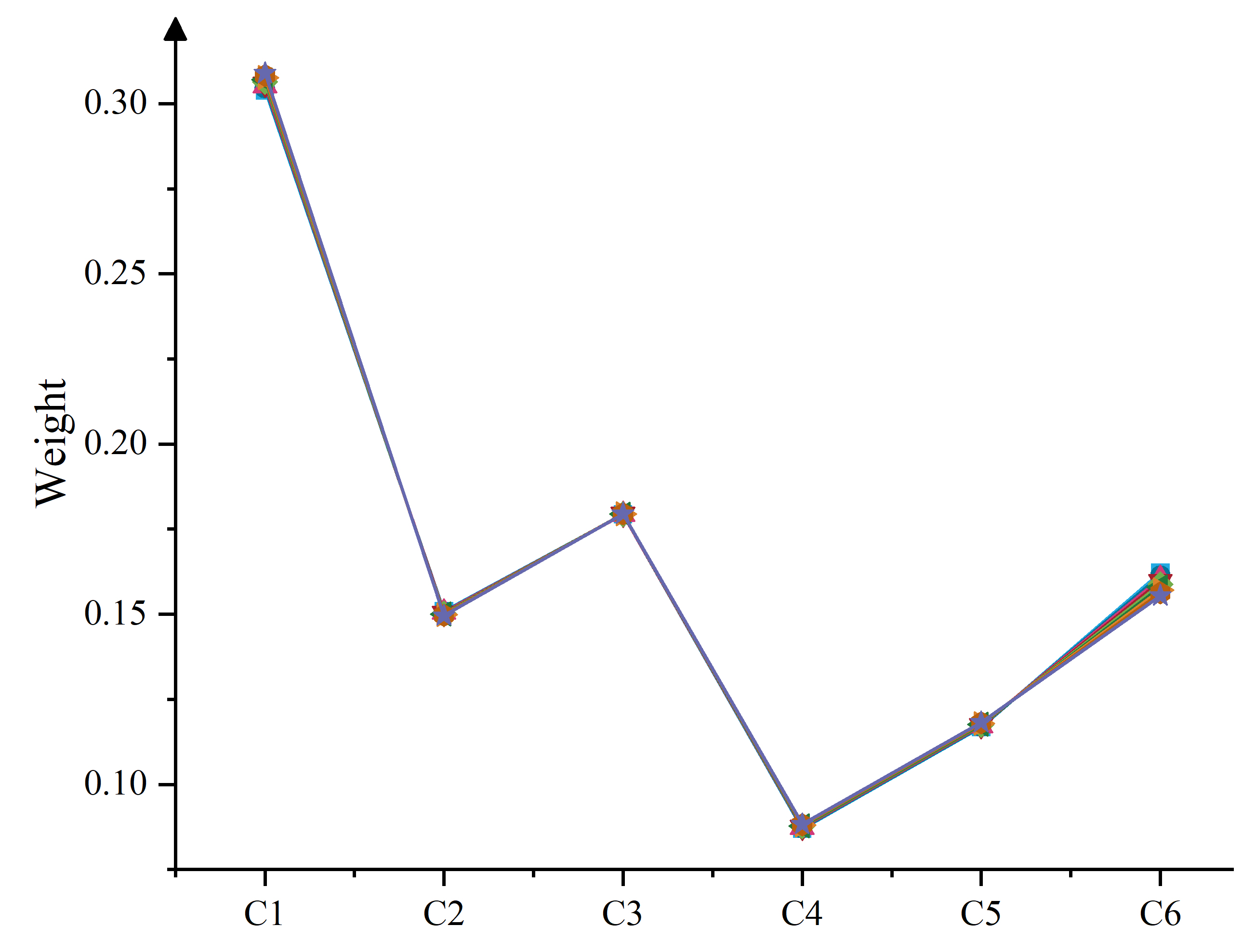}}
    \subfigure[Alternative weights of OPA-PR(N)] {
        \includegraphics[width=0.32\columnwidth]{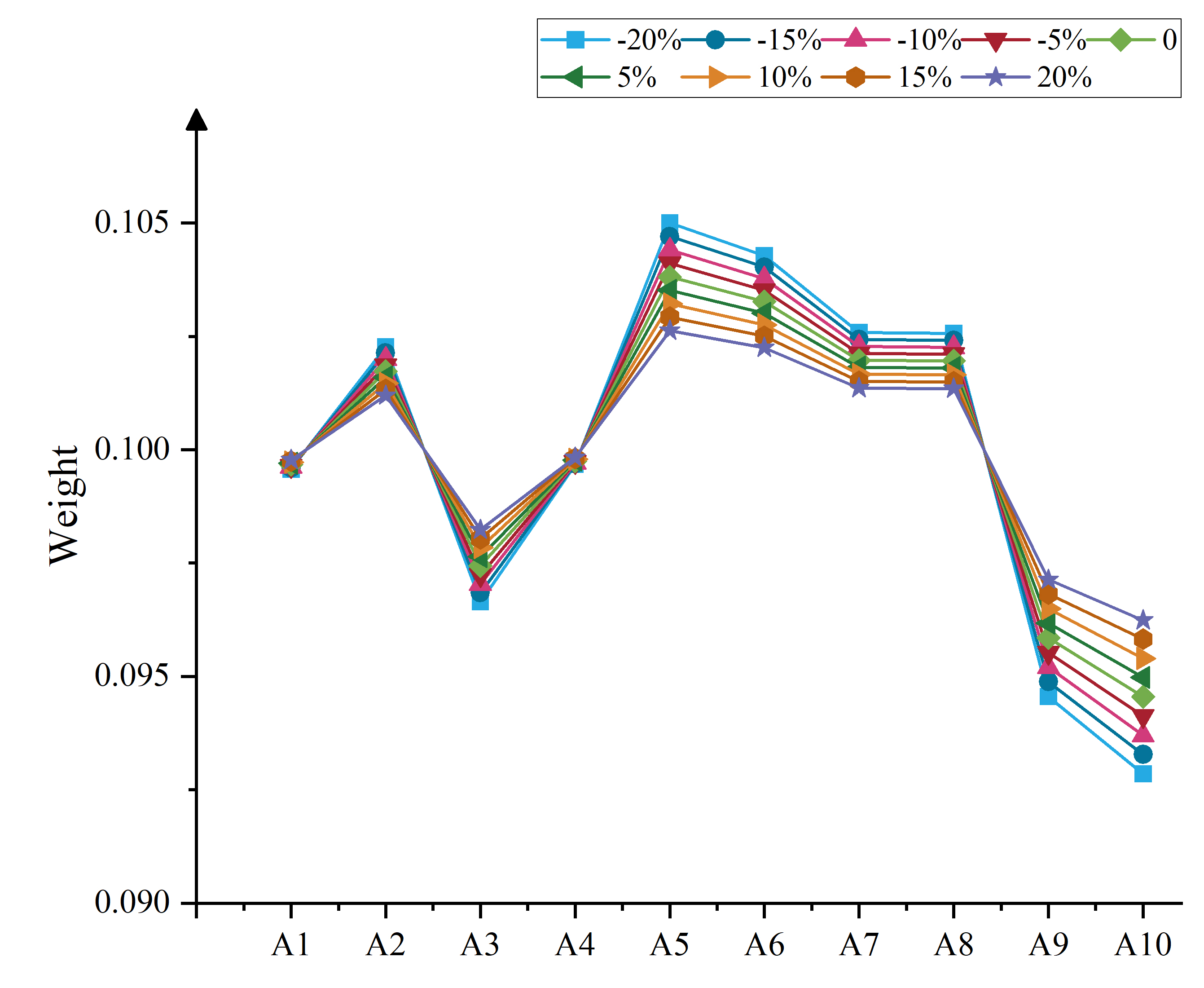}} \\
    \subfigure[Attribute weights of OPA-PR(B)] {
        \includegraphics[width=0.32\columnwidth]{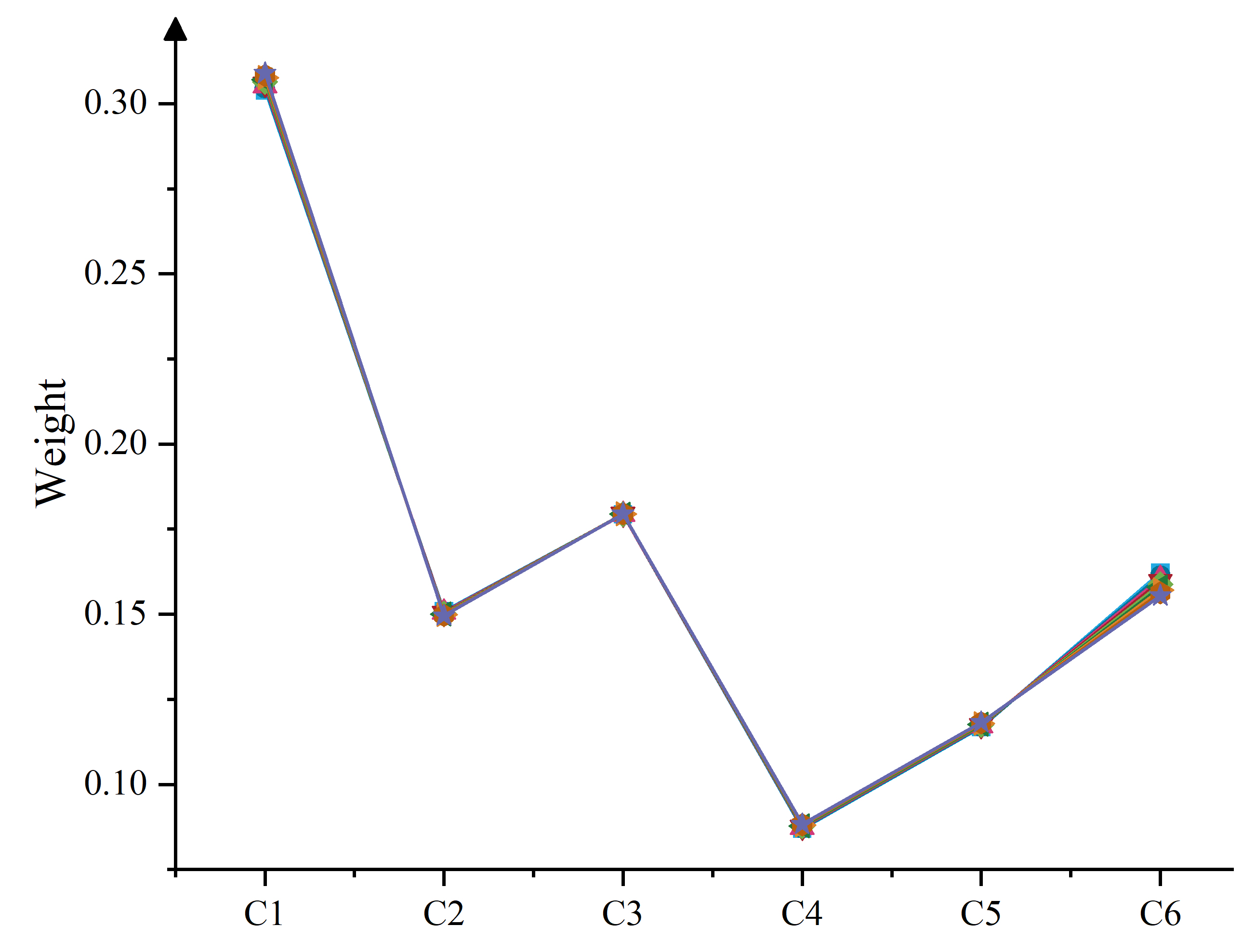}}
    \subfigure[Alternative weights of OPA-PR(B)] {
        \includegraphics[width=0.32\columnwidth]{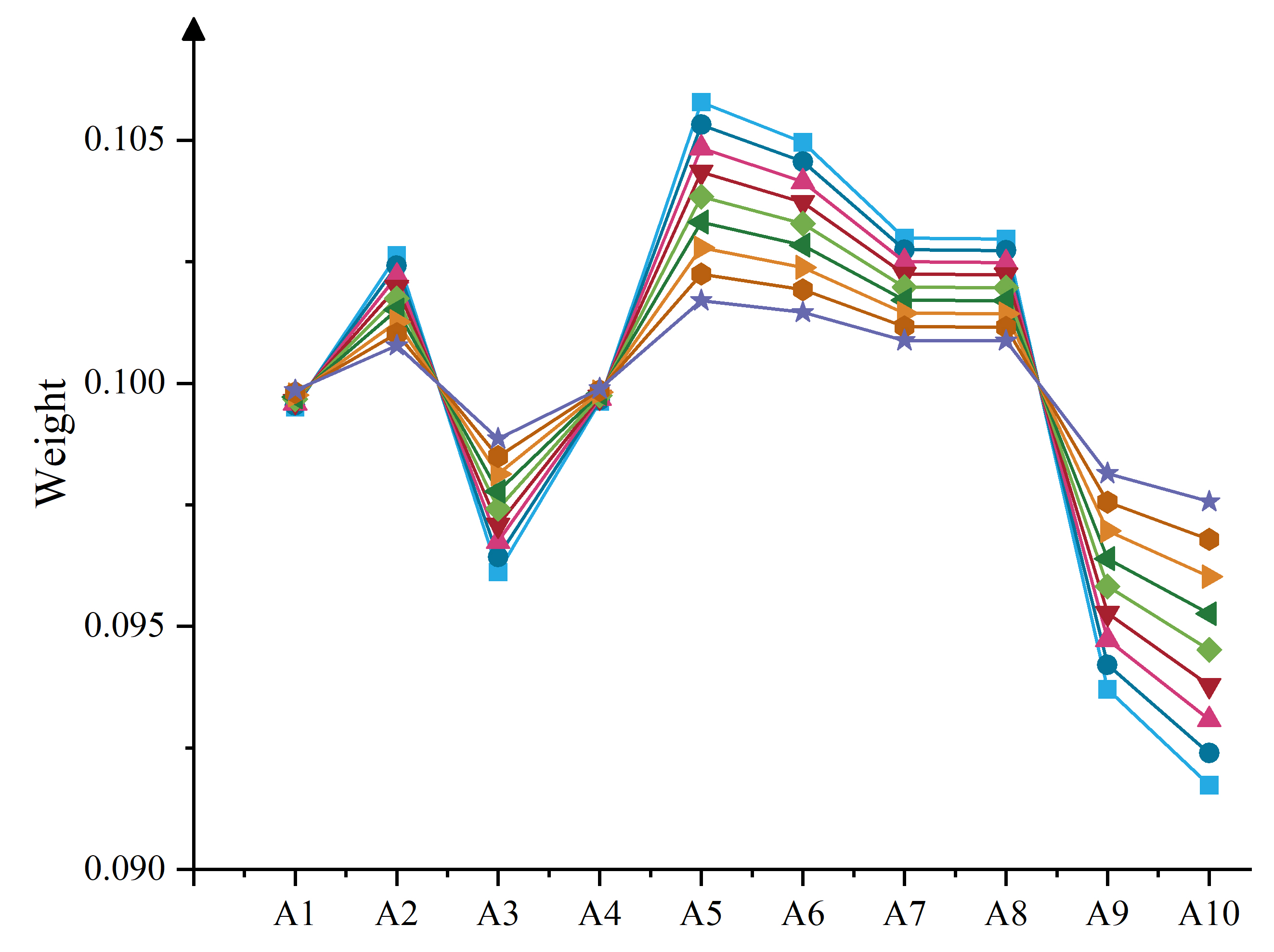}}        
    \caption{Weight results of of size parameter perturbations}
    \label{fig-06}
    \vspace{-1em}
\end{figure}

To evaluate the impact of the size parameters of attribtue ranking ambiguity sets on the weight results, we perturbate the size parameters by $\pm 5 \%$, $\pm 10 \%$, $\pm 15 \%$, and $\pm 20 \%$, with the results shown in Table \ref{tab-02} and Figure \ref{fig-06}. Table \ref{tab-02} shows that the weight disparities for both OPA-PR(N) and OPA-PR(B) decreases linearly as the size parameters of ambiguity sets increase, with OPA-PR(B) being more sensitive to size parameter perturbations than OPA-PR(N). The final ranking results for attributes and alternatives remain consistent across different perturbation scenarios. The attribute weights for both OPA-PR(N) and OPA-PR(B) are almost identical in all perturbation scenarios. However, the alternative weights for OPA-PR(N) and OPA-PR(B) reveal that as the ambiguity set size parameters decrease, variations in optimal weights for alternatives increase. An interesting phenomenon is observed, where the weights of higher-ranked alternatives increase, while those of lower-ranked alternatives decrease. This suggests that as uncertainty in attribute rankings reduces (i.e., the ambiguity set size decreases), the final weights between alternatives become more differentiated, providing more information to support less conservative decisions.

\begin{table}[h]
    \centering
    \renewcommand{\arraystretch}{1}
    \caption{Optimal weight disparities across perturbated size parameter}
    \label{tab-02}
    \begin{tabular}{@{}lllllllllll@{}}
    \toprule
    Model                      & Value & -20\%  & -15\%  & -10\%  & -5\%   & 0\%    & 5\%    & 10\%   & 15\%   & 20\%   \\ \midrule
    \multirow{2}{*}{OPA-PR(N)} & $z^{*}$     & 0.0033 & 0.0031 & 0.0029 & 0.0027 & 0.0025 & 0.0023 & 0.0021 & 0.0019 & 0.0017 \\
                               & Ratio & 131\%  & 123\%  & 116\%  & 108\%  & 100\%  & 92\%   & 84\%   & 77\%   & 69\%   \\
    \multirow{2}{*}{OPA-PR(B)} & $z^{*}$     & 0.0038 & 0.0035 & 0.0032 & 0.0029 & 0.0025 & 0.0022 & 0.0018 & 0.0015 & 0.0011 \\
                               & Ratio & 151\%  & 139\%  & 126\%  & 113\%  & 100\%  & 86\%   & 73\%   & 59\%   & 44\%   \\ \bottomrule
    \end{tabular}
\end{table}

\subsection{Comparison Analysis}
\label{section-4-4}

This section compares OPA-PR by examining its optimal weights for ranked alternatives across attributes against those of the original OPA. The OPA-PR outcomes are presented in three forms: the original OPA-PR (OPA-PR(B) and OPA-PR(N)), OPA-PR with expert-provided attribute rankings (OPA-PR(DR)), and OPA-PR with ROC weights for ranked alternatives (OPA-PR(B-ROC) and OPA-PR(N-ROC)). For the original OPA, expert-provided attribute rankings are directly used as input. The comparison with other MARS methods is omitted for two reasons: first, these methods are based on different assumptions and axioms; second, they use parameters with different meanings, making direct ranking comparisons unjustified. Comparing the proposed approach with the original model provides a clearer understanding of how its extension influences the results.

\begin{figure}[h]
    \centering
    \subfigure[C1]{
        \includegraphics[width=0.32\columnwidth]{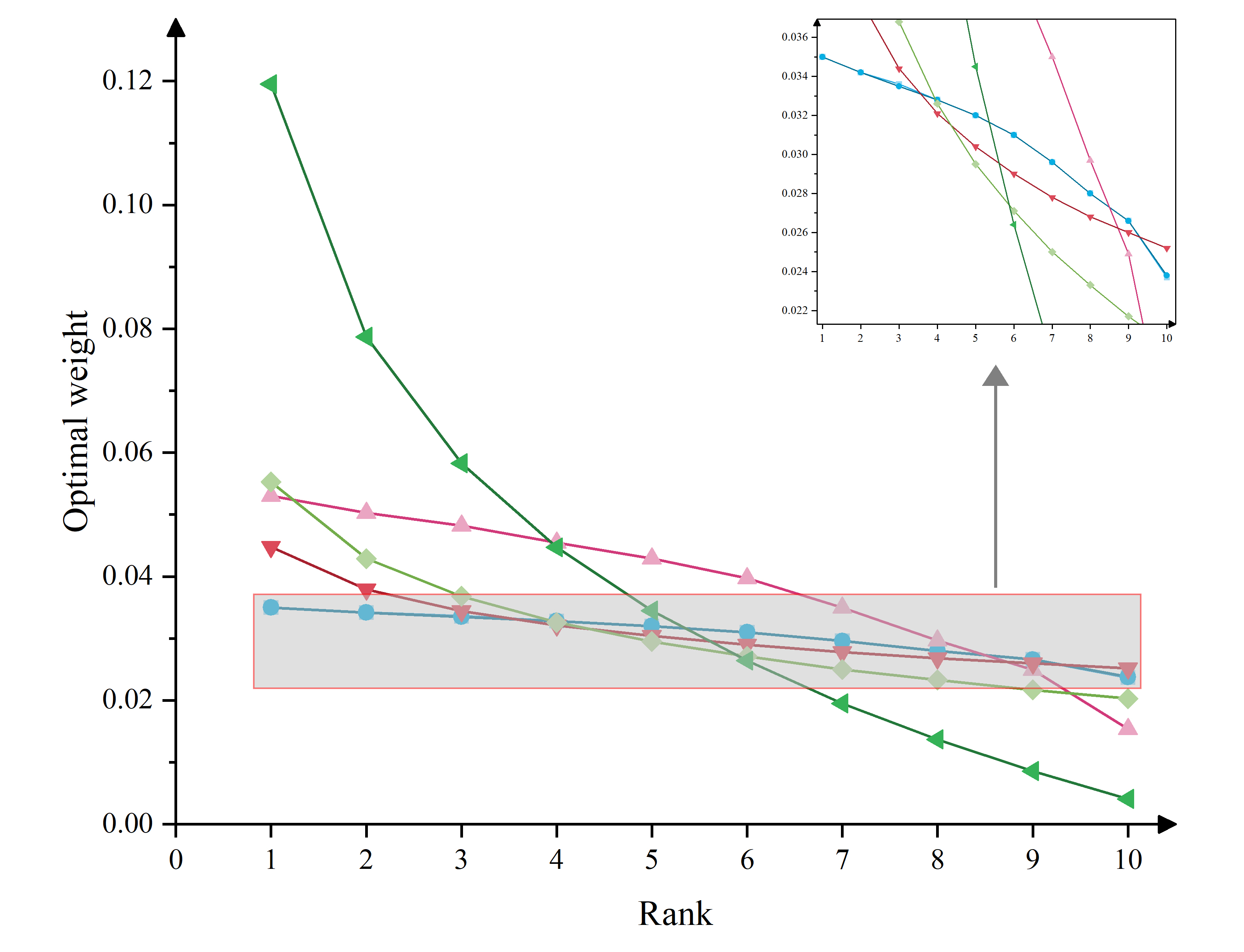}}
    \subfigure[C2]{
        \includegraphics[width=0.32\columnwidth]{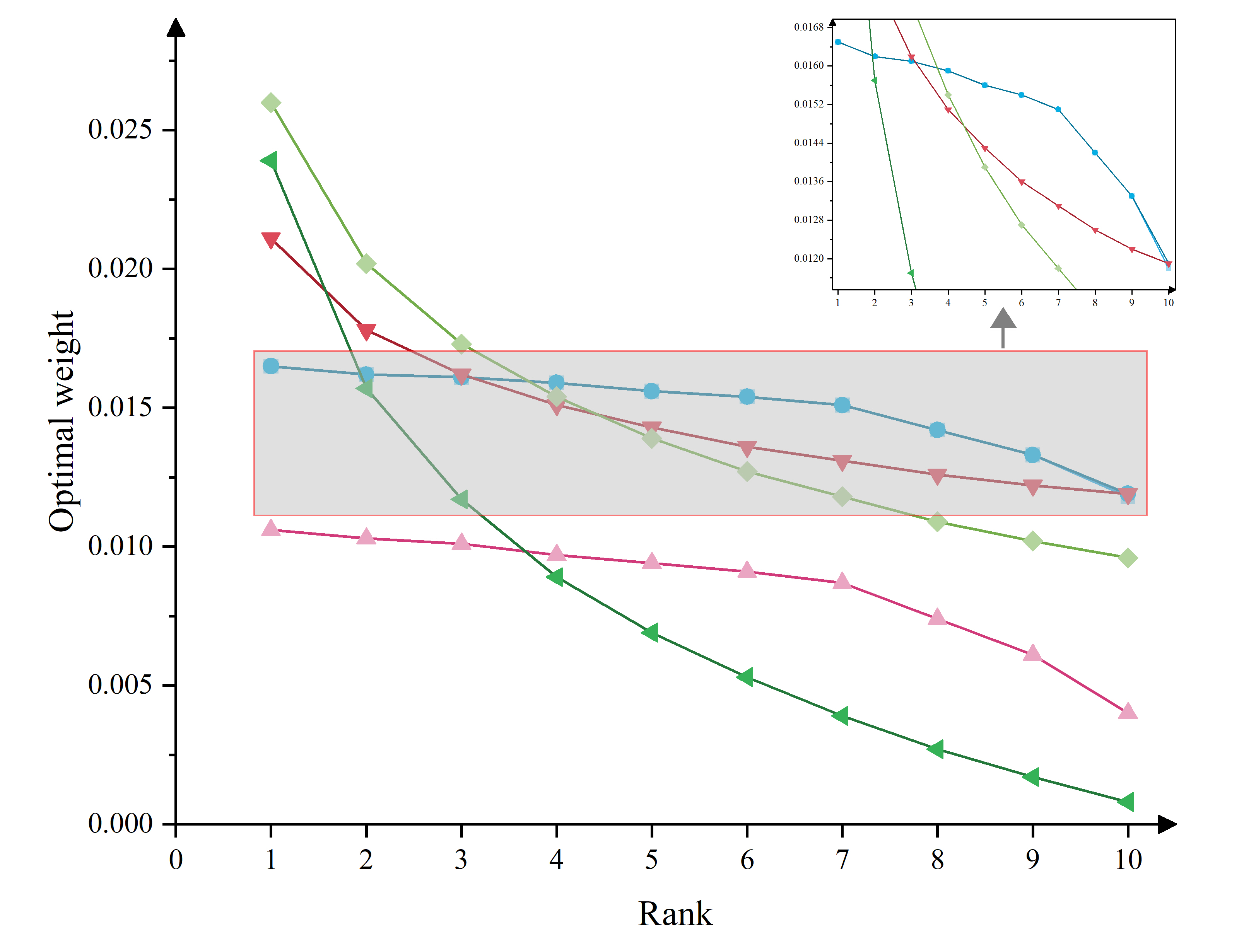}}
    \subfigure[C3]{
        \includegraphics[width=0.32\columnwidth]{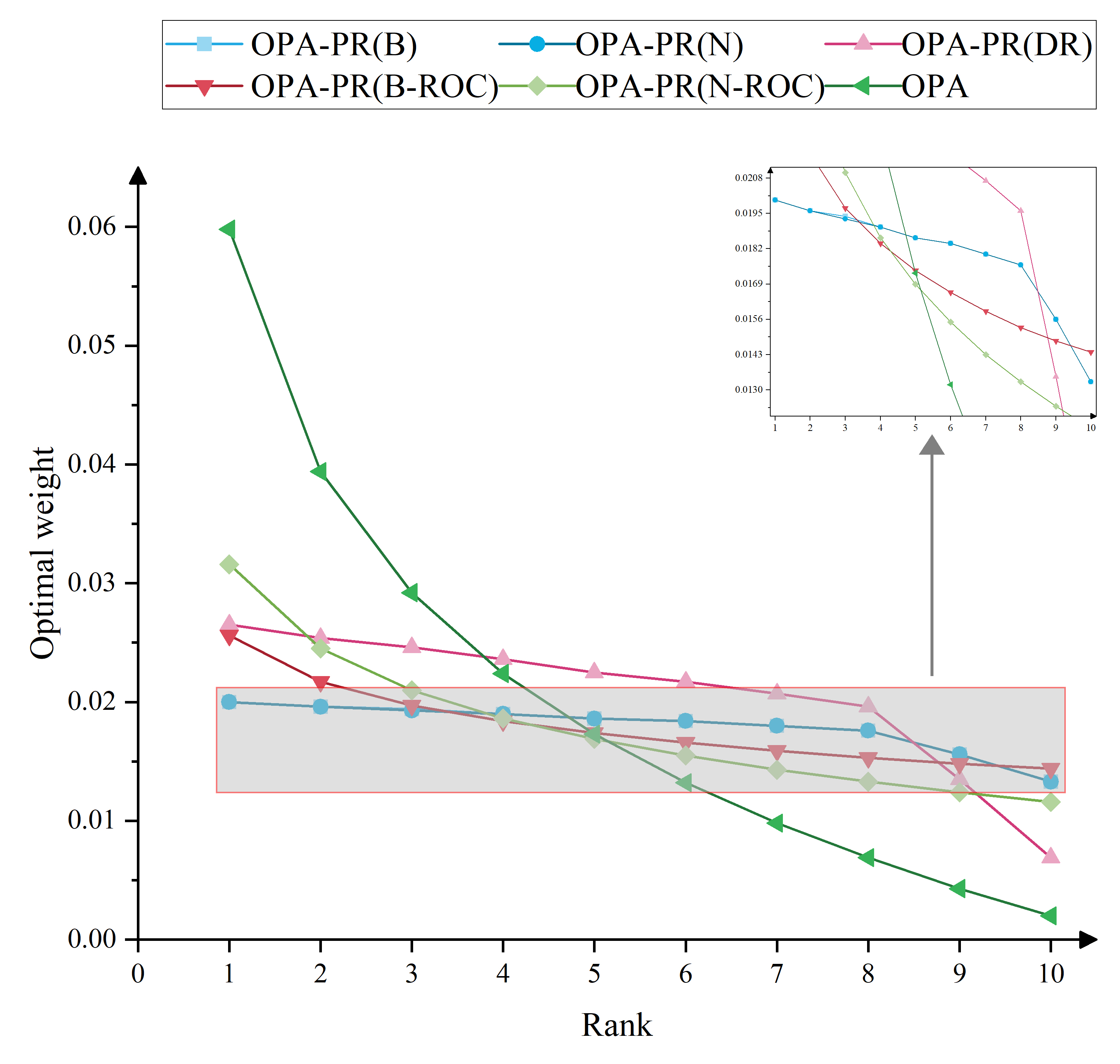}}        \\
    \subfigure[C4]{
        \includegraphics[width=0.32\columnwidth]{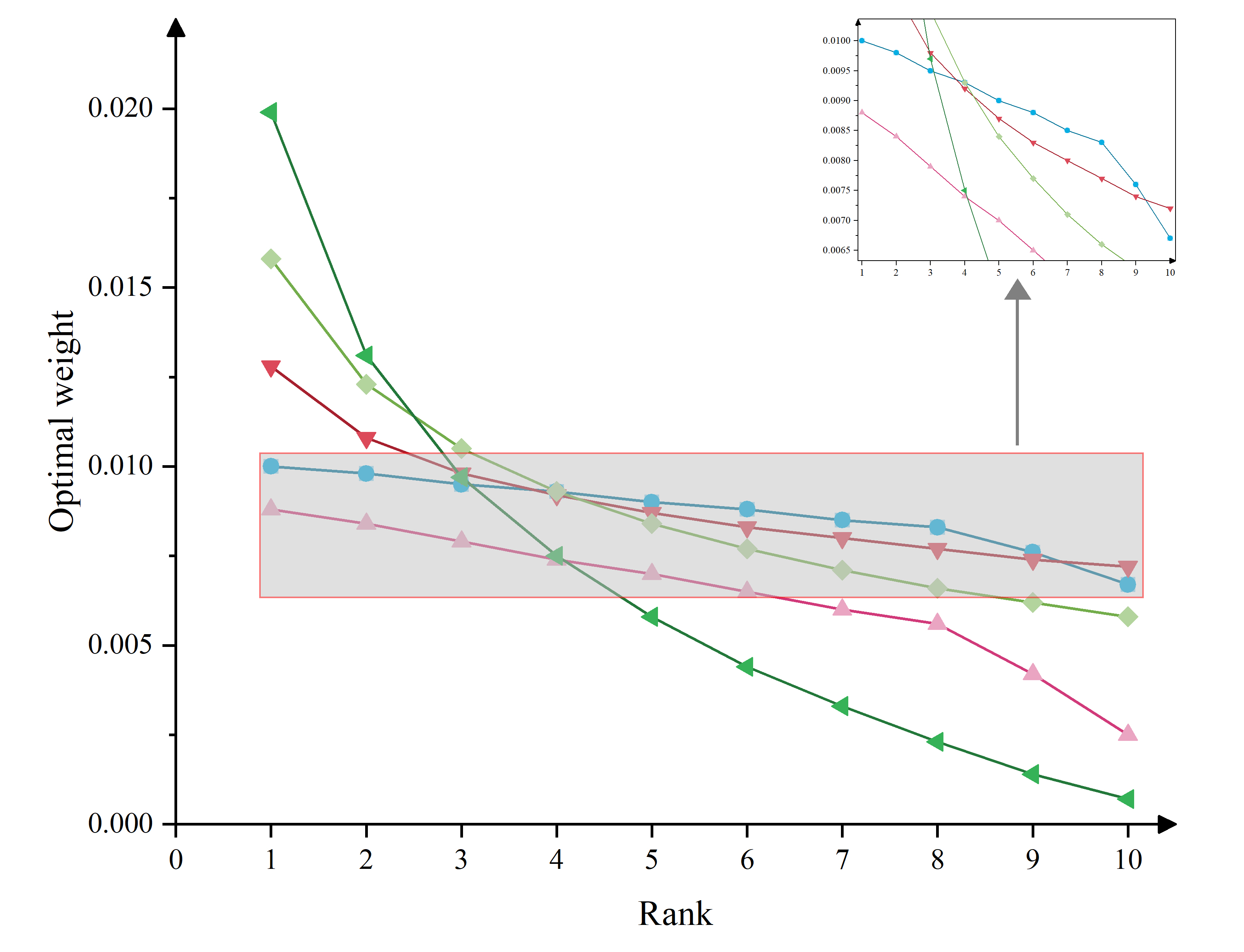}}
    \subfigure[C5]{
        \includegraphics[width=0.32\columnwidth]{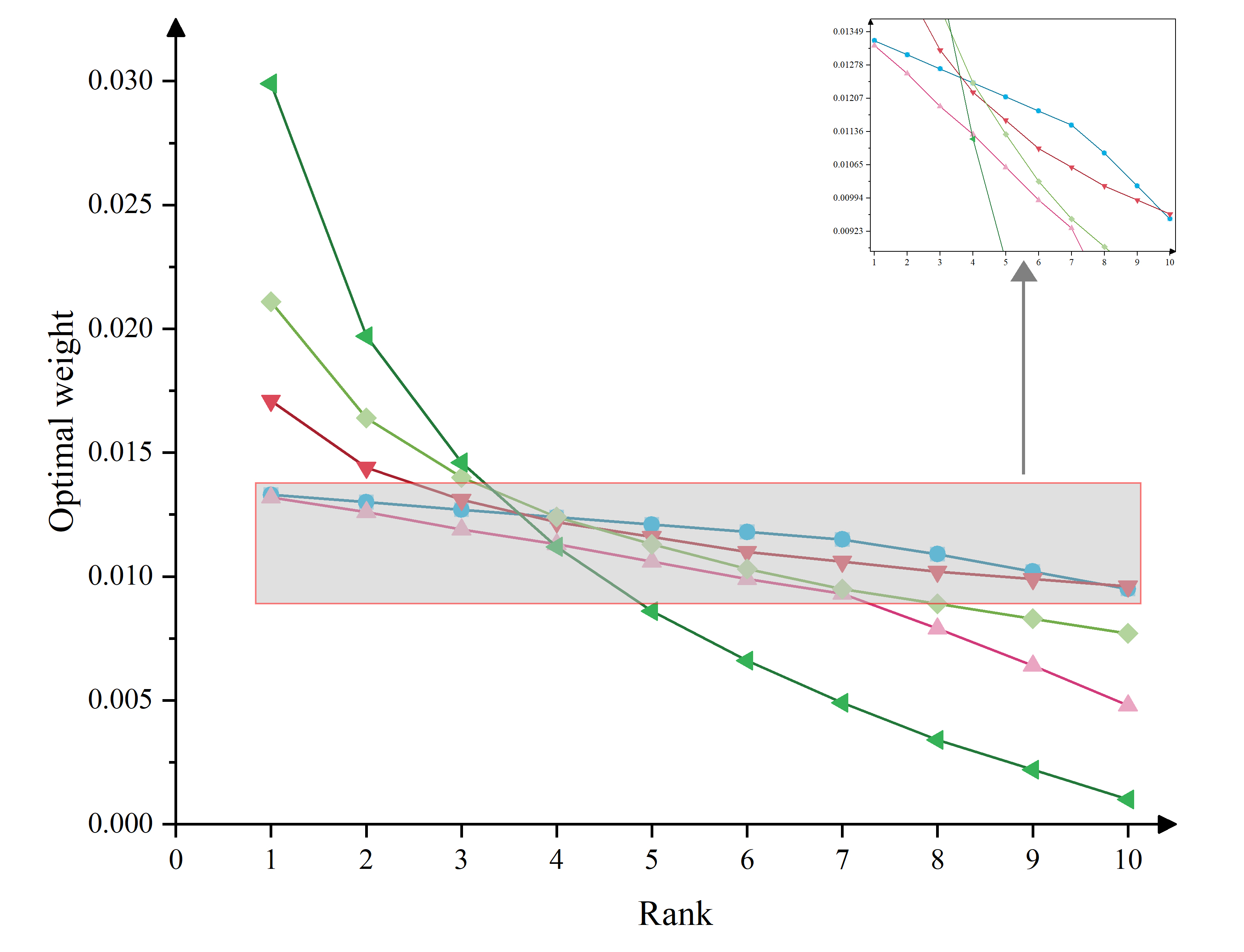}}
    \subfigure[C6]{
        \includegraphics[width=0.32\columnwidth]{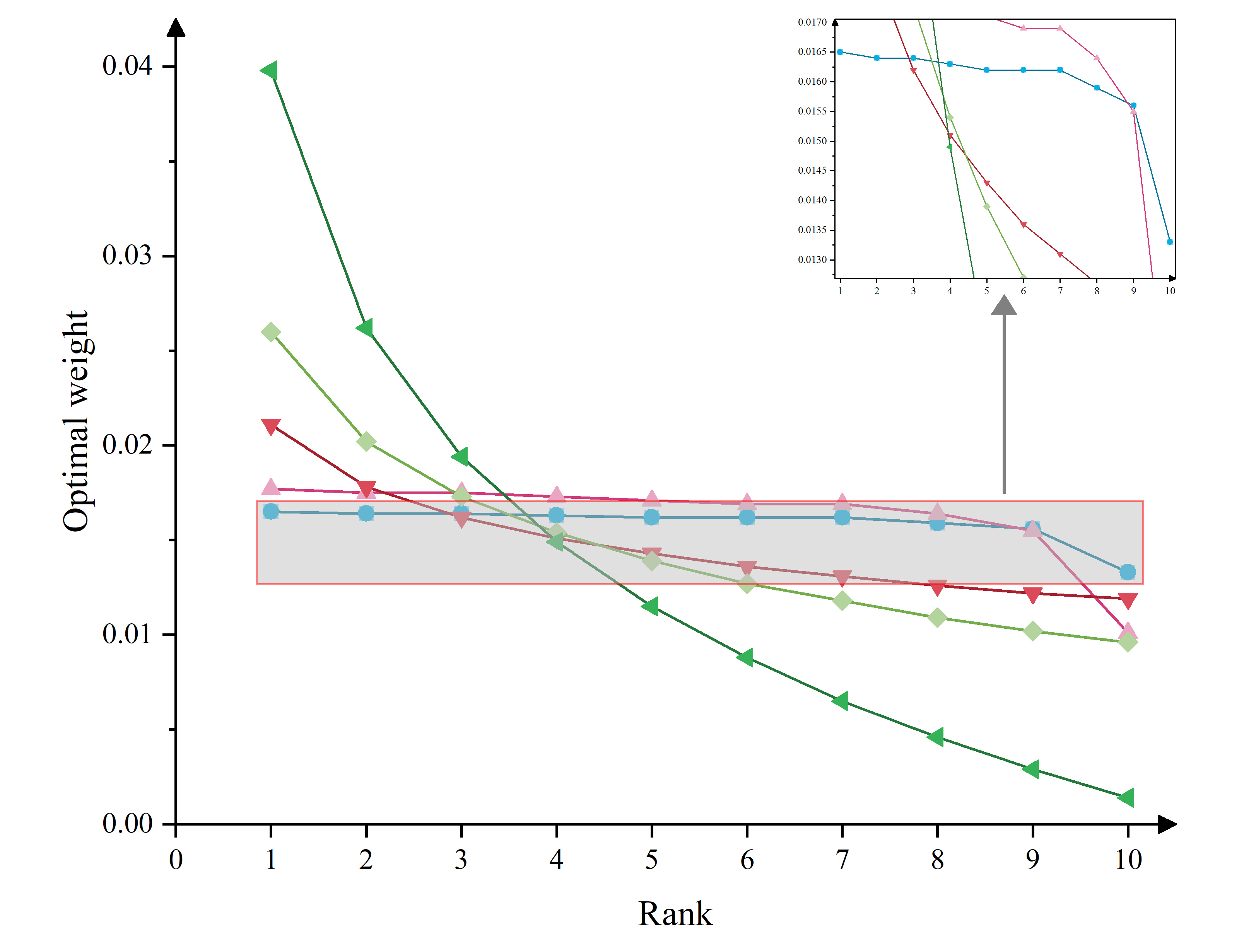}}        \\
    \caption{Optimal weights for ranked alternatives across attributes}
    \label{fig-05}
    \vspace{-1em}
\end{figure}

Figure \ref{fig-05} shows the optimal weight results for each method. The results indicate that the optimal weights for ranked alternatives in OPA span the widest range across all attributes. Additionally, the optimal weight curves for OPA-PR(B), OPA-PR(N), OPA-PR(B-ROC), and OPA-PR(N-ROC) are convex, as they use ROC weights for ranked alternatives, reflecting a risk-seeking preference. In contrast, methods based on elicited worst-case marginal utilities for ranked alternatives, such as OPA-PR(B), OPA-PR(N), and OPA-PR(DR), exhibit concave optimal weight curves, indicating risk aversion. We find that OPA-PR's robustness against preference and parametric uncertainty leads to a tendency toward balance. This suggests that when DM lacks sufficient and precise information, he/she is more likely to minimize attribute differences between alternatives, thus avoiding more aggressive decisions. We also calculated the standard deviations of the optimal weight results for each method, yielding the following outcomes: OPA (0.0207) $>$ OPA-PR(N-ROC) (0.0198) $>$ OPA-PR(DR) (0.0122) $>$ OPA-PR(B-ROC) (0.0079) $>$ OPA-PR(N) (0.0072) = OPA-PR(B) (0.0072). These results show that methods using ROC weights have higher standard deviations for optimized weights than those using elicited worst-case marginal utilities (OPA v.s. OPA-PR, OPA-PR(B-ROC) v.s. OPA-PR(B), and OPA-PR(N-ROC) v.s. OPA-PR(N)). Furthermore, under dual risk aversion toward preference and parametric uncertainty, OPA-PR(B-ROC) and OPA-PR(N) exhibit the most conservative results. Overall, even for the same problem, differences in risk preferences lead to varying outcomes, emphasizing the need to extend OPA in terms of preference and parametric uncertainty.

\section{Concluding Remarks}
\label{section-5}

By extending the traditional OPA to address parametric and preference uncertainty, the proposed OPA-PR offers a robust MARS framework. It employs an estimate-then-optimize two-stage procedure. In the first stage, it elicited the worst-case utility functions from ambiguity sets of utility preferences, incorporating properties such as monotonicity, normalization, and moment-based preference elicitation, thereby providing a more comprehensive representation of utility preferences. In the second stage, it optimizes decision weights under worst-case ranking parameters from various ambiguity sets (norm-, budget-, and CVaR-based), further strengthening the robustness of the decision-making process. We present a tractable reformulation of both stages. Notably, the PLA is utilized to approximate the utility preference ambiguity sets to derive worst-case marginal utility functions, and thus the error bounds for both stages of OPA-PR are analyzed. We demonstrate that PLA introduces no additional errors under specific preference elicitation methods (such as deterministic utility comparison and stochastic lottery comparison for moment-based preferences), which provides the foundation for designing the preference elicitation strategy. 

The effectiveness of the proposed approach is validated through numerical experiments, including case studies, sensitivity analysis, and comparison tests. The case results indicate that OPA-PR with norm- and budget-based ambiguity sets performs well, while OPA-PR with the CVaR-based ambiguity set is less effective at differentiating weights for ranked alternatives. Sensitivity analysis of alternative ranking perturbations, assuming normal distributions with the mean of expert-provided alternative rankings and standard deviations from perturbation ranges, shows that the final rankings converge to those derived from the unperturbed problem. Sensitivity analysis of the ambiguity set size parameter reveals the common trend that larger sets lead to more conservative weights. Comparison results demonstrate that the proposed approach outperforms the original OPA model in terms of robustness to parametric and preference ambiguity, resulting in more conservative weight outcomes.

While the current study focuses on a particular context of decision-making, we recognize the need for further exploration in diverse scenarios to validate the generalizability of the proposed approach. Future research could extend this work by incorporating different types of ambiguity sets for ranking parameters, expanding the applicability of OPA-PR in various domains. Additionally, a promising direction is to investigate how expert preferences can be elicited to design the size parameters for ambiguity sets more effectively. Further integration of advanced operations research techniques, such as distributionally robust optimization and stochastic contextual optimization, within the artificial intelligence framework, could offer even more robust decision-making tools. Finally, expanding OPA to accommodate model misspecifications presents a crucial avenue for future work, as this could provide valuable insights for real-world decision-making applications where expert judgements may not always hold true. 




\bibliographystyle{elsarticle-harv}
\bibliography{reference}

\newpage

\appendix

\section{Extended Formulation Considering Preference Inconsistency}
\label{Asection-01}

In this section, we discuss the potential preference inconsistency for the preference elicitation process, which may arise from factors such as misalignment with expected utility theory axioms or contaminated preference information. We examine three interpretable types of preference inconsistency, including weight disparity error and erroneous elicitation, along with their corresponding modified formulations.

First, we assume there are errors in weight disparity, indicating that the optimized weight disparity scalar does not satisfy all weight disparity constraints. 
For any $i \in \mathcal{I}$, the term $\bar{\mathbf{W}}_{i}\tilde{\bm{s}}_{i}$ is modified to $\bar{\mathbf{W}}_{i}\tilde{\bm{s}}_{i} + \bm{\gamma}_{i}$, where $\bm{\gamma}_{i} \in \mathbb{R}^{J \times R}$ denotes the error tolerance for expert $i$. An error budget $\Gamma_{i} \geq 0$ is introduced to regulate the total error, ensuring consistency within the model. The cumulative inconsistency is quantified by the total error term $\bm{e}^{\top} \bm{\gamma}_{i} \leq \Gamma$, which enforces the feasibility of the second-stage problem of OPA-PR. The modified formulation is as follows
\begin{equation}
    \begin{aligned}
        \max_{z_{i}, \bar{\bm{w}}_{i}, \bm{\gamma}_{i}} \text{ } & z_{i} \\
        \mathrm{s.t.} \text{ } & R \bm{U}_{i}^{*} z_{i} \leq \bar{\mathbf{W}}_{i}\tilde{\bm{s}}_{i} + \bm{\gamma}_{i}, && \quad \forall  \tilde{\bm{s}}_{i} \in \mathcal{V}_{i}, \\
        & \bm{e}^{\top} \bar{\bm{w}}_{i} = 1, \\
        & \bm{e}^{\top} \bm{\gamma}_{i} \leq \Gamma, \\
        & \bm{\gamma}_{i} \geq 0, \bar{\bm{w}}_{i} \geq 0,
    \end{aligned}
    \label{eq-A-1}
\end{equation}
where $\bm{\gamma}_{i}$ is treated as a decision variable. The tolerated total error can be designed based on the optimal value of $\tilde{\bm{z}}_{i}^{*}$ from the original formulation of OPA-PR. The tractable reformulation of Equation \eqref{eq-A-1} can be derived by the symmetric argument as the proofs of Propositions \ref{proposition-05}-\ref{proposition-07}. 

When addressing erroneous elicitation, we can selectively relax a portion of the preference elicitation constraints in the second-stage problem of OPA-PR. In this case, we assume that at least $1 - \theta_{ij}$ of the $R$ dominance relations under each expert and attribute are correct, implying that expert $i$ can give incorrect responses in at most $\theta_{ij}R$ relations on attribute $j$ and $\bm{\theta}_{i} \in \mathbb{R}^{J}$. We introduce binary variables $\bm{\vartheta}_{i} \in \{0,1\}^{J \times R}$ for any $i \in \mathcal{I}$, taking the value of 1 if expert $i$ is incorrect about the ranking relation $r$ on attribute $j$. The modified formulation is as follows
\begin{equation}
    \begin{aligned}
        \max_{z_{i}, \bar{\bm{w}}_{i}, \bm{\gamma}_{i}} \text{ } & z_{i} \\
        \mathrm{s.t.} \text{ } & R \bm{U}_{i}^{*} z_{i} \leq \bar{\mathbf{W}}_{i}\tilde{\bm{s}}_{i} + M \bm{\vartheta}_{i}, && \quad \forall  \tilde{\bm{s}}_{i} \in \mathcal{V}_{i}, \\
        & M (\bm{e}^{\top} - \bm{\vartheta}_{i}) + R \bm{U}_{i}^{*} z_{i} \geq \bar{\mathbf{W}}_{i}\tilde{\bm{s}}_{i}, && \quad \forall  \tilde{\bm{s}}_{i} \in \mathcal{V}_{i}, \\
        & \bm{e}^{\top} \bar{\bm{w}}_{i} = 1, \\
        & \bar{\mathbf{E}} \bm{\vartheta}_{i} \leq R \bm{\theta}_{i}, \\
        & \bm{\vartheta}_{i} \in \{0,1\}^{J \times R}, \bar{\bm{w}}_{i} \geq 0,
    \end{aligned}
    \label{eq-A-2}
\end{equation}
where $M$ is a large constant and 
\[
    \bar{\mathbf{E}} = \left[
        \begin{array}
        {ccc} \bar{\bm{e}}^{\top} & & \\
         & \ddots & \\
         & & \bar{\bm{e}}^{\top}
        \end{array}\right],
\]
and $\bar{\bm{e}}$ is an $R$-dimensional vector of all ones. The tractable reformulation of Equation \eqref{eq-A-2} can be derived by the symmetric argument as the proofs of Propositions \ref{proposition-05}-\ref{proposition-07} and then solved by the cutting-plane method.

\section{Technical Proofs}

\subsection{Proof of Proposition \ref{proposition-03}}

For any $(i,j) \in \mathcal{I} \times \mathcal{J}$, given that $\psi_{l_{ij}}(\tau)$ on $\tau \in \Theta$ is a step function with jumps at $\tau_{h_{ij}}$ for all $h_{ij} \in \mathcal{H}_{ij}$, we have
\[
    \int_{\Theta} \psi_{l_{ij}}(\tau)\mathrm{d}u_{ij}(\tau) = \sum_{h_{ij} \in \mathcal{H}_{ij}} \psi_{l_{ij}}(\tau_{h_{ij}}^\prime) (u_{ij}(\tau_{h_{ij}}) - u_{ij}(\tau_{h_{ij}-1})) = \int_{\Theta} \psi_{l_{ij}}(\tau)\mathrm{d}\tilde{u}_{ij}(\tau),
\]
where $\tau_{h_{ij}}^\prime \in (\tau_{h_{ij}-1}, \tau_{h_{ij}}]$ and the last equality follows from the fact that the integral only involves the value of $u_{ij}$ at these jumps. \qed

\subsection{Proof of Proposition \ref{proposition-04}}

We first divide the set of marginal utility functions based on their values at rankings on $\Theta$ for any $(i,j) \in \mathcal{I} \times \mathcal{J}$, denoted as $\mathcal{U}_{ij}(\bm{y}) := \{u_{ij}: u_{ij}(\tau_{h_{ij}}) = y_{h_{ij}}, \forall h_{ij} \in \mathcal{H}_{ij}\}$.
Then, we have $\mathcal{U}_{ij}(\bm{y}) \cap \mathcal{U}^{\text{conc}} \ne \emptyset$, which yields Equations \eqref{eq-12-02} and \eqref{eq-12-03}.
Equations \eqref{eq-12-02} and \eqref{eq-12-03} express the concavity property by the first order condition, $y_{h_{ij}}-y_{h_{ij}-1}=\mu_{h_{ij}-1}(\tau_{h_{ij}} - \tau_{h_{ij}-1})$ and $\mu_{h_{ij}}\leq \mu_{h_{ij}-1}$ for all $h_{ij} \in \mathcal{H}_{ij}$, where $\mu_{h_{ij}}$ represents the subgradient at turning point $h_{ij}$.
The left-hand side of Equation \eqref{eq-12-04} characterizes monotonic increasing properties of marginal utility function, while the right-hand side represents the Lipschitz continuity with the modulus being bounded by $G$.
Furthermore, $\mathcal{U}_{ij}(\bm{y}) \cap \mathcal{U}^{\text{nor}} \ne \emptyset$, which are represented by Equation \eqref{eq-12-07}.
Equation \eqref{eq-12-05} details the elicited moment-type preference information for each ranking points, which is derived from the fact that $\mathrm{d}\tilde{u}_{ij}(\tau_{h_{ij}-1}) = \mu_{h_{ij}-1}$ for all $h_{ij} \in \mathcal{H}_{ij}$.
Since $\tilde{u}_{ij}(h(\bm{x},\bm{\xi}_{e}))$ in the objective function is concave, non-decreasing, and affine in $\bm{x}$, we can apply the support function for increasing concave functions, as outlined in Lemma \ref{lemma-01}, to approximate $\tilde{u}_{ij}(h(\bm{x},\bm{\xi}_{e}))$. This leads to the objective function formulation in Equation \eqref{eq-12-01} and the constraints in Equations \eqref{eq-12-06} and \eqref{eq-12-08}. The worst-case marginal utility function in Equation \eqref{eq-13} for any $(i,j) \in \mathcal{I} \times \mathcal{J}$ is directly obtained by Proposition \ref{proposition-03}. \qed

\subsection{Proof of Proposition \ref{proposition-05}}
Consider the inequality constraint in Equation \eqref{eq-14} under Assumption \ref{assumption-08}, i.e., 
\[
    R \tilde{\bm{U}}_{i}^{*} \tilde{z}_{i} \leq \tilde{\mathbf{W}}_{i}\tilde{\bm{s}}_{i}, \forall \tilde{\bm{s}}_{i} \in \mathcal{V}_{i}^{\text{norm}} \Leftrightarrow R \tilde{\bm{U}}_{i}^{*} \tilde{z}_{i} \leq \min_{\tilde{\bm{s}}_{i} \in \mathcal{V}_{i}^{\text{norm}}} \left\{\tilde{\mathbf{W}}_{i}\tilde{\bm{s}}_{i}\right\}.
\]

Given any $\tilde{\bm{w}}_{i}$, we have
\[
\begin{aligned}
    & \min_{\tilde{\bm{s}}_{i} \in \mathcal{V}_{i}^{\text{norm}}} \tilde{\mathbf{W}}_{i}\tilde{\bm{s}}_{i}, \\
    = & \min_{\|\bm{\Sigma}^{-\frac{1}{2}}(\tilde{\bm{s}}_{i} - \bm{\mu}) \|_{2} \leq \delta_{i}} \tilde{\mathbf{W}}_{i}\tilde{\bm{s}}_{i}, \\
    = & \min_{\|\bm{\Sigma}^{-\frac{1}{2}}(\tilde{\bm{s}}_{i} - \bm{\mu}) / \delta_{i}\|_{2} \leq 1} \tilde{\mathbf{W}}_{i}\tilde{\bm{s}}_{i}.
\end{aligned} 
\] 

Let $\bm{\xi}_{i} = \bm{\Sigma}^{-\frac{1}{2}}(\tilde{\bm{s}}_{i} - \bm{\mu}) / \delta_{i}$, so $\tilde{\bm{s}}_{i} = \bm{\mu} + \delta_{i} \bm{\Sigma}^{\frac{1}{2}}\bm{\xi}_{i}$, which leads to
\[
\begin{aligned}
    & \min_{\|\bm{\xi}_{i}\|_{2} \leq 1} \tilde{\mathbf{W}}_{i}(\bm{\mu} + \delta_{i} \bm{\Sigma}^{\frac{1}{2}}\bm{\xi}_{i}), \\
    = & - \max_{\|\bm{\xi}_{i}\|_{2} \leq 1} \tilde{\mathbf{W}}_{i}(-\bm{\mu} - \delta_{i} \bm{\Sigma}^{\frac{1}{2}}\bm{\xi}_{i}), \\
    = & \tilde{\mathbf{W}}_{i}\bm{\mu} - \max_{\|\bm{\xi}_{i}\|_{2} \leq 1} \tilde{\mathbf{W}}_{i}(- \delta_{i} \bm{\Sigma}^{\frac{1}{2}}\bm{\xi}_{i}), \\
    = & (\tilde{\mathbf{W}}_{i})_{g}^{\top}\bm{\mu} - \delta_{i} \|(\bm{\Sigma}^{\frac{1}{2}})^{\top}\tilde{\bm{w}}_{i}\|_{2}^{*}, \\
    = & (\tilde{\mathbf{W}}_{i})_{g}^{\top}\bm{\mu} - \delta_{i} \|(\bm{\Sigma}^{\frac{1}{2}})^{\top}\tilde{\bm{w}}_{i}\|_{2},
\end{aligned}
\]
where the fourth equality follows from the dual norm definition, and the fifth equality is due to the fact that the dual norm of the $\ell_2$-norm is the $\ell_2$-norm itself. By reintegrating the inequality constraint, Equation \eqref{eq-14} under Assumption \ref{assumption-06} is equivalent to the following second order cone programming problem
\[
    \begin{aligned}
        \max_{\tilde{\bm{w}}_{i} \in \bar{\mathcal{W}}_{i}, \tilde{z}_{i}, \bm{\lambda}_{i}} \text{ } &  \tilde{z}_{i}, \\
        \mathrm{s.t.} \text{ } &  (\|(\bm{\Sigma}^{\frac{1}{2}})^{\top}\tilde{\bm{w}}_{i}\|_{2})\bm{e} \leq \frac{\tilde{\mathbf{W}}_{i}\bm{\mu} - R \tilde{\bm{U}}_{i}^{*} \tilde{z}_{i}}{\delta_{i}},
    \end{aligned}
\]
which gives the result in Proposition \ref{proposition-05}. \qed

\subsection{Proof of Proposition \ref{proposition-06}}

\begin{lemma}[Decompostion of support function]
    Let $\mathcal{S}_{1}, \mathcal{S}_{2}, \dots, \mathcal{S}_{n}$ be closed convex sets, such that $\cap_{i\in [n]} \mathrm{ri}(\mathcal{S}_{i}) \ne \emptyset$ and $\mathcal{S} = \cap_{i\in [n]} \mathcal{S}_{i}$. Then, we have
    \[
        \delta^*(x | S)=\min\left\{\sum_{i=1}^n\delta^*(y_i\mid S_i)\mid\sum_{i=1}^ny_i=x\right\},
    \]
    where $\delta^*(x | S)$ is the support function of $\mathcal{S}$.
    \label{lemma-03}
\end{lemma}

By Lemma \ref{lemma-03}, we are ready to derive the tractable reformulation of the second-stage problem of OPA-PR under Assumption \ref{assumption-07}. Consider the inequality constraint in Equation \eqref{eq-14} under Assumption \ref{assumption-07}, i.e., 
\[
    R \tilde{\bm{U}}_{i}^{*} \tilde{z}_{i} \leq \tilde{\mathbf{W}}_{i}\tilde{\bm{s}}_{i}, \forall \tilde{\bm{s}}_{i} \in \mathcal{V}_{i}^{\text{budget}} \Leftrightarrow R \tilde{\bm{U}}_{i}^{*} \tilde{z}_{i} \leq \min_{\tilde{\bm{s}}_{i} \in \mathcal{V}_{i}^{\text{budget}}} \tilde{\mathbf{W}}_{i}\tilde{\bm{s}}_{i} .
\]

Let 
\[
    \bm{\mu} = (\underbrace{\mu_{1}, \dots, \mu_{1}}_{R \text{ elements}} , \dots, \underbrace{\mu_{J}, \dots, \mu_{J}}_{R \text{ elements}})^{\top} \quad \text{and} \quad  \bm{\Delta}_{i} = \mathrm{diag}\Big(\underbrace{\frac{1}{\gamma_{i1}}, \dots, \frac{1}{\gamma_{i1}}}_{R \text{ elements}} , \dots, \underbrace{\frac{1}{\gamma_{iJ}}, \dots \frac{1}{\gamma_{iJ}}}_{R \text{ elements}}\Big).
\]

We can rewrite $\mathcal{V}_{i}^{\text{budget}}$ into
\[
    \mathcal{V}_{i}^{\text{budget}} := \left\{ \tilde{\bm{s}}_{i} \in \mathbb{R}^{JR}: \| \bm{\Delta}_{i}(\tilde{\bm{s}}_{i} - \bm{\mu}) \|_{\infty} \leq 1, \| \bm{\Delta}_{i}(\tilde{\bm{s}}_{i} - \bm{\mu}) \|_{1} \leq \Gamma \right\} := \mathcal{V}_{i}^{\infty} \cap \mathcal{V}_{i}^{1},
\]
where 
\[
\begin{aligned}
    & \mathcal{V}_{i}^{\infty} := \left\{ \tilde{\bm{s}}_{i} \in \mathbb{R}^{JR}: \| \bm{\Delta}_{i}(\tilde{\bm{s}}_{i} - \bm{\mu}) \|_{\infty} \leq 1 \right\}, \\
    & \mathcal{V}_{i}^{1} := \left\{ \tilde{\bm{s}}_{i} \in \mathbb{R}^{JR}: \| \bm{\Delta}_{i}(\tilde{\bm{s}}_{i} - \bm{\mu}) \|_{1} \leq \Gamma \right\}.
\end{aligned}
\]

Notice that $\mathcal{V}_{i}^{\infty}$ and $\mathcal{V}_{i}^{1}$ are closed convex sets and $\bm{0} \in \mathrm{ri}(\mathcal{V}_{i}^{\infty}) \cap \mathrm{ri} (\mathcal{V}_{i}^{1})$, then Lemma \ref{lemma-03} holds. Then, we have
\[
    \max_{\tilde{\bm{s}}_{i} \in \mathcal{V}_{i}^{\text{budget}}} - \tilde{\mathbf{W}}_{i}\tilde{\bm{s}}_{i} \leq - R \tilde{\bm{U}}_{i}^{*} \tilde{z}_{i} \Leftrightarrow \delta^*(- \tilde{\mathbf{W}}_{i} | \mathcal{V}_{i}^{\text{budget}}) \leq - R \tilde{\bm{U}}_{i}^{*} \tilde{z}_{i},
\]
where $\delta^*(- (\tilde{\mathbf{W}}_{i})_{g} | \mathcal{V}_{i}^{\text{budget}})$ is the support function of $\mathcal{V}_{i}^{\text{budget}}$. By Lemma \ref{lemma-03}, we have
\[
    \min \left\{ \delta^*(\bm{y}_{1} | \mathcal{V}_{i}^{\infty}) + \delta^*(\bm{y}_{2} | \mathcal{V}_{i}^{1}): \mathbf{Y}_{1} + \mathbf{Y}_{2} = - \tilde{\mathbf{W}}_{i} \right\} \leq - R (\tilde{\bm{U}}_{i}^{*})_{g} \tilde{z}_{i}, \quad \forall g \in \mathcal{G},
\]
where $\mathbf{Y}_{1} = \mathrm{diag}(y_{111}, \dots, y_{11R}, y_{121}, \dots, y_{1JR})$ and $\mathbf{Y}_{2} = \mathrm{diag}(y_{211}, \dots, y_{21R}, y_{221}, \dots, y_{2JR})$.

By the definition of support function and dual norm, we have
\[
\begin{aligned}
    \delta^*(\bm{y}_{1} | \mathcal{V}_{i}^{\infty}) & = \max_{\tilde{\bm{s}}_{i} \in \mathcal{V}_{i}^{\infty}} \bm{y}_{1}^{\top} \tilde{\bm{s}}_{i} = \max_{\| \bm{\Delta}_{i}(\tilde{\bm{s}}_{i} - \bm{\mu}) \|_{\infty} \leq \Gamma} \bm{y}_{1}^{\top} \tilde{\bm{s}}_{i} = \max_{\|\bm{d}_{1}\|_{\infty} \leq 1} \bm{y}_{1}^{\top} (\bm{\Delta}_{i}^{-1}\bm{d}_{1} + \bm{\mu}), \\
    & = \bm{y}_{1}^{\top}\bm{\mu} + \max_{\|\bm{d}_{1}\|_{\infty} \leq 1} \bm{y}_{1}^{\top} (\bm{\Delta}_{i}^{-1}\bm{d}_{1}) = \bm{y}_{1}^{\top}\bm{\mu} + \|\bm{\Delta}_{i}^{-1}\bm{y}_{1}\|_{\infty}^{*}, \\
    & = \bm{y}_{1}^{\top}\bm{\mu} + \|\bm{\Delta}_{i}^{-1}\bm{y}_{1}\|_{1},
\end{aligned}
\]
and 
\[
\begin{aligned}
    \delta^*(\bm{y}_{2} | \mathcal{V}_{i}^{1}) & = \max_{\tilde{\bm{s}}_{i} \in \mathcal{V}_{i}^{1}} \bm{y}_{2}^{\top} \tilde{\bm{s}}_{i} = \max_{\| \bm{\Delta}_{i}(\tilde{\bm{s}}_{i} - \bm{\mu}) \|_{1} \leq \Gamma} \bm{y}_{2}^{\top} \tilde{\bm{s}}_{i} = \max_{\|\bm{d}_{2}\|_{1} \leq 1} \bm{y}_{2}^{\top} (\Gamma \bm{\Delta}_{i}^{-1}\bm{d}_{2} + \bm{\mu}), \\
    & = \bm{y}_{2}^{\top}\bm{\mu} + \max_{\|\bm{d}_{2}\|_{1} \leq 1} \bm{y}_{2}^{\top} (\Gamma \bm{\Delta}_{i}^{-1}\bm{d}_{2}) = \bm{y}_{2}^{\top}\bm{\mu} + \Gamma \|\bm{\Delta}_{i}^{-1}\bm{y}_{2}\|_{1}^{*}, \\
    & = \bm{y}_{2}^{\top}\bm{\mu} + \Gamma \|\bm{\Delta}_{i}^{-1}\bm{y}_{2}\|_{\infty}.
\end{aligned}
\]

Notice that 
\[
    \left.\min_{\bm{y}_1}\|\bm{\Delta}_{i}^{-1}\bm{y}_{1}\|_1=\min_{\bm{y}_1}\sum_{g \in \mathcal{G}}\left|(\bm{\Delta}_{i}^{-1}\bm{y}_{1})_g\right|\Leftrightarrow\min_{\bm{y}_1,\lambda}\left\{\sum_{g \in \mathcal{G}}\lambda_g\left|
\begin{array}{l}
    \bm{\lambda}\geq-\bm{\Delta}_{i}^{-1}\bm{y}_{1} \\
    \bm{\lambda}\geq\bm{\Delta}_{i}^{-1}\bm{y}_{1}
\end{array}\right.\right.\right\},
\]
and
\[
    \left.\min_{\bm{y}_2}\|\bm{\Delta}_i^{-1}\bm{y}_2\|_\infty=\min_{\bm{y}_2}\max_{g \in \mathcal{G}}\left|(\bm{\Delta}_i^{-1}\bm{y}_2)_g\right|\Leftrightarrow\min_{\bm{y}_2,\nu}\left\{\nu\left|
\begin{array}{l}
    \nu \bm{e}\geq-\bm{\Delta}_i^{-1}\bm{y}_2 \\
    \nu \bm{e}\geq\bm{\Delta}_i^{-1}\bm{y}_2
\end{array}\right.\right.\right\}.
\]

By reintegrating the inequality constraint, Equation \eqref{eq-14} under Assumption \ref{assumption-07} is equivalent to the following linear programming problem
\[
    \begin{aligned}
        \max_{\tilde{\bm{w}}_{i} \in \bar{\mathcal{W}}_{i}, \tilde{z}_{i}, \bm{y}_{1}, \bm{y}_{2}, \bm{\lambda}, \nu} \text{ } &  \tilde{z}_{i}, \\
        \mathrm{s.t.} \text{ } & (\bm{y}_{1} + \bm{y}_{2})^{\top} \bm{\mu} + (\bm{\lambda}^{\top}\bm{e} + \Gamma \nu) \bm{e} \leq - R \tilde{\bm{U}}_{i}^{*} \tilde{z}_{i}, \\
        & \mathbf{Y}_{1} + \mathbf{Y}_{2} = -\tilde{\mathbf{W}}_{i}, \\
        & \bm{\lambda}\geq-\bm{\Delta}_{i}^{-1}\bm{y}_{1}, \\
        & \bm{\lambda}\geq\bm{\Delta}_{i}^{-1}\bm{y}_{1}, \\
        & \nu \bm{e}\geq-\bm{\Delta}_i^{-1}\bm{y}_2, \\
        & \nu \bm{e}\geq\bm{\Delta}_i^{-1}\bm{y}_2,
    \end{aligned}
\]
which gives the result in Proposition \ref{proposition-06}. \qed

\subsection{Proof of Proposition \ref{proposition-07}}

Consider the inequality constraint in Equation \eqref{eq-14} under Assumption \ref{assumption-08}, i.e., 
\[
    R \tilde{\bm{U}}_{i}^{*} \tilde{z}_{i} \leq \tilde{\mathbf{W}}_{i}\tilde{\bm{s}}_{i}, \forall \tilde{\bm{s}}_{i} \in \mathcal{V}_{i}^{\text{CVaR}}, \Leftrightarrow R \tilde{\bm{U}}_{i}^{*} \tilde{z}_{i} \leq \min_{\tilde{\bm{s}}_{i} \in \mathcal{V}_{i}^{\text{CVaR}}} \left\{\tilde{\mathbf{W}}_{i}\tilde{\bm{s}}_{i}\right\}.
\]

Given any $\tilde{\bm{w}}_{i}$, consider the minimization problem on the right-hand side
\[
    \Psi := \min_{\tilde{\bm{s}}_{i}, \bm{\eta}} \left\{ \tilde{\mathbf{W}}_{i} \tilde{\bm{s}}_{i} \left| \begin{array}{l} 
    \tilde{\bm{s}}_{i} = \sum_{i \in \mathcal{I}} \eta_{i} \bm{s}_{i} \\
    \bm{e}^{\top} \bm{\eta} = 1 \\
    \bm{\eta} \leq (\frac{1}{\alpha I}) \bm{e} \\
    \bm{\eta} \geq 0
    \end{array} \right. \right\},
\]
which is equivalent to 
\[
    \Psi := \min_{\bm{\eta}} \left\{ \tilde{\mathbf{W}}_{i} \left(\sum_{i \in \mathcal{I}} \eta_{i} \bm{s}_{i}\right) \left| \begin{array}{l} 
    \bm{e}^{\top} \bm{\eta} = 1 \\
    \bm{\eta} \leq \left(\frac{1}{\alpha I}\right) \bm{e} \\
    \bm{\eta} \geq 0
    \end{array} \right. \right\}.
\]

Let $\lambda \in \mathbb{R}$ and $\bm{\beta} \in \mathbb{R}^{I}$ be the dual variables. We have the following Lagrange function
\[
\begin{aligned}
    L(\bm{\eta}, \lambda, \bm{\beta}) & = \tilde{\mathbf{W}}_{i} \left(\sum_{i \in \mathcal{I}} \eta_{i} \bm{s}_{i}\right) + \lambda \left(1 - \bm{e}^{\top} \bm{\eta}\right) + \bm{\beta}^{\top} \left(\bm{\eta} - \left(\frac{1}{\alpha I}\right) \bm{e}\right), \\
    & =  \lambda -\left(\frac{1}{\alpha I}\right) \bm{\beta}^{\top} \bm{e} + \tilde{\mathbf{W}}_{i} \left(\sum_{i \in \mathcal{I}} \eta_{i} \bm{s}_{i}\right) - \lambda \bm{e}^{\top} \bm{\eta} + \bm{\beta}^{\top} \bm{\eta}, \\
\end{aligned}
\]
and 
\[
    \frac{\partial L(\bm{\eta}, \lambda, \bm{\beta})}{\partial \eta_{i}} = \tilde{\mathbf{W}}_{i}\bm{s}_{i} - \lambda + \beta_{i}.
\]

By Lagrangian duality, we have
\[
    \Psi = \min_{\bm{\eta} \geq 0} \max_{\lambda, \bm{\beta} \geq 0} L(\bm{\eta}, \lambda, \bm{\beta}).
\]

As is the case for linear programming, $\Psi$ is bounded below by
\[
    \Upsilon^{*} := \max_{\lambda, \bm{\beta} \geq 0} \min_{\bm{\eta} \geq 0} \lambda - \left(\frac{1}{\alpha I}\right) \bm{\beta}^{\top} \bm{e} + \tilde{\mathbf{W}}_{i} \left(\sum_{i \in \mathcal{I}} \eta_{i} \bm{s}_{i}\right) - \lambda \bm{e}^{\top} \bm{\eta} + \bm{\beta}^{\top} \bm{\eta} \leq \Psi.
\]

Let $g(\lambda, \bm{\beta}) = \min_{\bm{\eta} \geq 0} L(\bm{\eta}, \lambda, \bm{\beta})$. we have
\[
    g(\lambda, \bm{\beta}) = \left\{ \begin{array}{ll}
    \lambda - \left(\frac{1}{\alpha I}\right) \bm{\beta}^{\top} \bm{e}, & \tilde{\mathbf{W}}_{i}\bm{s}_{i} - \lambda + \beta_{i} \geq 0, \forall i \in \mathcal{I}, \\
    - \infty, & \text{otherwise},
    \end{array} \right.
\]
which can be redefined as the optimal value of the following dual problem
\[
\begin{aligned}
    \Upsilon^{*} := \max_{\lambda, \bm{\beta}} \text{ } & \lambda - \left(\frac{1}{\alpha I}\right) \bm{\beta}^{\top} \bm{e}, \\
    \mathrm{s.t.} \text{ } & \tilde{\mathbf{W}}_{i}\bm{s}_{i} - \lambda + \beta_{i} \geq 0, && \quad  \forall i \in \mathcal{I}, \\
    & \bm{\beta} \geq 0.
\end{aligned}
\]

The above dual problem is strictly feasible and the feasible set is bounded. Thus, the strong duality holds, i.e., $\Upsilon^{*} = \Psi$ and both primal and dual values are attained. By reintegrating the inequality constraint, Equation \eqref{eq-14} under Assumption \ref{assumption-08} is equivalent to the following linear programming problem
\[
    \begin{aligned}
        \max_{\tilde{\bm{w}}_{i} \in \bar{\mathcal{W}}_{i}, \tilde{z}_{i}, \lambda, \bm{\beta}} \text{ } &  \tilde{z}_{i}, \\
        \mathrm{s.t.} \text{ } & \left(\lambda - \frac{1}{\alpha I} \sum_{i^{\prime} \in \mathcal{I}} \beta_{i^{\prime}}\right)\bm{e} \geq R \tilde{\bm{U}}_{i}^{*} \tilde{z}_{i}, \\ 
        & \tilde{\mathbf{W}}_{i} \bm{s}_{i^{\prime}} + (- \lambda + \beta_{i^{\prime}})\bm{e} \geq 0,  &&\quad \forall i^{\prime} \in \mathcal{I}, \\
        & \bm{\beta} \geq 0,
    \end{aligned}
\]
which gives the result in Proposition \ref{proposition-08}. \qed

\subsection{Proof of Theorem \ref{theorem-01}}

We first briefly present the optimal solution to the second-stage problem of OPA-PR under the worst-case attribute rankings, following the proof of Theorem 2 in \citet{W24}. For any $i \in \mathcal{I}$, the second-stage problem of OPA-PR with the worst-case attribute rankings is a linear programming problem with the following constraints
\[
    \left\{ \begin{array}{l}
        R \tilde{\bm{U}}_{i}^{*} \tilde{z}_{i} - \tilde{\mathbf{W}}_{i} \tilde{\bm{s}}_{i}^{*} \leq 0, \\
        \bm{e}^{\top} \tilde{\bm{w}}_{i} = 1.
    \end{array} \right.
\]
The coefficient matrix of this system is full rank. Thus, any solution satisfying the KKT conditions is optimal. Assuming the inequality constraints are active, we have
\[
    R\tilde{U}_{ijr}^{*} \tilde{z}_{i}^{*} = \tilde{s}_{ij}^{*} \tilde{w}_{ijr}^{*}.
\]

Substituting these into the normalization constraint yields
\[
    \sum_{j \in \mathcal{J}} \sum_{r \in \mathcal{R}} \frac{R\tilde{U}_{ijr}^{*}}{\tilde{s}_{ij}^{*}} \tilde{z}_{i}^{*} = 1.
\]

By reintegrating the equality, we obtain
\[
    \tilde{z}_{i}^{*} = 1 \Big/ \left( \sum_{j \in \mathcal{J}} \sum_{r \in \mathcal{R}} \frac{R\tilde{U}_{ijr}^{*}}{\tilde{s}_{ij}^{*}} \right),
\]
and 
\[
    \tilde{w}_{ijr}^{*} = \frac{R \tilde{U}_{ijr}^{*} \tilde{z}_{i}^{*}}{\tilde{s}_{ij}^{*}}, \quad \forall (j,r) \in \mathcal{J} \times \mathcal{R}.
\]

It is straightforward to verify that $(\tilde{z}_{i}^{*}, \bm{\tilde{w}}_{i}^{*})$ satisfies the KKT conditions, confirming it as the optimal solution to the second-stage problem of OPA-PR.

We next demonstrate the decomposability of the optimal weights. Let 
\[
    z_{i}^{*} = 1 \Big/ \left( \sum_{j \in \mathcal{J}} \sum_{r \in \mathcal{R}} \frac{R}{\tilde{s}_{ij}^{*}} \right),
\]
which represents the weight disparity with normalized utilities across the ranked alternatives from the first-stage elicitation, where the sum of utilities of ranked alternatives equals 1. We have
\[
    \tilde{w}_{ijr}^{*} =\frac{R \tilde{U}_{ijr}^{*} z_{i}^{*}}{\tilde{s}_{ij}^{*}} \left(\frac{\tilde{z}_{i}^{*}}{z_{i}^{*}}\right) = \frac{1}{\tilde{s}_{ij}^{*} \sum_{j\in\mathcal{J}} \frac{1}{\tilde{s}_{ij}^{*}}} \left(\frac{\tilde{z}_{i}^{*}}{z_{i}^{*}}\tilde{U}_{ijr}^{*}\right) = w_{ij}^{WR} w_{ijr}^{WU},
\]
where 
\[
    w_{ij}^{WR} = \frac{1}{s_{ij}^{*} \sum_{j \in \mathcal{J}} \frac{1}{s_{ij}^{*}}}, \quad \forall j \in \mathcal{J},
\]
and 
\[
    w_{ijr}^{WU} = \left(\tilde{U}_{ijr}^{*} \sum_{j \in \mathcal{J}} \frac{1}{s_{ij}^{*}} \right) \Bigg/ \left(\sum_{j\in \mathcal{J}}\sum_{r \in \mathcal{R}} \frac{\tilde{U}_{ijr}^{*}}{s_{ij}^{*}} \right), \quad \forall (j, r) \in \mathcal{J} \times \mathcal{R}.
\]

Notice that $\sum_{j \in \mathcal{J}} w_{ij}^{WR} = 1$ and $\sum_{j \in \mathcal{J}} \sum_{r \in \mathcal{R}} w_{ijr}^{WU} = 1$, which gives the results in Theorem \ref{theorem-01}. \qed

\subsection{Proof of Proposition \ref{proposition-08}}

For any $(i,j) \in \mathcal{I} \times \mathcal{J}$, since $\tilde{\mathcal{U}}_{ij} \subset \mathcal{U}_{ij}$, then $\tilde{\rho}_{ij} = \underset{\tilde{u}_{ij} \in \tilde{\mathcal{U}}_{ij}}{\min} \mathbb{E}_{\mathbb{P}} [\tilde{u}_{ij}(h(\bm{x},\bm{\xi}))] \geq \underset{u_{ij} \in \mathcal{U}_{ij}}{\min} \mathbb{E}_{\mathbb{P}} [u_{ij}(h(\bm{x},\bm{\xi}))] = \rho_{ij}$.
It sufficices to show that $\tilde{\rho}_{ij} \leq \rho_{ij}$.
Let $\sigma$ be a sufficiently small positive number such that $\mathbb{E}_{\mathbb{P}} [u_{ij}^{\sigma}(h(\bm{x},\bm{\xi}))] \leq \rho_{ij} + \sigma$ with $u_{ij}^{\sigma} \in \mathcal{U}_{ij}$.
By Proposition \ref{proposition-03}, since $\psi_{l_{ij}}$ for all $l_{ij} \in \mathcal{L}_{ij}$ are step functions over $\Theta$, there exists a piecewise linear concave function $\tilde{u}_{ij}^{\sigma} \in \tilde{\mathcal{U}}_{ij}$ such that $\tilde{u}_{ij}^{\sigma} \leq u_{ij}^{\sigma}$ for all $\tau_{h_{ij}} \in \Theta$ and $h_{ij} \in \mathcal{H}_{ij}$, which implies $\mathbb{E}_{\mathbb{P}}[\tilde{u}_{ij}^{\sigma}(h(\bm{x},\bm{\xi}))] \leq \mathbb{E}_{\mathbb{P}}[u_{ij}^{\sigma}(h(\bm{x},\bm{\xi}))]$.
Therefore, we have $\tilde{\rho}_{ij} = \underset{\tilde{u} \in \tilde{\mathcal{U}}_{ij}}{\min} \mathbb{E}_{\mathbb{P}} [\tilde{u}_{ij}(h(\bm{x},\bm{\xi}))] \leq \mathbb{E}_{\mathbb{P}}[\tilde{u}_{ij}^{\sigma}(h(\bm{x},\bm{\xi}))] \leq \mathbb{E}_{\mathbb{P}}[u_{ij}^{\sigma}(h(\bm{x},\bm{\xi}))] \leq \rho_{ij} + \sigma$. Since $\sigma$ is arbitrarily small, it follows that $\tilde{\rho}_{ij} \leq \rho_{ij}$. It follows that $\tilde{\rho}_{ij} = \rho_{ij}$. \qed

\subsection{Proof of Proposition \ref{proposition-09}}

For any $(i,j) \in \mathcal{I} \times \mathcal{J}$, given the piecewise linear approximated $\tilde{u}_{ij}$ with jumps at $\tau_{h_{ij}}$ for all $h_{ij} \in \mathcal{H}_{ij}$, we have
\[
\begin{aligned}
    \int_{\Theta} \psi_{l_{ij}}(\tau) \mathrm{d}\tilde{u}_{ij}(\tau) & = \sum_{h_{ij} \in \mathcal{H}_{ij}} \int_{\tau_{h_{ij}-1}}^{\tau_{h_{ij}}} \psi_{l_{ij}}(\tau)\mathrm{d}\tilde{u}(\tau), \\
    & = \sum_{h_{ij} \in \mathcal{H}_{ij}} \frac{u_{ij}(\tau_{h_{ij}}) - u_{ij}(\tau_{h_{ij}-1})}{\tau_{h_{ij}} - \tau_{h_{ij}-1}} \int_{\tau_{h_{ij}-1}}^{\tau_{h_{ij}}} \psi_{l_{ij}}(\tau)\mathrm{d} \tau, \\
    & = \sum_{h_{ij} \in \mathcal{H}_{ij}} \left(u_{ij}(\tau_{h_{ij}}) - u_{ij}(\tau_{h_{ij}-1})\right) \frac{\int_{\tau_{h_{ij}-1}}^{\tau_{h_{ij}}} \psi_{l_{ij}}(\tau)\mathrm{d} \tau}{\tau_{h_{ij}} - \tau_{h_{ij}-1}}, \\
    & = \sum_{h_{ij} \in \mathcal{H}_{ij}} \left(u_{ij}(\tau_{h_{ij}}) - u_{ij}(\tau_{h_{ij}-1})\right) \psi_{l_{ij}}(\tau_{h_{ij}}^{\prime}),
\end{aligned}
\]
where $\tau_{h_{ij}}^{\prime} \in (\tau_{h_{ij}-1},\tau_{h_{ij}}]$ and $\psi_{l_{ij}}(\tau_{h_{ij}}^{\prime})$ is constant. 

On the other hand, consider the step-like approximation $\tilde{\psi}_{l_{ij}}$ for all $l_{ij} \in \mathcal{L}_{ij}$, which is discrete step function with jumps at $\tau_{h_{ij}}$ for all $h_{ij} \in \mathcal{H}_{ij}$. Then, we have
\[
\begin{aligned}
    \int_{\Theta} \tilde{\psi}_{l_{ij}}(\tau) \mathrm{d}u_{ij}(\tau) & = \sum_{h_{ij} \in \mathcal{H}_{ij}} \int_{\tau_{h_{ij}-1}}^{\tau_{h_{ij}}} \tilde{\psi}_{l_{ij}}(\tau) \mathrm{d}u_{ij}(\tau), \\
    & = \sum_{h_{ij} \in \mathcal{H}_{ij}} \int_{\tau_{h_{ij}-1}}^{\tau_{h_{ij}}} \psi_{l_{ij}}(\tau_{h_{ij}}^{\prime}) \mathrm{d}u_{ij}(\tau), \\
    & = \sum_{h_{ij} \in \mathcal{H}_{ij}} \psi_{l_{ij}}(\tau_{h_{ij}}^{\prime}) \left(u_{ij}(\tau_{h_{ij}}) - u_{ij}(\tau_{h_{ij}-1})\right).
\end{aligned}
\]
Thus, the step-like approximation of $\psi$ is equivalent to the piecewise linear approximation of $u$. \qed

\subsection{Proof of Lemma \ref{lemma-02}}

Given $\mathscr{F} := \left\{ f = I_{[0,\theta]}(\cdot) \right\}$, $d_{\mathscr{F}}(u_{1},u_{2})$ is well-defined for any $u_{1},u_{2} \in \mathcal{U}$.
For any $(i,j) \in \mathcal{I} \times \mathcal{J}$, since $u_{ij}^{*}$ is Lipschitz continuous with modulus $G$ and increasing on $\Theta$, we have
\[
    0 \leq \frac{u_{ij}^{*}(\tau_{h_{ij}}) - u_{ij}^{*}(\tau_{h_{ij}-1})}{\tau_{h_{ij}} - \tau_{h_{ij}-1}} \leq G, \quad \forall h_{ij} \in \mathcal{H}_{ij}.
\]
Then, by Proposition \ref{proposition-03}, for any $\tau_{h_{ij}} \in (\tau_{h_{ij}-1},\tau_{h_{ij}}]$, we have 
\[
\begin{aligned}
    \left| \tilde{u}_{ij}^{*}(\tau) - u_{ij}^{*}(\tau) \right| & = \left| u_{ij}^{*}(\tau_{h_{ij}-1}) + \frac{u_{ij}^{*}(\tau_{h_{ij}}) - u_{ij}^{*}(\tau_{h_{ij}-1})}{\tau_{h_{ij}} - \tau_{h_{ij}-1}}(\tau - \tau_{h_{ij}-1}) - u_{ij}^{*}(\tau) \right|, \\
    & = \left| \frac{\tau - \tau_{h_{ij}-1}}{\tau_{h_{ij}} - \tau_{h_{ij}-1}} \left( u_{ij}^{*}(\tau_{h_{ij}}) - u_{ij}^{*}(\tau) \right) + \frac{\tau_{h_{ij}} - \tau}{\tau_{h_{ij}} - \tau_{h_{ij}-1}} \left( u_{ij}^{*}(\tau_{h_{ij}-1}) - u_{ij}^{*}(\tau) \right) \right|, \\
    & \leq \left| \frac{\tau - \tau_{h_{ij}-1}}{\tau_{h_{ij}} - \tau_{h_{ij}-1}} \left( u_{ij}^{*}(\tau_{h_{ij}}) - u_{ij}^{*}(\tau) \right) \right| + \left| \frac{\tau_{h_{ij}} - \tau}{\tau_{h_{ij}} - \tau_{h_{ij}-1}} \left( u_{ij}^{*}(\tau_{h_{ij}-1}) - u_{ij}^{*}(\tau) \right) \right|, \\
    & \leq \frac{\tau - \tau_{h_{ij}-1}}{\tau_{h_{ij}} - \tau_{h_{ij}-1}} \left| G(\tau_{h_{ij}} - \tau_{h_{ij}-1})  \right| + \frac{\tau_{h_{ij}} - \tau}{\tau_{h_{ij}} - \tau_{h_{ij}-1}} \left| G(\tau_{h_{ij}} - \tau_{h_{ij}-1}) \right|, \\
    & = G(\tau_{h_{ij}} - \tau_{h_{ij}-1}), \\
    & \leq G \zeta_{ij},
\end{aligned}
\]
where $\zeta_{ij} = \max_{h_{ij} \in \mathcal{H}_{ij}} (\tau_{h_{ij}} - \tau_{h_{ij}-1})$, which gives Equation \eqref{eq-19}. Since $u_{ij}^{*}(0) = 0$ and $u_{ij}^{*}(R) = 1$, we have $G \geq 1/R$.
\qed

\subsection{Proof of Theorem \ref{theorem-02}}

By Proposition \ref{proposition-03}, since $u_{ij}^{*}(0) = \tilde{u}_{ij}^{*}(0) = 0$ such that $U_{ijr}^{*} = u_{ij}^{*}(R-r+1)$ and $\tilde{U}_{ijr}^{*} = \tilde{u}_{ij}^{*}(R-r+1)$. By Lemma \ref{lemma-02}, we have
\[
    \bm{U}_{i}^{*} - (G\zeta_{ij})\bm{e} \leq \tilde{\bm{U}}_{i}^{*} \leq \bm{U}_{i}^{*} + (G\zeta_{ij})\bm{e}, \quad \forall (i,j) \in \mathcal{I} \times \mathcal{J}.
\]

Consider the second-stage problem of OPA-PR with the unapproximated utilities
\[
    \max_{\bar{z}_{i}, \bar{\bm{w}}_{i}} \left\{\bar{z}_{i}: (\bar{z}_{i}, \bar{\bm{w}}_{i}) \in \mathcal{S}^{1} := \left\{ \bar{z}_{i} \in \mathbb{R}, \bar{\bm{w}}_{i} \in \mathbb{R}^{JR} \left| \begin{array}{l}
        R \bm{U}_{i}^{*} \bar{z}_{i} \leq \bar{\mathbf{W}}_{i} \tilde{\bm{s}}_{i}^{*} \\
        \bm{e}^{\top} \bar{\bm{w}}_{i} = 1 \\
        \bar{\bm{w}}_{i} \geq 0
        \end{array} \right. \right\} \right\},
\]
and with the PLA-approximated utilities
\[
    \max_{\tilde{z}_{i}, \tilde{\bm{w}}_{i}} \left\{\tilde{z}_{i}: (\tilde{z}_{i}, \tilde{\bm{w}}_{i}) \in \mathcal{S}^{2} := \left\{ \tilde{z}_{i} \in \mathbb{R}, \tilde{\bm{w}}_{i} \in \mathbb{R}^{JR} \left| \begin{array}{l}
    R (\bm{U}_{i}^{*} + G\zeta_{ij}\bm{e})\tilde{z}_{i} \leq \bar{\mathbf{W}}_{i} \tilde{\bm{s}}_{i}^{*} \\
    \bm{e}^{\top} \tilde{\bm{w}}_{i} = 1 \\
    \tilde{\bm{w}}_{i} \geq 0
    \end{array} \right. \right\} \right\}.
\]

Observe that $G\zeta_{ij}\tilde{z}_{i} \geq 0$ such that $\mathcal{S}^{2} \subseteq \mathcal{S}^{1}$, which implies that $\tilde{z}_{i}^{*} \leq \bar{z}_{i}^{*}$. Let $(\theta_{i}^{*}, \bm{\lambda}_{i}^{*})$ be the dual optimal solution of Equation \eqref{eq-20}. Assume that $(\tilde{z}_{i}, \tilde{\bm{w}}_{i})$ is feasible for the problem with PLA-approximated utilities, which is also feasible for the unapproximated problem. By strong duality, we have
\[
\begin{aligned}
    \bar{z}_{i}^{*} & \geq \tilde{z}_{i} + \theta_{i}^{*} (1 - \bm{e}^{\top} \tilde{\bm{w}}_{i}) + (\bm{\lambda}_{i}^{*})^{\top} \left(\bar{\mathbf{W}}_{i} \tilde{\bm{s}}_{i}^{*} - R \bm{U}_{i}^{*} \tilde{z}_{i}\right), \\
    & \geq \tilde{z}_{i} + RG\zeta_{ij}(\bm{\lambda}_{i}^{*})^{\top}\bm{e} \tilde{z}_{i}, \\
    & = (1 + RG\zeta_{ij}(\bm{\lambda}_{i}^{*})^{\top}\bm{e}) \tilde{z}_{i},
\end{aligned}
\]
where the second inequality follows from $\bar{\mathbf{W}}_{i} \tilde{\bm{s}}_{i}^{*} - R \bm{U}_{i}^{*} \tilde{z}_{i} \geq RG\zeta_{ij}(\bm{\lambda}_{i}^{*})^{\top}\bm{e} \tilde{z}_{i}$ for any feasible solution for the approximated problem, $\bm{\lambda}_{i}^{*} \geq 0$, and $\bm{e}^{\top} \tilde{\bm{w}}_{i} = 1$. For any feasible $\tilde{z}_{i}$, we have
\[
    \tilde{z}_{i} \leq \frac{1}{(1 + RG\zeta_{ij}(\bm{\lambda}_{i}^{*})^{\top}\bm{e})} \bar{z}_{i}^{*},
\]
which gives a upper bound for $\tilde{z}_{i}^{*}$, i.e., 
\[
    \tilde{z}_{i}^{*} \leq \frac{1}{(1 + RG\zeta_{ij}(\bm{\lambda}_{i}^{*})^{\top}\bm{e})} \bar{z}_{i}^{*}.
\]

Consider the approximated problem with the following feasible set
\[
    \mathcal{S}^{3} := \left\{ \tilde{z}_{i} \in \mathbb{R}, \tilde{\bm{w}}_{i} \in \mathbb{R}^{JR} \left| \begin{array}{l}
        R (\bm{U}_{i}^{*} - G\zeta_{ij}\bm{e})\tilde{z}_{i} \leq \bar{\mathbf{W}}_{i} \tilde{\bm{s}}_{i}^{*} \\
        \bm{e}^{\top} \tilde{\bm{w}}_{i} = 1 \\
        \tilde{\bm{w}}_{i} \geq 0
        \end{array} \right. \right\}.
\]
which yields $\mathcal{S}^{1} \subseteq \mathcal{S}^{3}$ and $\bar{z}_{i}^{*} \leq \tilde{z}_{i}^{*}$. Assume that $(\bar{z}_{i}, \bar{\bm{w}}_{i})$ is feasible for the unapproximated problem, which is also feasible for the approximated problem. Similarly, by strong duality, we have
\[
\begin{aligned}
    \tilde{z}_{i}^{*} \geq \bar{z}_{i}^{*} & \geq \bar{z}_{i} + \theta_{i}^{*} (1 - \bm{e}^{\top} \bar{\bm{w}}_{i}) + (\bm{\lambda}_{i}^{*})^{\top} \left(\bar{\mathbf{W}}_{i} \tilde{\bm{s}}_{i}^{*} - R \bm{U}_{i}^{*} \bar{z}_{i} \right), \\
    & \geq \bar{z}_{i} - RG\zeta_{ij}(\bm{\lambda}_{i}^{*})^{\top}\bm{e} \bar{z}_{i}, \\
    & = (1 - RG\zeta_{ij}(\bm{\lambda}_{i}^{*})^{\top}\bm{e}) \bar{z}_{i},
\end{aligned}
\]
which gives a lower bound for $\tilde{z}_{i}^{*}$, i.e, 
\[
    (1 - RG\zeta_{ij}(\bm{\lambda}_{i}^{*})^{\top}\bm{e}) \bar{z}_{i}^{*} \leq \tilde{z}_{i}^{*}.
\]

By the proof of Theorem \ref{theorem-01}, for all $(j,r) \in \mathcal{J} \times \mathcal{R}$, we have 
\[
    \tilde{w}_{ijr}^{*} = \frac{R \tilde{U}_{ijr}^{*} \tilde{z}_{i}^{*}}{\tilde{s}_{ij}^{*}} \Leftrightarrow \frac{(1 - RG\zeta_{ij}(\bm{\lambda}_{i}^{*})^{\top}\bm{e})R \tilde{U}_{ijr}^{*}}{\tilde{s}_{ij}^{*}}\bar{z}_{i}^{*} \leq \tilde{w}_{ijr}^{*} \leq \frac{R \tilde{U}_{ijr}^{*}}{(1 + RG\zeta_{ij}(\bm{\lambda}_{i}^{*})^{\top}\bm{e})\tilde{s}_{ij}^{*}}\bar{z}_{i}^{*},
\]
which gives the results in Theorem \ref{theorem-02}. \qed

\end{document}